\documentclass[12pt]{amsart}
\usepackage{amscd,amsmath,amssymb,amsfonts}
\usepackage[cmtip, all]{xy}
\begin{document}
\newtheorem{thm}{Theorem}[subsection]
\newtheorem{prop}[thm]{Proposition}
\newtheorem{lem}[thm]{Lemma}
\newtheorem{cor}[thm]{Corollary}
\newtheorem{conj}[thm]{Conjecture}
\theoremstyle{definition}
\newtheorem{defin}[thm]{Definition}
\theoremstyle{remark}
\newtheorem{remk}[thm]{Remark}
\newtheorem{remks}[thm]{Remarks}
\newtheorem{exm}[thm]{Example}
\newtheorem{exms}[thm]{Examples}
\newtheorem{notat}[thm]{Notation}
\numberwithin{equation}{subsection}

\newenvironment{enum}{\makeatletter
\renewcommand{\theenumi}{\roman{enumi}}
\renewcommand{\labelenumi}{(\theenumi)}
\makeatother
\begin{enumerate}}%
{\end{enumerate}}

\newenvironment{enumera}{\makeatletter
\renewcommand{\theenumi}{\alph{enumi}}
\renewcommand{\labelenumi}{(\theenumi)}
\makeatother
\begin{enumerate}}%
{\end{enumerate}}

\newenvironment{Def}{\begin{defin}}%
{\hfill$\square$\end{defin}}

\newenvironment{rem}{\begin{remk}}%
{\hfill$\square$\end{remk}}

\newenvironment{rems}{\begin{remks}}%
{\hfill$\square$\end{remks}}

\newenvironment{ex}{\begin{exm}}%
{\hfill$\square$\end{exm}}

\newenvironment{exs}{\begin{exms}}%
{\hfill$\square$\end{exms}}

\newenvironment{note}{\begin{notat}}%
{\hfill$\square$\end{notat}}

\newcommand{\thmref}{Theorem~\ref}
\newcommand{\propref}{Proposition~\ref}
\newcommand{\corref}{Corollary~\ref}
\newcommand{\defref}{Definition~\ref}
\newcommand{\lemref}{Lemma~\ref}
\newcommand{\exref}{Example~\ref}
\newcommand{\remref}{Remark~\ref}
\newcommand{\chref}{Chapter~\ref}
\newcommand{\apref}{Appendix~\ref}
\newcommand{\secref}{Section~\ref}
\newcommand{\conjref}{Conjecture~\ref}

\newcommand{\sA}{{\mathcal A}}
\newcommand{\sB}{{\mathcal B}}
\newcommand{\sC}{{\mathcal C}}
\newcommand{\sD}{{\mathcal D}}
\newcommand{\sE}{{\mathcal E}}
\newcommand{\sF}{{\mathcal F}}
\newcommand{\sG}{{\mathcal G}}
\newcommand{\sH}{{\mathcal H}}
\newcommand{\sI}{{\mathcal I}}
\newcommand{\sJ}{{\mathcal J}}
\newcommand{\sK}{{\mathcal K}}
\newcommand{\sL}{{\mathcal L}}
\newcommand{\sM}{{\mathcal M}}
\newcommand{\sN}{{\mathcal N}}
\newcommand{\sO}{{\mathcal O}}
\newcommand{\sP}{{\mathcal P}}
\newcommand{\sQ}{{\mathcal Q}}
\newcommand{\sR}{{\mathcal R}}
\newcommand{\sS}{{\mathcal S}}
\newcommand{\sT}{{\mathcal T}}
\newcommand{\sU}{{\mathcal U}}
\newcommand{\sV}{{\mathcal V}}
\newcommand{\sW}{{\mathcal W}}
\newcommand{\sX}{{\mathcal X}}
\newcommand{\sY}{{\mathcal Y}}
\newcommand{\sZ}{{\mathcal Z}}
\newcommand{\A}{{\mathbb A}}
\newcommand{\B}{{\mathbb B}}
\newcommand{\C}{{\mathbb C}}
\newcommand{\D}{{\mathbb D}}
\newcommand{\E}{{\mathbb E}}
\newcommand{\F}{{\mathbb F}}
\newcommand{\G}{{\mathbb G}}
\renewcommand{\H}{{\mathbb H}}
\newcommand{\I}{{\mathbb I}}
\newcommand{\J}{{\mathbb J}}
\newcommand{\M}{{\mathbb M}}
\newcommand{\N}{{\mathbb N}}
\renewcommand{\P}{{\mathbb P}}
\newcommand{\Q}{{\mathbb Q}}
\newcommand{\R}{{\mathbb R}}
\newcommand{\T}{{\mathbb T}}
\newcommand{\U}{{\mathbb U}}
\newcommand{\V}{{\mathbb V}}
\newcommand{\W}{{\mathbb W}}
\newcommand{\X}{{\mathbb X}}
\newcommand{\Y}{{\mathbb Y}}
\newcommand{\Z}{{\mathbb Z}}
\renewcommand{\1}{{\mathbb{S}}}

\newcommand{\fm}{{\mathfrak m}}
\newcommand{\an}{{\rm an}}
\newcommand{\alg}{{\rm alg}}
\newcommand{\cl}{{\rm cl}}
\newcommand{\Alb}{{\rm Alb}}
\newcommand{\CH}{{\rm CH}}
\newcommand{\mc}{\mathcal}
\newcommand{\mb}{\mathbb}
\newcommand{\surj}{\twoheadrightarrow}
\newcommand{\red}{{\rm red}}
\newcommand{\codim}{{\rm codim}}
\newcommand{\rank}{{\rm rank}}
\newcommand{\Pic}{{\rm Pic}}
\newcommand{\Div}{{\rm Div}}
\newcommand{\Hom}{{\rm Hom}}
\newcommand{\im}{{\rm im}}
\newcommand{\Spec}{{\rm Spec \,}}
\newcommand{\sing}{{\rm sing}}
\newcommand{\Char}{{\rm char}}
\newcommand{\Tr}{{\rm Tr}}
\newcommand{\Gal}{{\rm Gal}}
\newcommand{\Min}{{\rm Min \ }}
\newcommand{\Max}{{\rm Max \ }}
\newcommand{\supp}{{\rm supp}\,}
\newcommand{\0}{\emptyset}
\newcommand{\sHom}{{\mathcal{H}{om}}}

\newcommand{\Nm}{{\operatorname{Nm}}}
\newcommand{\NS}{{\operatorname{NS}}}
\newcommand{\id}{{\operatorname{id}}}
\newcommand{\Zar}{{\text{\rm Zar}}} 
\newcommand{\Ord}{{\mathbf{Ord}}}
\newcommand{\FSimp}{{{\mathcal FS}}}
\newcommand{\inj}{{\text{\rm inj}}}  
\newcommand{\Sch}{{\operatorname{\mathbf{Sch}}}} 
\newcommand{\cosk}{{\operatorname{\rm cosk}}} 
\newcommand{\sk}{{\operatorname{\rm sk}}} 
\newcommand{\subv}{{\operatorname{\rm sub}}}
\newcommand{\bary}{{\operatorname{\rm bary}}}
\newcommand{\Comp}{{\mathbf{SC}}}
\newcommand{\IComp}{{\mathbf{sSC}}}
\newcommand{\Top}{{\mathbf{Top}}}
\newcommand{\holim}{\mathop{{\rm holim}}}
\newcommand{\Holim}{\mathop{{\it holim}}}
\newcommand{\op}{{\text{\rm op}}}
\newcommand{\<}{\mathopen<}
\renewcommand{\>}{\mathclose>}
\newcommand{\Sets}{{\mathbf{Sets}}}
\newcommand{\del}{\partial}
\newcommand{\fib}{{\operatorname{\rm fib}}}
\renewcommand{\max}{{\operatorname{\rm max}}}
\newcommand{\bad}{{\operatorname{\rm bad}}}
\newcommand{\Spt}{{\mathbf{Spt}}}
\newcommand{\Spc}{{\mathbf{Spc}}}
\newcommand{\Sm}{{\mathbf{Sm}}}
\newcommand{\cofib}{{\operatorname{\rm cofib}}}
\newcommand{\hocolim}{\mathop{{\rm hocolim}}}
\newcommand{\Glu}{{\mathbf{Glu}}}
\newcommand{\can}{{\operatorname{\rm can}}}
\newcommand{\Ho}{{\mathbf{Ho}}}
\newcommand{\GL}{{\operatorname{\rm GL}}}
\newcommand{\sq}{\square}
 \newcommand{\Ab}{{\mathbf{Ab}}}
\newcommand{\Tot}{{\operatorname{\rm Tot}}}
\newcommand{\loc}{{\operatorname{\rm s.l.}}}
\newcommand{\HZ}{{\operatorname{\sH \Z}}}
\newcommand{\Cyc}{{\operatorname{\rm Cyc}}}
\newcommand{\cyc}{{\operatorname{\rm cyc}}}
 \newcommand{\RCyc}{{\operatorname{\widetilde{\rm Cyc}}}}
 \newcommand{\Rcyc}{{\operatorname{\widetilde{\rm cyc}}}}
\newcommand{\Sym}{{\operatorname{\rm Sym}}}
\newcommand{\fin}{{\operatorname{\rm fin}}}
\newcommand{\SH}{{\operatorname{\sS\sH}}}
\newcommand{\Anna}{{$\mathcal{ANNA\ SARAH\ LEVINE}$}}
\newcommand{\Wedge}{{\Lambda}}
\newcommand{\eff}{{\operatorname{\rm eff}}}
\newcommand{\rcyc}{{\operatorname{\rm rev}}}
\newcommand{\DM}{{\operatorname{\mathcal{DM}}}}
\newcommand{\GW}{{\operatorname{\rm{GW}}}}
\newcommand{\sSets}{{\mathbf{Spc}}}
\newcommand{\Nis}{{\operatorname{Nis}}}
\newcommand{\et}{{\text{\'et}}}
\newcommand{\Cat}{{\mathbf{Cat}}}
\newcommand{\ds}{{/\kern-3pt/}}
\newcommand{\bD}{{\mathbf{D}}}
\newcommand{\res}{{\operatorname{res}}}

\newcommand{\qf}{{\operatorname{q.fin.}}}

\newcommand{\fil}{\phi}
\newcommand{\barfil}{\sigma}

\newcommand{\Fac}{{\mathop{\rm{Fac}}}}
\newcommand{\Fun}{{\mathbf{Func}}}

\newcommand{\cp}{\coprod}
\newcommand{\ess}{\text{\rm{ess}}}
\newcommand{\dimrel}{\operatorname{\text{dim-rel}}}

\def\un#1{\underline{#1}}
\def\arr#1{\xrightarrow{#1}}
\def\arrn#1{\xrightarrow[#1]{}}
\def\arl#1{\xleftarrow{#1}}

\newcommand{\gen}{{\text{gen}}}
\newcommand{\Gr}{{\text{\rm Gr}}}
\newcommand{\Ind}{{\operatorname{Ind}}}
\newcommand{\Supp}{{\operatorname{Supp}}}
\newcommand{\Cone}{{\operatorname{Cone}}}
\newcommand{\Tor}{{\operatorname{Tor}}}
\newcommand{\cok}{{\operatorname{coker}}}
\newcommand{\Br}{{\operatorname{Br}}}
\newcommand{\sm}{{\operatorname{sm}}}
\newcommand{\Proj}{{\operatorname{Proj}}}
\newcommand{\Pos}{{\operatorname{\mathbf{Pos}}}}
\newcommand{\colim}{\mathop{\text{colim}}}
\newcommand{\ps}{{\operatorname{\text{p.s.s.}}}}

\addtocounter{section}{-1}
\title{Chow's moving lemma and the homotopy coniveau tower}

\author{Marc Levine}

\address{
Department of Mathematics\\
Northeastern University\\
Boston, MA 02115\\
USA}

\email{marc@neu.edu}


\thanks{The author gratefully acknowledges the support of the Humboldt Foundation through the Wolfgang Paul Program, and support of the NSF via grant DMS 0140445}

\renewcommand{\abstractname}{Abstract}
\begin{abstract}  We show how the classical moving lemma of Chow extends to give functoriality for the specta arising in the homotopy coniveau tower.\\
\end{abstract}

\maketitle
\tableofcontents
\section{Introduction} This paper is part of a program to understand and generalize  the papers \cite{BL, FriedSus}, which construct the spectral sequence from motivic cohomology to $K$-theory. We have started a study of this approach in our paper \cite{Loc}, where, relying on some of the results of Friedlander-Suslin, we gave a method for globalizing the Bloch-Lichtenbaum spectral sequence to schemes of finite type over a Dedekind domain. This method is a bit different from the one used by Friedlander-Suslin, and is somewhat better adapted to localization properties, while being rather less convenient regarding questions of functoriality.

   It turns out that the approach of \cite{Loc} can be applied to a wide variety of cohomology theories on algebraic varieties, including all the cohomology theories arising from objects in the Morel-Voevodsky $\A^1$-stable homotopy category. Just as for $K$-theory, localization properties follow rather easily from the existing technology. The purpose of this paper is to handle the problem of functoriality in a fairly general setting.

The functoriality obtained here is used in our paper \cite{HomotopyConiveau} to extend the main results of  \cite{BL, FriedSus} to a  general setting and, at the same time, give a new proof of the main results of  \cite{BL}. This yields, we think, a much more conceptual approach to the Bloch-Lichtenbaum-Friedlander-Suslin spectral sequence.

In somewhat more detail, one can consider a functor $E:\Sm/B^\op\to\Spt$ from 
smooth $B$-schemes to spectra (for example $X\mapsto K(X)$, the $K$-theory spectrum of a $B$-scheme $X$). For such a functor, and for an $X$ in $\Sm/B$, we construct the {\em homotopy coniveau tower}
\begin{equation}\label{HomConTowerIntro}
\ldots\to E^{(p+1)}(X,-)\to  E^{(p)}(X,-)\to\ldots\to E^{(0)}(X,-)
\end{equation}
where the $E^{(p)}(X,-)$ are simplicial spectra with $n$-simplices the limit of the spectra with support $E^W(X\times\Delta^n)$, where $W$ is a closed subset of codimension $\ge p$ in ``good position" with respect to the faces of $\Delta^n$. This is just the evident extension of the tower used by Friedlander-Suslin in \cite{FriedSus}. One can consider this tower as the algebraic analog of the one in topology formed by applying a cohomology theory to the skeletal filtration of a CW complex.

Our main purpose in this paper is construct a functorial model of the presheaf  $E^{(p)}(X_\Nis,-)$, i.e., the presheaf of spectra $U\mapsto E^{(p)}(U,-)$ on $X_\Nis$. Our construction is based on first a  ``local functoriality", which is essentially a generalization of the  classical method used (e.g. in \cite{Roberts}) to prove Chow's moving lemma for cycles modulo rational equivalence. This local functoriality yields a weak version of  functoriality to the homotopy category of presheaves on $X$, for varying $X$, which we promote to give a full functoriality in the homotopy category of presheaves on  $\Sm/B$ (with respect to Nisnevic-local weak equivalences) by a formal rectification process. Naturally, we need to impose some addtional (but quite natural) hypotheses on our functor $E$ to achieve this.

One consequence of the general theory discussed in \cite{HomotopyConiveau} is, as mentioned above, a localization property for the spectra $E^{(p)}(X,-)$; this yields a Mayer-Vietoris property for the presheaf $E^{(p)}(X_\Nis,-)$. In particular (by the well-known theorem of Thomason \cite{Jardine1, Thomason}), the functoriality for  the presheaf $E^{(p)}(X_\Nis,-)$  up   to local weak equivalence yields functoriality for the spectra $E^{(p)}(X,-)$ up to pointwise weak equivalence.

The homotopy coniveau tower \eqref{HomConTowerIntro} leads directly to the spectral sequence
\begin{equation}\label{eqn:HomConSS0}
E_1^{p,q}:=\pi_{-p-q}(E^{(p/p+1)}(X,-))\Longrightarrow \pi_{-p-q}(E^{(0)}(X,-)),
\end{equation}
which, in the case of the $K$-theory spectrum was the main object of interest in the papers  \cite{BL, FriedSus}; here $E^{(p/p+1)}(X,-)$ is the homotopy cofiber of the map $E^{(p+1)}(X,-)\to  E^{(p)}(X,-)$.
 One major goal of this paper is to make this spectral sequence functorial on $\Sm/B$. In fact, our results yield directly the functoriality of an associated spectral sequence, where one replaces the spectra 
$E^{(p)}(X,-)$ with fibrant models for $E^{(p)}(X_\Nis,-)$.  The Mayer-Vietoris property for  $E^{(p)}(X_\Nis,-)$ mentioned in the preceding paragraph allows us to recover the functoriality of 
\eqref{eqn:HomConSS0} from this. See \S\ref{sec:HomConSS} for further details.

This paper also has roots in our preprint \cite{KThyMotCoh}. In another paper, we will complete our detailed investigation of the spectral sequence from $K$-theory to motivic cohomology, relying on the results of this paper for the functoriality. To keep some contact with our original problem, we will occasionally illustrate our results by refering to the example $E=K$.

We have had a great deal of help in developing the techniques that went in to this paper. Fabien Morel played a crucial role in numerous discussions on the $\A^1$-stable homotopy category and related topics. Conversations with Bruno Kahn, Jens Hornbostel and Marco Schlichting were very helpful, as were the lectures of Bjorn Dundas and Vladamir Voevodsky at the Sophus Lie Summer Workshop in $\A^1$-stable homotopy theory. I would also like to thank Paul Arne Ostvar for his detailed comments on an earlier version of this manuscript. I also want to thank Eckart Viehweg for his comments on section 10, which greatly simplified my earlier arguments. Finally, I am very grateful to the Humboldt Foundation and the Universit\"at Essen for their support of this reseach, especially my colleagues H\'el\`ene Esnault and Eckart Viehweg.
  
\hbox{}
\newpage

\section{Preliminaries}
\subsection{Triangulations and homotopies}\label{Triang} We fix some notation and recall the standard construction of simplicial homotopies of maps of simplicial spaces.

We have the ordered sets $[p]:=\{0<\ldots,<p\}$. $\Ord$ is the category with objects $[p]$, $p=0,1,\ldots$, and with maps the order-preserving maps of sets. For a category $\sC$, the category of  functors $F:\Ord\to \sC$ is  the category of  {\em cosimplicial} objects of $\sC$, and that of functors $F:\Ord^\op\to\sC$ the category of  {\em simplicial} objects  of $\sC$. Let $\Ord_{\le N}$ denote the full subcategory of the category $\Ord$, with objects  
$[p]$, $0\le p\le N$. We have the category of $N$-truncated cosimplicial objects of $\sC$, namely, the functors $F;\Ord_{\le N}\to \sC$, and similarly the category of $N$-truncated simplicial objects of $\sC$.

The category $\sSets$ of simplicial sets is called the category of {\em spaces}; if $F:\Ord^\op\to\sSets$ is a functor, we have the geometric realization $|F|\in\sSets$. This is defined as the quotient simplicial set 
\[
|F|:=\coprod_n F(n)\times\Delta[n]/\sim
\]
where $\Delta[n]:= m\mapsto \Hom_\Ord([m],[n])$, and $\sim$ is the standard relation
\[
(x,g_*(y))\sim (F(g)(x),y)
\]
for $g:[n]\to[n']$ in $\Ord$, $x\in F([n'])$.   We have as well the category of pointed spaces $\sSets*$ and geometric realization functors from pointed simplicial spaces  to $\sSets*$.

We give the product
$[p]\times[q]$ the product partial order
\[
(a,b)\le(a',b')\Leftrightarrow a\le a'\text{ and }b\le b'.
\]
In what follows, we work in the the category $\Pos$ of partially ordered sets, i.e. maps are order-preserving maps of partially ordered sets.
 
Let $F:\Ord^\op \to \sSets$, $G:\Ord^\op \to \sSets$ be functors, giving the geometric realizations $|F|$ and $|G|$.
Suppose we have, for each $h:[p]\to[q]\times[1]$,   a morphism
\[
H(h):G([q])\to F([p])
\]
such that, for   $g:[q]\to [r]$, $f:[s]\to[p]$, we have 
\[
H(h\circ f)=F(f)\circ H(h),\ H((g\times\id)\circ h)=H(h)\circ G(g).
\]

Restricting the morphisms  $h$ to those of the form $i_\epsilon^p:[p]\to [p]\times[1]$,
\[
i_\epsilon^p(j)=(j,\epsilon),\ \epsilon=0,1,
\]
the maps $H(i^p_\epsilon):G([p])\to F([p])$ define the maps of  simplicial spaces
\[
H_0, H_1:G\to F.
\]
Also, the maps $H(h)$ define the map in $\sSets$
\[
H:|G|\times I\to |F|,
\]
with 
\[
H_{||G|\times0}= |H_0|,\ H_{||G|\times1}= |H_1|.
\]
Indeed, $|G|\times I$ is the geometric realization of the  simplicial space 
\[
(G\times I)([p])=\coprod_{[q]}G([q])\times\Hom_\Pos([p],[q]\times[1])/\sim,
\]
where $\sim$ is the equivalence relation $(G(g)(x),h)\sim(x,(g\times\id)\circ h)$, with the evident map on
morphisms; to see this, one first considers the case $G=\Delta[n]$, from which the general case follows using the canonical isomorphism
\[
G\cong\colim_{[n], \sigma:\Delta[n]\to G}\Delta[n].
\]
Our assertion follows directly from this.

\subsection{Presheaves of simplicial sets}\label{sec:PreshSimp} For a category $\sC$, we have the category $\ps\sC$ of functors $F:\sC^\op\to \sSets$ (presheaves of simplicial sets), and $\ps_*\sC$ of functors $F:\sC^\op\to \sSets*$ (presheaves of pointed simplicial sets). 

We give $\sSets$ and $\sSets*$ the standard model structures: cofibrations are (pointed) monomorphisms, weak equivalences are weak equivalences on the geometric realization, and fibrations are detemined by the RLP with respect to trivial cofibrations; the fibrations are then exactly the Kan fibrations. We let  $|A|$ denote the geometric realization, and $[A,B]$  the homotopy classes of (pointed) maps $|A|\to |B|$.

We give  $s\sC$  and  $s_*\sC$ the model structure of functor categories described by Bousfield-Kan \cite{BousfieldKan}. That is, the cofibrations and weak equivalences are the pointwise ones, and the fibrations are determined by the RLP with respect to trivial cofibrations. We let $\sH s\sC$ and $\sH s_*\sC$ denote the associated homotopy cateogies (see \cite{Hovey} for details).

\subsection{Presheaves of spectra} Let $\Spt$ denote the category of
spectra. To fix ideas, a spectrum will be a sequence of pointed simplicial sets
$E_0, E_1,\ldots$ together with maps of pointed simplicial sets $\epsilon_n:S^1\wedge E_n\to E_{n+1}$.  Maps of spectra are maps of the underlying simplicial sets which are compatible with the attaching maps $\epsilon_n$. The stable homotopy groups $\pi_n^s(E)$ are defined by
\[
\pi_n^s(E):=\lim_{m\to\infty} [S^{m+n},E_m].
\]

The category $\Spt$ has the following model structure: Cofibrations are maps $f:E\to F$ such that $E_0\to F_0$ is a cofibration, and for each $n\ge0$, the map
\[
E_{n+1}\coprod_{S^1\wedge E_n}S^1\wedge F_n\to F_{n+1}
\]
is a cofibration.  Weak equivalences are the stable weak equivalences, i.e., maps $f:E\to F$ which induce an isomorphism on $\pi_n^s$ for all $n$. Fibrations are characterized by having the RLP with respect to trivial cofibrations.

Let $\sC$ be a category. We say that a natural transformation $f:E\to E'$ of functors  $\sC^\op\to\Spt$ is a weak equivalence if $f(X):E(X)\to E'(X)$ is a stable weak equivalence for all $X$.  

We will assume that the category $\sC$ has an initial object $\0$ and admits finite coproducts over $\0$, denoted $X\amalg Y$. A functor $E:\sC^\op\to \Spt$ is called {\em additive} if for each $X, Y$ in $\sC$, the canonical map
\[
E(X\amalg Y)\to E(X)\oplus E(Y)
\]
is a weak equivalence. An additive functor  $E:\sC^\op\to \Spt$ is called a {\em presheaf of spectra} on $\sC$. This forms a full subcategory of the functor category.

We use the following model structure on the category of presheaves of spectra (see \cite{Jardine1}): Cofibrations and weak equivalences are given pointwise, and fibrations are characterized by having the RLP with respect to trivial cofibrations. We denote this model category by   $\Spt(\sC)$, and the associated homotopy category by $\sH\Spt(\sC)$.   We write $\SH$ for the homotopy category of $\Spt$.

For the remainder of the paper, we let $B$ denote  a noetherian separated scheme of finite Krull dimension. We let $\Sm/B$ denote the category of smooth $B$-schemes of finite type over $B$. We often write $\Spt(B)$ and $\sH\Spt(B)$ for $\Spt(\Sm/B)$ and $\sH\Spt(\Sm/B)$.

For $Y\in\Sm/B$, a subscheme $U\subset Y$ of the form $Y\setminus \cup_\alpha F_\alpha$, with $\{F_\alpha\}$ a possibly infinite set of closed subsets of $Y$, is called {\em essentially smooth over } $B$; the category of essentially smooth $B$-schemes is denoted $\Sm^\ess$.

\subsection{Local model structure} 

If the category $\sC$ has a topology, there is often another model structure on $\Spt(\sC)$ which takes this into account. We consider the case of the small Nisnevic site $X_\Nis$ on a scheme $X$ (assumed to be noetherian, separated and of finite Krull dimension), and on the big Nisnevic sites $\Sm/B_\Nis$ or $\Sch/B_\Nis$, as well as the Zariski versions $X_\Zar$, etc. We describe the Nisnevic version below; the definitions and results for the Zariski topology are exactly parallel.

\begin{Def} A map $f:E\to F$ of presheaves of   spectra on $X_\Nis$ is a {\em local weak equivalence} if the induced map on the Nisnevic sheaf of stable homotopy groups $f_*:\pi_m^s(E)_\Nis\to \pi_m^s(F)_\Nis$ is an isomorphism of sheaves for all $m$. A map
$f:E\to F$ of presheaves of   spectra on $\Sm/B_\Nis$ or $\Sch/B_\Nis$ is a local weak equivalence if the restriction of $f$ to $X_\Nis$ is a local weak equivalence for all $X\in \Sm/B$ or $X\in\Sch/B$.
\end{Def}

The (Nisnevic) local model structure on the category of presheaves of spectra on $X_\Nis$ has cofibrations given pointwise, weak equivalences the local weak equivalences and  fibrations are characterized by having the RLP with respect to trivial cofibrations. We denote this model structure $\Spt_\Nis(X)$ and the associated homotopy category $\sH\Spt_\Nis(X)$. The same definitions yield the Nisnevic local model structures $\Spt_\Nis(\Sm/B)$ and $\Spt_\Nis(\Sch/B)$, with respective homotopy categories $\sH\Spt_\Nis(\Sm/B)$, $\sH\Spt_\Nis(\Sch/B)$. For details, we refer the reader to \cite{Jardine1}.

\subsection{Simplicial spectra}

Let $E:\Ord^\op\to \Spt$ be a simplicial spectrum. Writing $E=(E_0, E_1,\ldots)$ with structure maps $\epsilon_n:S^1\wedge E_n\to E_{n+1}$, we may take the geometric realizations $|E_n|$ and form the {\em total spectrum} of $E$ as the spectrum $(|E_0|, |E_1|,\ldots)$ with structure maps $|\epsilon_n|:S^1\wedge|E_n|=|S^1\wedge E_n|\to |E_{n+1}|$. 

We  use the geometric realization functor to define weak equivalences and cofibrations of  simplicial spectra, giving us the model category of simplicial spectra, as well as the local model structure for the Nisnevic or Zariski topology.

When the context makes the meaning clear, we will often omit the separate notation for the total spectrum, and freely pass between a simplicial spectrum and its associated total spectrum.

\section{The homotopy coniveau tower}
\subsection{Spectra with support} In order to axiomatize the situation, we introduce some categories which describe ``spectra with support in codimension $\ge q$". 

\begin{Def} Let $B$ be a scheme and fix an integer $q\ge0$. Let $\Sm/B^{(q)}$ be the category with objects pairs $(Y,W)$, where $Y$ is in $\Sm/B$ and $W$ is a closed subset of $Y$ with $\codim_YW\ge q$. A morphism $f:(Y,W)\to (Y',W')$ is a $B$-morphism $f:Y\to Y'$ such that
 $W\supset f^{-1}(W')$.
\end{Def}

We have the model category of presheaves $\Spt(\Sm/B^{(q)})$ and the homotopy category $\sH\Spt(\Sm/B^{(q)})$, denoted $\Spt(B^{(q)})$ and $\sH\Spt(B^{(q)})$. For $E\in \Spt(B^{(q)})$, we often write $E^W(Y)$ for $E(Y,W)$.

Take $E\in \Spt(B^{(q)})$, $(Y,W)\in \Sm/B^{(q)}$ and $W'\subset W$ a closed subset. Then $(Y,W')$ is also in $\Sm/B^{(q)}$ and the identity on $Y$ gives the morphism $\id_Y:(Y,W)\to (Y,W')$. We thus have the  map of spectra
\[
i_{W,W'}:=\id_Y^*:E^{W'}(X)\to E^W(X)
\]
with $i_{W,W'}\circ i_{W',W''}=i_{W,W''}$ for $W''\subset W'\subset W$.

Take integers $q'\ge q$.   Sending $(Y,W)\in \Sm/B^{(q')}$ to  $(Y,W)\in  \Sm/B^{(q)}$ identifies $\Sm/B^{(q')}$ with a full subcategory of  $\Sm/B^{(q)}$; we denote the inclusion functor by
\[
\rho_{q,q'}: \Sm/B^{(q')}\to  \Sm/B^{(q)}.
\]

\subsection{Examples} Our primary source of spectra with supports is a presheaf of spectra  on $\Sm/B$. Given $E\in \Spt(B)$, and $(Y,W)\in \Sm/B^{(q)}$, we set
\[
E^W(Y):=\fib(E(Y)\xrightarrow{j^*}E(Y\setminus W)),
\]
where ``$\fib$" is the homotopy fiber and $j:Y\setminus W\to Y$ is the inclusion. This clearly defines functors $\Spt(B)\to \Spt(B^{(q)})$, compatible in $q$.  As a particular example, we can take $E$ to be the $K$-theory spectrum $X\mapsto K(X)$.

Another source of examples comes from algebraic cycles. For $(Y,W)$ in $\Sm/B^{(q)}$, let $z^q_W(Y)$ denote the group of codimension $q$ algebraic cycles on $Y$ with support contained in $W$. Sending $(Y,W)$ to $z^q_W(Y)$ defines a presheaf of abelian groups  on $\Sm/B^{(q)}$; taking the associated Eilenberg-Maclane spectrum gives us the presheaf $z^q\in \Spt(B^{(q)})$.  
 
\begin{rem}
It may well be that the ``spectrum with supports" functor $\Spt(B)\to \Spt(B^{(q)})$ defined above is in fact essentially surjective, so that the introduction of the category $\Spt(B^{(q)})$ is superfluous, but we have not investigated this issue.
\end{rem}

\subsection{Support conditions}   Let $X:\Ord\to\Sm/B$ be a smooth cosimplicial $B$-scheme. For $E:\Sm/B^\op\to \Spt$ in $\Spt(B)$, we may evaluate $E$ on $X$, forming the simplicial spectrum $E(X):\Ord^\op\to \Spt$.  This is our basic method for constructing simplicial spectra; we will also need an extension of this construction by introducing {\em support conditions}.

Let $X$ be a smooth cosimplicial $B$-scheme. We have the simplicial set $\sS_{X}:\Ord^\op\to\Sets$ with $\sS_{X}(n)$ the set of closed subsets of $X(n)$. 

\begin{Def} (1) A {\em  support condition} on $X$ is a simplicial subset $\sS$ of $\sS_{X}$ such that, for each $n$, 
\begin{enumera}
\item $\sS(n)$ is closed under finite union
\item If $W$ is in $\sS(n)$ and $W'\subset W$ is a closed subset, then $W'$ is in $\sS(n)$. 
\end{enumera} 
(2) An assignment $n\mapsto W_n$, with $W_n\in\sS_{X}(n)$ is a {\em simplicial closed subset} of $X$ if for each $g:[m]\to [n]$ in $\Ord$, we have $X(g)^{-1}(W_n)\subset W_m$.  \\
(3) For a  support condition $\sS$ on a smooth cosimplicial $B$-scheme $X$, and a simplicial closed subset $W$ of $X$, we write $W\in \sS$ if $W_n$ is in $\sS(n)$ for all $n$.\\
(4) If $\codim_{X(n)}W\ge q$ for all $W\in \sS(n)$, we say that $\sS$ is {\em supported in codimension $\ge q$}. 
  \end{Def}

\begin{rems}\label{rems:support} Fix a smooth cosimplicial $B$-scheme $X$. \\
 (1) A collection of closed subsets $W_n\subset X(n)$ forms a simplicial closed subset of $X$ exactly when the complements $U(n):=X(n)\setminus W_n$ form an open cosimplicial subscheme of $X$.\\
(2) Suppose we have closed subsets $T_1\subset X(n_1)$, $\ldots$, $T_s\subset X(n_s)$. Then there is a unique minimal (in the sense of containment of closed subsets) simplicial closed subset $W$ of $X$ with $T_j\subset W_{n_j}$ for all $j=1,\ldots, s$; this is clear since the intersection of simplicial closed subsets is a simplicial closed subset. We call $W$ the simplicial closed subset generated by the $T_j$. The explicit formula for $W_m$ is
\[
W_m=\cup_{i=1}^s\cup_{g:[m]\to[n_i]}X(g)^{-1}(T_i);
\]
note that the union is finite, so $W_m$ is indeed a closed subset of $X(m)$.\\
(3) Let $n\mapsto W_n$ be a simplicial closed subset of $X$, supported in codimension $\ge q$. Then $n\mapsto (X(n),W_n)$ defines a cosimplicial object of $\Sm/B^{(q)}$.
 \end{rems}

\begin{lem}\label{lem:SimpSupp} Let $\sS$ be a  support condition on a smooth cosimplicial $B$-scheme $X$, and let $T_i$ be in $\sS(n_i)$, $i=1,\ldots, s$. Then the simplicial closed subset $W$ generated by the $T_i$ satisfies  
\begin{enumerate}
\item $T_i\subset W_{n_i}$ for $i=1,\ldots, s$.
\item $W_m\in\sS(m)$ for all $m$.
\end{enumerate}
\end{lem}

\begin{proof} Obvious from Remark~\ref{rems:support}(2).
\end{proof}

\begin{Def} Let $W$ be a  simplicial closed subset  of $X$ with $\codim_{X(n)}W_n\ge q$ for all $n$, and  take $E$  in $\Spt(B^{(q)})$. Let $E^W(X)$ denote the simplicial spectrum
\[
n\mapsto E^{W_n}(X(n)).
\]
\end{Def}

Let $\sS$ be a  support condition on $X$,  supported in codimension $\ge q$. We make $\sS$ into a small category by ordering the $W\in\sS$ by inclusion.  For $E\in\Spt(B^{(q)})$, and $W\subset W'$ an inclusion of simplicial closed subsets of $X$, we have the map of simplicial spectra $i_{W',W}:E^W(X)\to E^{W'}(X)$, giving the functor $W\mapsto E^W(X)$ from  $\sS$  to simplicial spectra. We set
\[
E^\sS(X):=\hocolim_{W\in\sS}E^W(X).
\]

By Lemma~\ref{lem:SimpSupp}, we have the canonical homotopy equivalence
\[
E^\sS(X)(n)\cong \hocolim_{T\in\sS(n)}E^T(X(n)).
\]

\subsection{The homotopy coniveau tower} We introduce our main object of study.

We have the cosimplicial scheme $\Delta^*$, with
\[
\Delta^r=\Spec(\Z[t_0,\ldots, t_r]/\sum_jt_j-1).
\]
The {\em vertices} of $\Delta^r$ are the closed subschemes $v^r_i$ defined by
$t_i=1$, $t_j=0$ for $j\neq i$. A {\em face} of $\Delta^r$ is a closed subscheme
defined by equations of the form $t_{i_1}=\ldots = t_{i_s}=0$.

\begin{Def}\label{Def:SupportCond1}For $X$ in $\Sm/B$, we let $\sS_X^{(q)}(p)$ denote the set of closed subsets
$W$ of $X\times\Delta^p$ such that 
\[
\codim_{X\times F}W\cap(X\times F)\ge q
\]
for all faces $F$ of $\Delta^p$. 
\end{Def}

Clearly, sending $p$ to $\sS_X^{(q)}(p)$ defines a  support condition $\sS_X^{(q)}$ on $X\times\Delta^*$, supported in codimension $\ge q$. We let $X^{(q)}(p)$ be the set of codimension $q$ points $x$ of $X\times\Delta^p$ with closure $\bar{x}\in\sS_X^{(q)}(p)$.

\begin{Def} Take $E\in \Spt(B^{(q)})$.
We  set
\begin{align*}
&E^{(q)}(X,p):=E^{\sS_X^{(q)}(p)}(X\times\Delta^p)\\
&E^{(q)}(X,-):=E^{\sS_X^{(q)}}(X\times\Delta^*).  
\end{align*}
Similarly, if $W$ is a simplicial closed subset of $X\times\Delta^*$, supported in codimension $\ge q$, we write $E^W(X,-)$ for $E^W(X\times\Delta^*)$.  
\end{Def}

Take $q'\ge q$ and let  $E$ be in $\Spt(B^{(q)})$.  We may form the restriction $E\circ \rho_{q,q'}\in \Spt(B^{(q')})$. Similarly, for $X\in\Sm/B$, we have  the simplicial spectrum $(E\circ \rho_{q,q'})^{(q')}(X,-)$, which we write simply as $E^{(q')}(X,-)$.

The inclusions $\sS_X^{(q')}(p)\subset \sS_X^{(q)}(p)$ and the maps
 \begin{align*}
i_{W,W'}:&E^{W'}(X\times\Delta^p)\to E^W(X\times\Delta^p);\\
 &W'\in \sS_X^{(q')}(p), 
W\in \sS_X^{(q)}(p), W'\subset W
\end{align*}
 define the natural map of simplicial spectra
\[
\iota_{q,q'}:E^{(q')}(X,-)\to E^{(q)}(X,-)
\]
with $\iota_{q,q'}\circ\iota_{q',q''}=\iota_{q,q''}$ for $q\le q'\le q''$. This gives us the {\em homotopy coniveau tower} for $E$:
\begin{equation}\label{eqn:HCTower}
\ldots\to E^{(q+r+1)}(X,-)\to E^{(q+r)}(X,-)\to\ldots\to E^{(q)}(X,-).
\end{equation}
This is our main object of study.

\subsection{First properties}\label{subsec:FirstProperties} We give a list of elementary properties of the spectra $E^{(q)}(X,-)$ and the tower \eqref{eqn:HCTower}
\begin{enumerate}
\item Sending $X$ to $E^{(q)}(X,-)$ is functorial for equi-dimensional (e.g. flat) maps $Y\to X$ in $\Sm/B$; the tower \eqref{eqn:HCTower} is natural with respect to equi-dimensional maps.
\item Sending $E$ to $E^{(q)}(X,-)$ is functorial in $E$.
\item The functor $E\mapsto E^{(q)}(X,-)$ sends weak equivalences to weak equivalences, 
and send homotopy (co)fiber sequences to homotopy (co)fiber sequences.
\end{enumerate}

We let $\Sm\ds B$ denote the subcategory of $\Sm/B$ with the same objects, but with morphisms $Y\to X$ being the {\em smooth} $B$-morphisms. We let $E^{(q)}_\sm$ be the presheaf on 
$\Sm\ds B$: $E^{(q)}_\sm(X):=E^{(q)}(X,-)$. We thus have the tower of presheaves on 
 $\Sm\ds B$

\begin{equation}\label{eqn:HCTowerSm}
\ldots\to E^{(q+r+1)}_\sm\to E^{(q+r)}_\sm\to\ldots\to E^{(q)}_\sm
\end{equation}
which is natural in $E\in\Spt(\Sm/B^{(q)})$.

\subsection{Good position} Chow's moving lemma (Theorem~\ref{thm:main}) asserts that the spectrum $E^{(q)}(X,-)$ defined by the support conditions $\sS_X^{(q)}$  is weakly equivalent to another spectrum $E^{(q)}(X,-)_\sC$ defined with respect to a somewhat finer   support condition $\sS_{X,\sC,e}^{(q)}\subset \sS_X^{(p)}$, which we define in this section. We assume from now on that the base-scheme $B$ is regular, noetherian, separated and has Krull dimension at most one. For a $B$-scheme $f:X\to B$ and a point $b\in B$, we set $X_b:=f^{-1}(b)$.

\begin{Def}\label{Def:GoodPosition} Let $X$  be in $\Sm/B$, let $\sC$ be a
finite set of locally closed subsets of $X$, and let
$e:\sC\to\N$ be a function. We assume that each $C\in \sC$ has pure codimension in $X$.   Let $\sS_{X,\sC,e}^{(q)}(p)$ be the subset of $\sS_{X}^{(q)}(p)$
consisting of those $W$ such that, for each $C\in\sC$ and each face $F$ of $\Delta^p$,  we have
\[
\codim_{C\times F}(C\times F)\cap W\ge q-e(C).
\]
\end{Def}

It is more or less obvious that $p\mapsto \sS_{X,\sC,e}^{(q)}(p)$ defines a support condition on $X\times\Delta^*$. We write  $\sS_{X,\sC}^{(q)}(p)$ for $\sS_{X,\sC,0}^{(q)}(p)$.

Given $E\in\Spt(B^{(q)})$, we set
\[
E^{(q)}(X,-)_{\sC,e}:=E^{\sS_{X,\sC,e}^{(q)}}(X\times\Delta^*).
\]
The inclusion $\sS_{X,\sC,e}^{(q)}\to \sS_{X}^{(q)}$ defines the natural map
\[
E^{(q)}(X,-)_{\sC,e}\to E^{(q)}(X,-);
\]
similarly, if $\sC\subset \sC'$ and $e'\le e$, we have the natural map
\[
E^{(q)}(X,-)_{\sC',e'}\to E^{(q)}(X,-)_{\sC,e}.
\]
We write $E^{(q)}(X,-)_{\sC}$ for $E^{(q)}(X,-)_{\sC,0}$.

Given a function
$e:\sC\to\N$, we let $e-1:\sC\to \N$ be the function $(e-1)(C)=\max(e(C)-1,0)$.

Let $f:Y\to X$ be a morphism in $\Sm/B$, $\sC$ a finite set of locally closed subsets of $X$ such that $f^{-1}(C)$ has pure codimension on $Y$ for each $C\in\sC$. We let $f^*\sC$ be the set of irreducible components of the locally closed subsets $f^{-1}(C)$, $C\in \sC$ and for $D\in f^*\sC$, define $f^*(e)(D)$ to be the minimum of the numbers $e(C)$, with $D$ an irreducible component of $f^{-1}(C)$. If $f$ happens to be the identity morphism, then, even though $f^*\sC$ may not be equal to $\sC$, the support conditions $\sS_{X,\sC,e}$ and $\sS_{X,f^*\sC,f^*e}$ are equal.

If $f$ is smooth, the pull-backs by $f\times\id_{\Delta^p}$ defines the map of simplicial sets
\[
f^{-1}:\sS_{X,\sC,e}^{(q)}\to \sS_{Y,f^*\sC,f^*e}^{(q)}
\]
and the map of simplicial spectra
\[
f^*:E^{(q)}(X,-)_{\sC,e}\to E^{(q)}(Y,-)_{f^*\sC,f^*e}.
\]
In particular, we have the presheaf of simplicial spectra (for the Nisnevic topology) on $X$, 
$E^{(p)}(X_\Nis,-)_{\sC,e}$.

We now state the main result of this paper; for a description of the axioms \eqref{eqn:Ax1},
\eqref{eqn:Ax2} and \eqref{eqn:Ax3} we refer the reader to \S\ref{sec:axioms}:

\begin{thm}\label{thm:main} Let $B$ be a regular  noetherian separated scheme of Krull dimension at most one and $X\to B$ a smooth $B$-scheme of finite type. Let $\sC$ be a finite set of irreducible locally closed subsets of $X$, $e:\sC\to\N$ a function and $q\ge0$ an integer. Take $E\in\Spt(B^{(q)})$ satisfying axioms  \eqref{eqn:Ax1},
\eqref{eqn:Ax2} and \eqref{eqn:Ax3}. Suppose that each $C\in\sC$ dominates an irreducible component of $B$. \\

(1) Let $x$ be a point of $X$. Then there is a basis of Nisnevic neighborhoods $p:U\to X$ of $x$, depending only on $x$ and $X$,  such that $E^{(q)}(U,-)_{p^*\sC,p^*e}\to E^{(q)}(U,-)$ is a weak equivalence. In particular,  the map of presheaves 
\[
E^{(q)}(X_\Nis,-)_{\sC,e}\to E^{(q)}(X_\Nis,-)
\]
is a local weak equivalence.\\

(2) Suppose that $B=\Spec k$, $k$ a field, and let $j:U\to X$ be an  affine open subscheme. Then the map of spectra
\[
E^{(q)}(U,-)_{j^*\sC,j^*e}\to E^{(q)}(U,-)
\]
is a weak equivalence.  In particular,  the map of presheaves 
\[
E^{(q)}(X_\Zar,-)_{\sC,e}\to E^{(q)}(X_\Zar,-)
\]
is a local weak equivalence.\\
\end{thm}

\section{Easy moving}\label{sec:EasyMov} The proof of Chow's moving lemma is in two steps, just as the classical case ({\it cf.} \cite{Roberts}): One needs a ``moving by translation" result to handle the case of affine $n$-space, and then one uses the method of the projecting cone to handle the general case. In this section, we deal with the moving by translation portion. This result also allows us to prove the homotopy invariance of the spectra $E^{(p)}(X,-)$, following Bloch's argument for the cycle complexes \cite{AlgCyc}.  In this section, we assume that the base scheme $B$ is of the form $B=\Spec A$, with $A$ a semi-local PID.

\subsection{Some added generality}\label{sec:ExtraGen}
In fact, neither the moving-by-trans\-lation technique nor the proof of homotopy invariance require smoothness, so both are applicable to presheaves of spectra on other categories besides $\Sm/B^{(q)}$.    To adapt to this more general setting, we first need to give a definition of the ``dimension" of a finite-type $B$-scheme.

\begin{Def} Let $B$ be a regular irreducible  scheme of Krull dimension $\le 1$. Let $f:X\to B$ be an irreducible finite type $B$-scheme. If $B=\Spec k$, $k$ a field, set $\dim X=\dim_BX$. If $B$ has Krull dimension 1, let $\eta$ be the generic point of $B$. Set
\[
\dim X:=\begin{cases} \dim_{k(\eta)}X_\eta+1&\quad\text{if }X_\eta\neq\0\\
\dim_{k(b)}X&\quad\text{if } f(X)=b\in B^{(1)}. 
\end{cases}
\]
We extend the definition to the case of reducible $B$ in the obvious way.

If $W\subset X$ is a closed subset of some $X\in\Sch/B$, we write $\dim W\le n$ if each irreducible component $W'$ of $W$ has $\dim W'\le n$.
\end{Def}
Note that, if $W'\subset W\subset X$ are closed subsets of $X$, then
\[
\dim W\le n\Longrightarrow \dim W'\le n.
\]

Let $X$ be in $\Sch/B$, $Y\in\Sm/B$ and $W$ a closed subset of $X\times_BY$. We write $\dimrel_YW\le q$ if for each irreducible component $Y_j$ of $Y$, we have
\[
\dim W\cap X\times_BY_j \le q+\dim_B Y_j.
\]

\begin{Def} Let $X$ be a finite type $B$-scheme. Fix an integer $q$. Let $X/B_{(q)}$ be the category of pairs $(Y,W)$, with $Y$ in $\Sm/B$, quasi-projective over $B$, and $W$ a closed subset of $X\times_BY$ such that $\dimrel_YW\le q$.

A morphism $f:(Y,W)\to (Y',W')$ is a morphism $f:Y\to Y'$ in $\Sm/B$ such that $W\supset(\id_X\times f)^{-1}(W')$.
\end{Def}
\begin{rem}
If  $B$ has pure Krull dimension zero,  then $B/B_{(q)}$ is the full subcategory of $\Sm/B^{(-q)}$ consisting of $(Y,W)$ with $Y$ quasi-projective. If $B$ has pure Krull dimension one, then $B/B_{(q)}$ is the similarly defined full subcategory of $\Sm/B^{(1-q)}$.
\end{rem}

We have the category $\Spt(X/B_{(q)})$ of presheaves of spectra on $X/B_{(q)}$ and the homotopy category $\sH\Spt(X/B_{(q)})$. As before, we write $E^W(Y)$ for $E(Y,W)$ if $E$ is in $\Spt(X/B_{(q)})$ and $(Y,W)$ is in $X/B_{(q)}$.

Let $C$ be a constructible subset of $X\in\Sch/B$ with closure $\bar C$. If $\bar C$ has irreducible components $\bar C_1,\ldots, \bar C_s$, we call the constructible sets $C\cap\bar C_1,\ldots, C\cap\bar C_s$ the irreducible components of $C$ and say that $C$ is irreducible if $s=1$, i.e. if $\bar C$ is irreducible. 

If $C$ is irreducible define $\dim C:=\dim\bar C$; similarly, we  say that $C$ has {\em pure codimension $d$} on $X$ if for each irreducible component $X_i$ of $X$ containing $C$, we have $\codim_{X_i}\bar C=d$. We write this as $\codim_XC=d$. If $C$ has irreducible components $C_1,\ldots, C_s$, we say $\dim C\le m$ if $\dim C_i\le m$ for each $i$. In particular, for $W$ a closed subset of $X$, the statement $\dim W\cap C\le m$ makes sense. Finally, if $C$ is an irreducible constructible subset of $X$, we say that $C$ dominates a component of $B$ if this is so for $\bar C$.

\begin{Def} (1) Let $X$ be in $\Sch/B$. Let $\sS^X_{(q)}(p)$ be the set of closed subsets $W\subset X\times\Delta^p_B$ such that
\[
\dimrel_{F} W\cap(X\times F)\le q
\]
for each face $F$ of $\Delta^p_B$.\\

(2)  Let $X$ be in $\Sch/B$,   let $\sC$ be a finite set of  constructible subsets of $X$, and $e:\sC\to\N$ a function. We assume that each $C\in \sC$ has pure codimension on $X$. \\
\\
Let $\sS^{X,\sC,e}_{(q)}(p)$ be the subset of $\sS^{X}_{(q)}(p)$
consisting of those $W$ such that, for each $C\in\sC$ and each face $F$ of $\Delta^p_B$,  we have
\[
\dimrel_{F}(C\times F)\cap W\le q+e(C) -\codim_XC.
\]
\end{Def}

\begin{Def} Take $X\in\Sch/B$, $E$ in $\Spt(X/B_{(q)})$, $q\in\Z$. Set
\[
E_{(q)}(X,-):=E^{\sS^{X}_{(q)}}(X\times_B\Delta^*).
\]
If we have a finite set $\sC$ of constructible subsets of $X$ (with pure codimension) and a function $e:\sC\to\N$, set
\[
E_{(q)}(X,-)^{\sC,e}:=E^{\sS^{X,\sC,e}_{(q)}}(X\times\Delta^*).
\]
\end{Def}

\begin{rems} \label{rem:LocSuppCond} (1) Suppose that each irreducible component $X_i$ of $X$ has $\dim X_i=d$; we call such an $X$ {\em equi-dimensional of dimension $d$}. Then we can index by codimension,  as in Definitions~\ref{Def:SupportCond1} and \ref{Def:GoodPosition}, and we have
\[
\sS^{X,\sC,e}_{(q)}(p)=\sS_{X,\sC,e}^{(d-q)}(p).
\]

(2) Suppose that $X$ is irreducible and that $X\to B$ is dominant, so we can index using codimension. Let $\sC$ be a finite set of irreducible locally closed subsets of  $X$, $e:\sC\to\N$ a function and $W$ an irreducible  closed subset of $X\times\Delta^n$. Let $i:\eta\to B$ be the generic point of $B$. Suppose that $W\to B$ is dominant and  $i^{-1}(C)$ has pure codimension on $X_b$ for each $C\in\sC$ and each closed point $b\in B$. Then $W$ is in $\sS_{X,\sC,e}^{(q)}(n)$ if  and only if $W_\eta$ is in 
 $\sS_{X_\eta,i^*\sC,i^*(e)}^{(q)}(n)$ and $W_b$ is in $\sS_{X,i_b^*\sC,i_b^*(e)}^{(q-1)}(n)$ for all codimension one points $b$ of $B$.  

Similarly, if  each $C\in \sC$ dominates a component of $B$ and $W\subset X\times\Delta^n$ is arbitrary, then $W$ is in $\sS_{X,\sC,e}^{(q)}(n)$ if  and only if $W_\eta$ is in 
 $\sS_{X,i^*\sC,i^*(e)}^{(q)}(n)$ and $W_b$ is in $\sS_{X,i_b^*\sC,i_b^*(e)}^{(q-1)}(n)$ for all codimenison one points $b$ of $B$.  
 
 (3) The extension from locally closed subsets to constructible subsets of $X$ is only a notational convenience, at least if $X$ is quasi-projective over $B$. Indeed, if $C$ is a constructible subset of $X$ of pure codimension, we can write $C$ as a finite disjoint union of locally closed subsets $C_i\subset X$ of pure codimension $d_i\ge d$; if we set $e(C_i):=e(C)+d_i-d$, we can replace  $\sC=\{C\}, e$ with $\sD:=\{C_i\ |\ C\in \sC\}, e$ and $\sS_{X,\sC,e}= \sS_{X,\sD,e}$. Similarly if $X$ is irreducible and an irreducible constructible   subset $C$ dominates a component of $B$, then we can find a decomposition $C=\amalg_iC_i$ into locally closed subsets such that each  $C_i$ dominates a component of $B$, and we may then apply the comments of (2). The same holds if $X=Y\times_BZ$, with $Z\to B$ quasi-projective, $\sC=\{Y\times C\ |\ C\in \sC_0\}$, $e(Y\times C)=e_0(C)$ for some collection $\sC_0$ of constructible subsets of $Z$ and function $e_0:\sC_0\to\N$.
\end{rems}

\begin{note} For $X$, $E$, $\sC$ and $e$ as above, let $Y$ be in $\Sm/B$ of pure dimension $d$ over $B$. Sending $W\subset (X\times_BY)\times_BT$ to $W\subset X\times_B(Y\times_BT)$ defines the functor
\[
p_Y:X\times_BY/B_{(q+d)}\to X/B_{(q)}
\]
and by pull-back the functor
\[
p_Y^*:\Spt(X/B_{(q)})\to \Spt(X\times_BY/B_{(q+d)}).
\]
We write $E_{(q)}(X\times_BY,-)^{\sC,e}$ for $(p_Y^*E)_{(q+d)}(X\times_BY,-)^{p_2^*\sC,p_2^*e}$. The projection $X\times_BY\to X$ defines the ``base extension" map
\[
p_2^*:E_{(q)}(X,-)^{\sC,e}\to E_{(q)}(X\times_BY,-)^{\sC,e}.
\]

We occasionally need to extend the assignment $Y\mapsto E_{(q)}(X_Y,-)^{\sC,e}$ ($X_Y:=X\times_BY$) to essentially smooth subschemes $U\subset Y$,  i.e., $Y\setminus U$ is a possibly infinite union of closed subsets.  We define 
\[
E_{(q)}(X_U,-)^{\sC,e}:= \hocolim_FE_{(q)}(X_{Y\setminus F},-)^{\sC,e}
\]
where  $F$ runs over closed subsets  of $Y$ contained in $Y\setminus U$.  

In case $p:B'\to B$ is in $\Sm/B^\ess$, and $B'$ is also a regular scheme of Krull dimension one, we have
\[
E_{(q)}(X_{B'},-)^{\sC,e}=E_{(q)}(X_{B'},-)^{p_2^*\sC,p_2^*(e)},
\]
where on the right-hand side, we consider $X_{B'}$ as in $\Sch/B'$, and we extend $E$ to a presheaf on $\Sm/B'_{(q)}$ by taking limits as above. Also, we have the base-extension map
\[
p^*:E_{(q)}(X,-)^{\sC,e}\to E_{(q)}(X_{B'},-)^{\sC,e}.
\]
\end{note}

\begin{rem} In this section $E$ will be a functor $E:X/B_{(q)}^\op\to\Spt$. In case the residue fields of $B$ are  infinite, the results of this section are valid for every such $E$. In the case of finite residue fields, one needs to impose an additional hypothesis on $E$ which allows one to  increase the size of the residue fields at will. As this is a technical point, we postpone the description of this  additional hypothesis (axiom   \eqref{eqn:Ax3})  until \S\ref{sec:axioms}.
\end{rem}

\begin{exs} Our primary examples are:\\
(1) $X$ is in $\Sm/B$ of dimension $d$ over $B$, and $E$ is in $\Sm/B^{(d-q)}$. Then we have $(Y,W\subset X\times_BY)\mapsto E^W(X\times_BY)$, giving the functor
\[
E_X:X/B_{(q)}^\op\to\Spt.
\]
As particular examples, we have the functor $E^{(d-q)}\in\Spt(B^{(d-q)}$ induced by a functor $E\in\Spt(B)$, e.g. $E$ the $K$-theory spectrum $Y\mapsto K(Y)$.\\
(2) $X$ is a quasi-projective $B$-scheme and $E$ is the $G$-theory spectrum with supports. More precisely, let $Y$ be a quasi-projective $B$-scheme and let $\sM_X(Y)$ be the exact category of coherent sheaves on $X\times_BY$ which are flat as $\sO_Y$-modules. We let $G_X(Y)$ denote the $K$-theory spectrum of the exact category $\sM_X(Y)$. Let $G(X\times_BY)$ denote the $K$-theory spectrum of the category of coherent sheaves on $X\times_BY$.

For $Y$ quasi-projective over $B$, each $\sF\in \sM_X(Y)$ admits a surjection $\sE\to \sF$ with $\sE$ a locally free coherent sheaf on $X\times_BY$; if $Y$ is smooth over $B$ and has Krull dimenison $\le d$, then each $\sF$ admits a resolution
\[
0\to\sK\to\sE_d\to\ldots\to\sE_0\to\sF\to0
\]
with the $\sE_j$ locally free coherent sheaves and $\sK$ a flat $\sO_Y$-module. By Quillen's resolution theorem \cite[\S 4, Theorem 3]{Quillen}, this implies that the canonical map $G_X(Y)\to G(X\times_BY)$ is a weak equivalence.

By the flatness condition, $Y\mapsto \sM_X(Y)$ is an exact pseudo-functor on $\Sm/B$. By the usual trick (see e.g. \cite[\S 7, proof of Proposition 2.2]{Quillen}), we replace  $\sM_X(Y)$ with an equivalent exact category and change notation so that $Y\mapsto \sM_X(Y)$ is an exact functor on $\Sm/B$. Thus, we have the functor $Y\mapsto G_X(Y)$ on $\Sm/B$, and the weak equivalence $G_X(Y)\to G(X\times_BY)$, natural in $Y$ for flat morphisms. In particular, if $W\subset X\times_BY$ is a closed subset, Quillen's localization theorem \cite[Thm. 5, \S 5]{Quillen} yields a weak equivalence
\[
G_X^W(Y)\to G^W(X\times_BY)\sim G(W).
\]
We have thus succeeded in extending  (up to weak equivalence)  the assignment 
\[
(Y,W\subset X\times_BY)\mapsto G(W)
\]
to a functor
\[
G_{X,(q)}:X/B_{(q)}\to \Spt.
\]
\ \\
(3) Let $X$ be a finite type $B$-scheme. For $(Y,W)\in X/B_{(q)}$, $Y$ irreducible, let $z_q(Y,W)$ be the free abelian group on the irreducible components $W_i$ of $W$ with
\[
\dim W_i=q+\dim_BY.
\]
We extend to arbitrary $Y$ by taking the direct sum over the irreducible components of $Y$. 

Let $f:Y'\to Y$ be a morphism in $\Sm/B$. Since $f$ has finite Tor-dimenison, the same is true for $\id\times f:X\times_BY'\to X\times_BY$. For $(Y,W), (Y',W')\in X/B_{(q)}$ such that $f$ extends to $f:(Y',W')\to (Y,W)$, and $W_i$ an irreducible component of $W$ as above, each component of $(\id\times f)^{-1}(W_i)$ has the proper dimension, and so the cycle-theoretic pull-back
\[
(\id\times f)^*:z_q(Y,W)\to z_q(Y',W')
\]
is well-defined. We thus have the functor
\[
z_q:X/B_{(q)}\to \Ab;
\]
using the associated Eilenberg-Maclane spectrum gives us the functor
\[
z_q:X/B_{(q)}\to \Spt.
\]
\end{exs}

\subsection{The homotopy}
 Let $v_i^p$ be the vertex
$t_j=0$, $j\ne i$ of $\Delta^p$, and let
$g=(g_1,g_2):[n]\to[p]\times[q]$ be an order-preserving map. We have the affine linear map
$T(g):\Delta^n\to\Delta^p\times\Delta^q$ with $T(g)(v_i^n)=v_{g_1(i)}^p\times v_{g_2(i)}^q$. If $g$ is injective, then $T(g)$ is a closed embedding.
A {\em face} $F$ of $\Delta^p\times\Delta^q$ is a subscheme of the form $T(g)(\Delta^n)$ for an injective $g$.

\begin{Def} Let $\sC$ be a finite set of constructible subsets of $X$, with each $C\in \sC$ having pure codimension on $X$, $e:\sC\to\N$ a function. We let $\sS^{X,\sC,e}_{(q),h}(n)\subset \sS^{X\times\Delta^1}_{(q)}(n)$ be the set of $W\in \sS^{X\times\Delta^1}_{(q)}(n)$ such that
\begin{enumerate}
\item For each face $F$ of $\Delta^1\times\Delta^n$, we have
\[
\dimrel_F W\cap X\times F\le q
\]
\item For each face $F$ of $\Delta^1\times\Delta^n$ and each $C\in \sC$, we have
\[
\dimrel_F W\cap C\times F\le q-\codim_XC+e(C).
\]
\end{enumerate}
\end{Def}
$n\mapsto \sS^{X,\sC,e}_{(q),h}(n)$ defines a support condition on $X\times\Delta^1$; we set
\[
E_{(q)}(X\times\Delta^1,-)^{\sC,e}_h:=E^{\sS^{X,\sC,e}_{(q),h}}(X\times\Delta^1\times\Delta^*).
\]

\begin{rems} (1) Letting $0=\delta_0(\Delta^0), 1=\delta_1(\Delta^0)\in\Delta^1$, the subschemes $0\times \Delta^n$ and $1\times\Delta^n$ are faces of $\Delta^1\times\Delta^n$. Thus the pull-back maps
\[
\delta_0^*,\delta_1^*:E_{(q)}(X\times\Delta^1,-)^{\sC,e}_h\to
E_{(q)}(X,-)^{\sC,e}
\]
are defined.\\
\\
(2) The projection $p_1:X\times\Delta^1\to X$ has
\[
p_1^{-1}(\sS^{X,\sC,e}_{(q)})\subset \sS^{X,\sC,e}_{(q),h}.
\]
Thus we have a well-defined pull-back map
\[
p_1^*:E_{(q)}(X,-)^{\sC,e}\to E_{(q)}(X\times\Delta^1,-)^{\sC,e}_h.
\]
\end{rems}

\begin{lem} \label{lem:Homotopy0}The maps $T(g)$ define a homotopy of 
\[
\delta_0^*, \delta_1^*:E_{(q)}(X\times\Delta^1,-)^{\sC,e}_h\to
E_{(q)}(X,-)^{\sC,e}, 
\]
natural in  $(\sC,e)$ and $E$, and natural with respect to base-change $X\mapsto X_{B'}$ for $B$-schemes $B'\to B$ in $\Sm/B^\ess$.
\end{lem}
\begin{proof} This follows immediately from the homotopy construction of  \S\ref{Triang}.
\end{proof}

\subsection{Action by a group-scheme}\label{subsec:Action} The ``easy" moving lemma requires a
transitive action by a connected linear algebraic group. We will only need the case of $\A^n$, with the action of translation on itself.

Let $\sG$ be a group scheme over $B$, and let $X\in\Sch/B$ be a $B$-scheme with a (right) $\sG$-action
$\rho:X\times_B\sG\to X$. Let $B'\to B$ be in $\Sm/B^\ess$, $\psi:\A^1_{B'}\to\sG$ a $B$-morphism with
$\psi(0)=\id_{\sG}$. We let $\phi:X_{B'}\times\Delta^1\to X$ be the composition
\[
X_{B'}\times\Delta^1\xrightarrow{\sim} X\times_B\A^1_{B'}\xrightarrow{\id\times\psi}
X\times_B\sG\xrightarrow{\rho}X,
\]
where the isomorphism $X_{B'}\times\Delta^1\xrightarrow{\sim} X\times_B\A^1_{B'}$ is induced by the
isomorphism $(t_0,t_1)\mapsto t_0$ of $\Delta^1$ with $\A^1$. We have the composition
\[
X_{B'}\times\Delta^p 
\arr{\cdot\psi(1)\times\id}X_{B'}\times\Delta^p \arr{\pi\times\id}
X\times\Delta^p 
\]
which we denote by $\phi(1,p)$. Note that
\begin{align*}
&\phi(1,p)=\phi\circ(\id\times\delta_1\times \id)\\
&\pi=\phi\circ(\id\times\delta_0\times \id)
\end{align*}
where $\pi:X_{B'}\times\Delta^p\to X\times\Delta^p$ is the projection.

We will take $\sG$ to be the additive group $\A^n$, and $\psi$ to be the linear map
$\psi_x$, $\psi_x(t)=t\cdot x$, where $x:B'\to\A^n$ is a $B$-morphism. We take $X=Y\times\A^n$, $Y\in\Sch/B$, with $\sG$ acting on $\A^n$ by translation and on $Y$ trivially.

\begin{note}
Let $A(x_1,\ldots, x_n)$ be the localization of the polynomial ring $A[x_1,\ldots, x_n]$ with respect to the multiplicatively closed set of $f=\sum_Ia_Ix^I$ such that the ideal in $A$ generated by the coefficients $a_I$ is the unit ideal. 

The inclusion $A\to A(x_1,\ldots, x_n)$ is a faithfully flat extension of semi-local PID's. We set
\[
B(x):=\Spec A(x_1,\ldots, x_n), 
\]
and let $x:B(x)\to \A^n$ be the morphism with $x^*:A[x_1,\ldots,  x_n]\to A(x_1,\ldots,x_n)$ the evident inclusion. Note that $B(x)$ is in $\Sm/B^\ess$.
\end{note}

\begin{lem}\label{GenPosLem1}  Let $\sC_0$ be a finite set of irreducible constructible subsets of $\A^n_B$, $e:\sC_0 \to\N$ a
function and let $Y$ be in $\Sch/B$. Let $X=Y\times_B\A^n$, $\sC=p_2^*\sC_0$ and $e=p_2^*(e_0)$.  Let 
$W$ be  in $\sS_X^{(q)}(p)$, and let $\psi:=\psi_x:B(x)\to\A^n$. Suppose that each $C\in\sC_0$ dominates an irreducible component of $B$. Then
\begin{enumerate}
\item $\phi(1)^{-1}(W)$ is in
$\sS^{B(x)\times_BX,p_2^*\sC}_{(q)}(p)$.
\item Suppose that $W$ is in $\sS^{X,\sC,e}_{(q)}(p)$. Then $\phi^{-1}(W)$ is in $\sS^{X_{B(x)},\sC,e}_{(q),h}(p)$.
\end{enumerate}
\end{lem}

\begin{proof} By Remark~\ref{rem:LocSuppCond}(3), we may assume that $\sC_0$ is a finite set of irreducible, locally closed subsets of $\A^n_B$, such that each $C\in\sC_0$ dominates an irreducible component of $B$.  Since the translation action preserves all the irreducible components of $X$, we may assume that $X$ is irreducible. If $d=\dim X$, then $\sS^X_{(q)}=\sS_X^{(d-q)}$, etc., so we may index with codimension. 

In case $A$ is a field, or if $X=X_b$ for some $b\in B^{(1)}$, the result follows directly from \cite[Lemma 2.2]{AlgCyc}.

Suppose that $A$ is not a field and that $X\to B$ is dominant. Let $\eta$ be the generic point of $B$.  Since the inclusion $i:\eta\to B$ is flat, $W_\eta$ is in $\sS_{X_\eta}^{(q)}$, and if $W$ is in $\sS_{X,\sC,e}^{(q)}$, then $W_\eta$ is in $\sS_{X_\eta,i^*\sC,i^*(e)}^{(q)}$.  From the case of a field, it follows that (1) and (2) are valid with $\eta$ replacing $B$, and $i^*\sC, i^*(e)$ replacing $\sC, e$.

For a codimension one point $i_b:b\to B$, we have the set of irreducible constructible subsets $i_b^*\sC$ of $X_b$ and the  function $i_b^*e:i_b^*\sC\to\N$. 

 To complete the proof it suffices by Remark~\ref{rem:LocSuppCond}  to show that, for each codimension one
point $i_b:b\to B$ of $B$ that
\begin{enumerate}
\item   $\phi(1)_b^{-1}(W_b)$ is in
$\sS_{X_{b(x)},i_b^*\sC}^{(q-1)}(p)$.
\item Suppose that $W$ is in  $\sS_{X,\sC,e}^{(q)}(p)$.
Then $\phi^{-1}(W_b)$ is in \\$\sS_{X_{b(x)},i_b^*\sC,i_b^*(e)}^{(q-1),h}(p)$.
\end{enumerate}
By Remark~\ref{rem:LocSuppCond}
\[
W\in\sS_{X,\sC,e}^{(q)}\Longrightarrow W_b\in \sS_{X_b,i_b^*\sC,i_b^*e}^{(q-1)}
\]
and similarly without the $(\sC,e)$.  We note that  $B(x)\to B$ is a bijection on points and that the $b(x)$-valued point $x\times_Bb$ of $\A^n_b$ has $\A^n_b$ as Zariski closure.  Thus (1) and (2)  follow from the case of fields.  
\end{proof}

We  let $[E_{(q)}/E_{(q)}^{\sC, e}](X)$ be the homotopy cofiber of the map of 
simplicial spectra 
\[
E_{(q)}(X,-)^{\sC,e}\to E_{(q)}(X,-).
\]
We similarly let $[E_{(q)}/E_{(q)}^{\sC, e}]_h(X\times\Delta^1)$ be the homotopy cofiber of the map of 
simplicial spectra 
\[
E_{(q)}(X\times\Delta^1,-)^{\sC,e}_h\to E_{(q)}(X\times\Delta^1,-)_h.
\]
Let $\pi:X_{B(x)}\to X$ denote the projection.

\begin{lem}\label{lem:TransMovLem} Let $X$, $\sC$ and $e$  be as in Lemma~\ref{GenPosLem1}.    Then, for each $q\in\Z$, the base-change map
\[
\pi^*: [E_{(q)}/E_{(q)}^{\sC, e}](X) \to [E_{(q)}/E_{(q)}^{\sC, e}](X_{B(x)})
\]
is homotopic to zero.
\end{lem}

\begin{proof}   Consider the commutative diagram, which is well-defined by Lemma~\ref{GenPosLem1}(2),
\[
\xymatrix{
[E_{(q)}/E_{(q)}^{\sC, e}](X)\ar[r]^-{\phi^*}\ar@<-3pt>[d]_-{\pi^*}\ar@<3pt>[d]^{\phi(1)^*}&[E_{(q)}/E_{(q)}^{\sC, e}]_h(X\times\Delta^1_{B(x)})\ar@<-3pt>[d]_{\delta_0^*}\ar@<3pt>[d]^{\delta_1^*}\\
[E_{(q)}/E_{(q)}^{\sC, e}](X_{B(x)})\ar@{=}[r]&[E_{(q)}/E_{(q)}^{\sC, e}](X_{B(x)}).
}
\]
By Lemma~\ref{lem:Homotopy0},  $\phi(1)^*$ and $\pi^*$ are homotopic, but by Lemma~\ref{GenPosLem1}(1), $\phi(1)^*$ is homotopic to zero.
\end{proof}

\begin{prop}[Moving by translation] \label{TransMovProp} Let $\sC_0$ be a finite set of irreducible constructible subsets of $\A^n_B$, $e:\sC_0 \to\N$ a
function and let $Y$ be in $\Sch/B$. Let $X=Y\times_B\A^n$, $\sC=p_2^*\sC_0$ and $e=p_2^*(e_0)$.    In case $B$ has a point $b$ with finite residue field $k(b)$, we assume that $E$ satisfies axiom   \eqref{eqn:Ax3}. Suppose further that each $C\in\sC_0$ dominates a component of $B$. Then, for each $q\in\Z$, the map
\[
E_{(q)}(X,-)^{\sC,e}\to E_{(q)}(X,-) 
\]
is a weak equivalence.
\end{prop}

\begin{proof} It suffices to show that
$[E_{(q)}/E_{(q)}^{\sC, e}](X)$ is weakly equivalent to a point. By Lemma~\ref{lem:TransMovLem}, it suffices to show that the base-change map
\begin{equation}\label{TrPropEqn1}
\pi^*:[E_{(q)}/E_{(q)}^{\sC, e}](X)\to [E_{(q)}/E_{(q)}^{\sC, e}](X_{B(x)}) 
\end{equation}
is injective on homotopy groups.

For this, we note that the scheme $B(x)$ is a filtered inverse limit of open subschemes
$U$ of
$\A^n_B$, with $U$ faithfully flat and of finite type over $B$. As
\[
[E_{(q)}/E_{(q)}^{\sC, e}](X_{B(x)}) =\hocolim_{B(x)\subset U\subset \A^n_B}
[E_{(q)}/E_{(q)}^{\sC, e}](X_U)
\]
by definition, it suffices to show that the map
\[
p_2^*:[E_{(q)}/E_{(q)}^{\sC, e}](X)\to [E_{(q)}/E_{(q)}^{\sC,e}](X_U) 
\]
is injective on homotopy groups for each such $U$. 

For this, it suffices to show that, for each $U$ and each collection $(W',W,V',V)$ of simplicial closed subsets with
\begin{enumerate}
\item $W'\subset W$, $V'\subset V$, $p_2^*W'\subset V'$ and $p_2^*W\subset V$
\item $W'\in \sS_{(q)}^{X,\sC,e}$, $W\in  \sS_{(q)}^X$,
$V'\in \sS_{(q)}^{X_U, \sC,e}$, $V\in  \sS_{(q)}^{X_U}$
\end{enumerate}
the map
\begin{equation}\label{TrPropEqn2}
p_2^*:[E^W(X\times\Delta^*)/E^{W'}(X\times\Delta^*)]\to [E^{V'}(X_U\times\Delta^*)/E^V(X_U\times\Delta^*)]
\end{equation}
has the property that, if $p_2^*(x)=0$ for $x\in \pi_m([E^W(X\times\Delta^*)/E^{W'}(X\times\Delta^*)])$, then $x$ goes to zero in $\pi_m([E_{(q)}/E_{(q)}^{\sC, e}](X))$. We may also assume that each of $W',W,V',V$ are generated by finitely many closed subsets $W_i',W_i,V_i',V_i$ of $X\times\Delta^{n_i}$, $X_U\times\Delta^{n_i}$, $i=1,\ldots, s$.

Suppose first that each residue field $k(b)$, $b\in B$, is infinite. As $B$ is semi-local, and $U$ is an open subscheme of $\A^n_B$ with $\A^n_{k(b)}\cap U\neq\0$ for all $b$, the projection $U\to B$ admits a section $s:B\to U$ such that 
\begin{align*}
&(s\times\id)^{-1}(V')\in  \sS_{(q)}^{\sC,e,X}\\
&(s\times\id)^{-1}(V)\in  \sS_{(q)}^X
\end{align*}
Thus we have the map  
\[
s^*:[E^{V'}(X_U\times\Delta^*)/E^V(X_U\times\Delta^*)]\to  [E_{(q)}/E_{(q)}^{\sC,e}](X_U)
\]
 with $s^*\circ p_2*$ the canonical map $[E^W(X\times\Delta^*)/E^{W'}(X\times\Delta^*)]\to [E_{(q)}/E_{(q)}^{\sC,e}](X)$, which suffices for our purpose.

If some $k(b)$ is finite, then $E$ satisfies axiom  \eqref{eqn:Ax3}.  Take distinct primes $p_1$, $p_2$, an integer $n\ge1$ and let $k(b)\to k(b_1)$, $k(b)\to k(b_2)$ be field extensions of respective degree $p_1^n$, $p_2^n$, for each $b$ with $k(b)$ finite. There exist finite Galois \'etale $B$-schemes $f_1:B_1\to B$, $f_2:B_2\to B$, of  relatively prime degrees $n_1$, $n_2$
over $B$, such that, for each $b$ such that  $k(b)$ is finite,  each point $c\in B_j\times_Bb$ has $k(c)\supset k(b_j)$  (see \cite[Lemma 6.1]{Loc}).  By axiom  \eqref{eqn:Ax3}, the maps $E \to f_{j*}^{G_j}f_j^*E$ are weak homotopy equivalences after inverting $n_j$. One can then argue as above to show that $p_2^*$ is injective on homotopy groups.

By Lemma~\ref{lem:TransMovLem}, $p_2^*$ is zero on homotopy groups, whence the result.
\end{proof}

\begin{thm}[Homotopy invariance] Let $X$ be in $\Sch/B$ and let  $E:X/B_{(q)}^\op\to\Spt$ be in $\Spt(X/B_{(q)})$. If $k(b)$ is finite for some $b\in B$, we assume in addition that $E$ satisfies axiom \eqref{eqn:Ax3}. Then for $X\in\Sch/B$ the map
\[
p^*:E_{(q)}(X,-)\to E_{(q)}(X\times\A^1,-)
\]
is a weak equivalence.
\end{thm}

\begin{proof}  One  repeats the argument used to prove \cite[Theorem 2.1]{AlgCyc}, replacing homology with homotopy, chain homotopy equivalence with homotopy equivalence and using 
Lemma~\ref{lem:TransMovLem} in place of   \cite[Lemma 2.3]{AlgCyc}. For the reader's convenience, we give a sketch.

The canonical map $E_{(q)}(X\times\A^1,-)_{\{X\times0\}}\to E_{(q)}(X\times\A^1,-)$ is a weak equivalence by Proposition~\ref{TransMovProp}, and $p^*$ factors through 
\[
p^*:E_{(q)}(X,-)\to E_{(q)}(X\times\A^1,-)_{\{X\times0\}}. 
\]
As this map is split by $i_0^*$, $p^*$ is injective on $\pi_n$, so we need only show surjectivity.

One first shows:
\begin{lem} Let $i_0, i_1:X\to X\times\A^1$ be the 0- and 1-sections, respectively. Then the maps
\[
i_0^*, i_1^*:E_{(q)}(X\times\A^1,-)_{\{X\times0, X\times1\}}\to
E_{(q)}(X,-)
\]
induce the same map on homotopy groups $\pi_n$, $n\in\Z$.
\end{lem}

\begin{proof}[proof of lemma]    Let $\A^1$ act on $\Delta^1$ by translation: $t\cdot(t_0,t_1)=(t_0+t,t_1-t)$. Let $\psi:\A^1_{B(x)}\to \A^1:=\Spec A[t]$ be the map with $\psi^*(t)=xt$. As in the proof of Lemma~\ref{lem:TransMovLem}, the maps
\[
\phi(1)^*,\pi^*:E_{(q)}(X\times\A^1,-)_{\{X\times0, X\times1\}}\to 
E_{(q)}(X\times\A^1_{B(x)},-)_{\{X\times0, X\times1\}}
\]
are homotopic ($\pi^*$ the base-change morphism), hence $\pi^*\circ i_\epsilon^*$ and $ i_\epsilon^*\circ\phi(1)^*$ are homotopic, $\epsilon=0,1$.

Identifying $\A^1$ with $\Delta^1$ by $(t_0,t_1)\mapsto t_0$,  one shows as in the proof of  Lemma~\ref{GenPosLem1} that $\phi(1)^{-1}$ maps $\sS^{X\times\Delta^1}_{(q)}$ into $\sS^{X_{B(x)}}_{(q),h}$, giving us the map
\[
\phi(1)^*:E_{(q)}(X\times\A^1,-)_{\{X\times0, X\times1\}}\to 
E_{(q)}(X\times\Delta^1_{B(x)},-)_h
\]
By Lemma~\ref{lem:Homotopy0}, the maps $\delta_0^*\circ\phi(1)^*$ and $\delta_1^*\circ\phi(1)^*$ are homotopic. Thus $\pi^*\circ i_0^*$ and $\pi^*\circ i_1^*$ are homotopic. But $\pi^*$ is injective on $\pi_n$ ({\it cf.} proof of Proposition~\ref{TransMovProp}), whence the lemma.
\end{proof}

The next step is to consider the multiplication map
\[
\mu:\A^1\times\A^1\to \A^1
\]
$\mu(x,y)=xy$; note that $\mu$ is flat. We have the commutative diagram
\[
\xymatrix{
E_{(q)}(X\times\A^1,-)_{X\times0}\ar[r]^-{\mu^*}\ar@<3pt>[d]^\iota\ar@<-3pt>[d]_{p^*i_0^*}&E_{(q)}(X\times\A^1\times\A^1,-)_{X\times\A^1\times\{0,1\}}\ar@<-3pt>[d]_{i_0^*}\ar@<3pt>[d]^{i_1^*}\\
E_{(q)}(X\times\A^1,-)\ar@{=}[r]&E_{(q)}(X\times\A^1,-)}
\]
where $\iota:E_{(q)}(X\times\A^1,-)_{X\times0}\to E_{(q)}(X\times\A^1,-)$ is the canonical map. Using the lemma above and Proposition~\ref{TransMovProp}, it follows that $p^*i_0^*$ is an isomorphism on $\pi_n$, hence $p^*$ is surjective.
\end{proof}

\section{The axioms}\label{sec:axioms}  We have seen that the moving-by-translation results are valid for an arbitrary functor $E:X/B_{(q)}^\op\to\Spt$ (except in the case of finite residue fields, which requires the axiom  \eqref{eqn:Ax3}). To proceed further, however, we will need to have some additional properties of $E$: localization and Nisnevic excision. In this section, we list the axioms we need.

\subsection{Localization and excision}

Let $E$ be in $\Spt(B^{(q)})$. We  consider the following conditions on $E$:\\
\begin{equation}\label{eqn:Ax1}
\vbox{\ \\
Axiom 1. (Localization) 
\begin{enumera}
\item Let $\Sm/B^\0$ be the full subcategory of $\Sm/B^{(q)}$ with objects $(Y,\0)$. The restriction of $E$ to  $\Sm/B^\0$ is contractible. We fix a contraction $h$.
\item Let $(Y,W)$ be in $\Sm/B^{(q)}$, let $W'\subset W$ be a closed subset and let $j:Y\setminus W'\to Y$ be the inclusion. The sequence
\[
E^{W'}(Y)\xrightarrow{i_{W',W}}E^W(Y)\xrightarrow{j^*}E^{W\setminus W'}(Y\setminus W')
\]
together with $h$ gives a canonical morphism 
\[
\phi_{W',W,Y}:E^{W\setminus W'}(Y\setminus W')\to \cofib(i_{W',W})
\]
Then $\phi_{W',W,Y}$ is a weak equivalence.
\end{enumera}}
\end{equation}
To explain: We have the commutative diagram
\[
\xymatrix{
E^{W'}(Y)\ar[r]^{i_{W',W}}\ar[d]_{j^*}&E^W(Y)\ar[d]^{j^*}\\
E^\0(Y)\ar[r]_{i_{\0,W\setminus W'}}&E^{W\setminus W'}(Y\setminus W')}
\]
The contraction $h$ gives a contraction of $j^*\circ i_{W',W}$, which in turn induces the map $\phi_{W',W,Y}$. In particular, $\phi_{W',W,Y}$ is natural with respect to the category of morphisms $(Y,W)\to (Y,W')$ in $\Sm/B^{(q)}$.

\begin{equation}\label{eqn:Ax2}
\vbox{\ \\
Axiom 2. (Nisnevic excision) Let $f:Y\to X$ be an \'etale map in $\Sm/B$, $f:(Y,W')\to (X,W)$ an extension of $f$ to a morphism in $\Sm/B^{(q)}$  with $W'=f^{-1}(W)$. Suppose that $f:W'\to W$ is an isomorphism (of schemes with the reduced structure). Then
\[
f^*:E^W(X)\to E^{W'}(Y)
\]
is a weak equivalence.}
\end{equation}

\begin{lem}\label{lem:Add}  Let $E$ be  in $\Spt(B^{(q)})$, let $(Y,W)$ be in $\Sm/B^{(q)}$ such that $W$ is a disjiont union of closed subsets, $W=W'\amalg W''$. Suppose that $E$ satisfies axiom 2. Then the natural map
\[
i_{W',W}+i_{W'',W}:E^{W'}(Y)\oplus E^{W''}(Y)\to E^W(Y)
\]
is a weak equivalence.
\end{lem}

\begin{proof} We have the \'etale morphism $p:(Y\setminus W') \amalg (Y\setminus W'')\to Y$, induced by the inclusions. Clearly $p^{-1}(W)\to W$ is an isomorphism, so
\[
p^*:E^W(Y)\to E^{p^{-1}(W)}((Y\setminus W') \amalg (Y\setminus W''))
\]
is a weak equivalence by axiom 2. As 
\[
((Y\setminus W') \amalg (Y\setminus W''), p^{-1}(W))=(Y\setminus W', W'') \amalg (Y\setminus W'', W')
\]
and $E$ is a presheaf, the map
\[
E^{W''}(Y\setminus W') \oplus E^{W'}(Y\setminus W'')\to
E^{p^{-1}(W)}((Y\setminus W') \amalg (Y\setminus W''))
\]
is a weak equivalence. Using axiom 2 again, the restriction map
\[
E^{W''}(Y) \oplus E^{W'}(Y)\to
E^{W''}(Y\setminus W') \oplus E^{W'}(Y\setminus W'')
\]
is a weak equivalence. The commutativity of the diagram
\[
\xymatrix{
E^W(Y)\ar[r]&E^{p^{-1}(W)}((Y\setminus W') \amalg (Y\setminus W''))\\
E^{W''}(Y) \oplus E^{W'}(Y)\ar[r]\ar[u]&
E^{W''}(Y\setminus W') \oplus E^{W'}(Y\setminus W'')\ar[u]}
\]
completes the proof.
\end{proof}
We let $p_{W'}:E^W(Y)\to E^{W'}(Y)$ be the composition in $\SH$
\[
E^W(Y)\xrightarrow{(i_{W',W}+i_{W'',W})^{-1}}E^{W'}(Y)\oplus E^{W''}(Y)\xrightarrow{p_1}
E^{W'}(Y).
\]

We introduce one final axiom, to handle the case of finite residue fields. Suppose we have a finite  \'etale Galois cover $\pi:B'\to B$ with group $G$. Given $E\in\Spt(B^{(q)})$, define $\pi_*^G\pi^*E$ by
\[
\pi_*^G\pi^*E^W(X):= (E^{W_{B'}}(X_{B'}))^G,
\]
where $(-)^G$ denotes a functorial model for the $G$ homotopy fixed point spectrum. We have as well the natural transformation
\[
\pi^*:E\to \pi_*^G\pi^*E.
\] 

\begin{equation}\label{eqn:Ax3}
\vbox{\ \\
Axiom 3.  Suppose that $B$ has some closed point $b$ with finite residue field $k(b)$. Let $\pi:B'\to B$ be a finite Galois \'etale cover of $B$ with group $G$. Then after inverting $|G|$, the natural transformation $\pi^*:E\to \pi^G_*\pi^* E$ is a weak equivalence.
}
\end{equation}

\section{Moving by the projecting cone} \label{sec:Cone} The method of moving by translation takes care of the case $X=\A^n$; for a general smooth affine $B$-scheme, we need to apply the classical method of the projecting cone. Since most of the known results ({\it cf.} \cite{Roberts}, \cite[Part I, Chap. II, \S3]{MixMot}) are restricted to the case of varieties over a field, we will need to prove some technical results which enable the extension to base-schemes $B=\Spec A$, $A$ a semi-local PID.

\subsection{The projecting cone}\label{ProjCone} Although we will eventually need to assume that $B=\Spec A$, with $A$ a semi-local PID, we will work in a somewhat more general setting for a while, assuming only that $A$ is semi-local and noetherian.

 Let $i:X\to\A^n_B$ be a closed subscheme of $\A^n_B$ with $X\in\Sm/B$. We assume that $X$ is equi-dimensional over $B$; let $m=\dim_BX$. Let  $\bar X\subset\P^n_B$ be the closure of $X$ in $\P^n_B$ and let $\P^{n-1}_B(\infty)=\P^n_B\setminus\A^n_B$. Using homogeneous coordinates $X_0,\ldots, X_n$ for $\P^n$, $\P^{n-1}(\infty)$ is defined by $X_0=0$.

Let $M_{n\times m}$ be the affine space $\Spec B[\{X_{ij}\ |\ 1\le i\le n, 1\le j\le m\}]$; we view $M_{n\times m}\times\A^m$ as the parameter space for affine-linear maps $\pi:\A^n\to \A^m$ by defining
\[
\pi_{(x_{ij}), (v_j)}(y_1,\ldots, y_n)= (y_1,\ldots,y_n)\cdot (x_{ij})+(v_1,\ldots, v_m).
\]

Let $\sV\subset M_{n\times m}$ be the matrices of rank $m$. Each $B'$-valued point $x:=(x_{ij})$ of $\sV$ yields the codimension $m$ linear subspace $L(x)\subset \P^{n-1}_{B'}(\infty)$ defined by the equations
\[
\sum_ix_{i1}X_i=\ldots=\sum_ix_{im}X_i=0.
\]
We let $\sV_X\subset \sV$ be the open subscheme with $L(x)\cap \bar{X}=\0$ for $x\in \sV_X$ and set
\[
\sU_X:=\sV_X\times\A^m.
\]
For $\pi\in \sU_X$, let $\pi_X:X\to\A^m$ be the restriction of $\pi$ to $X$.

For $x\in X$, we let $\sU_X(x)$ be the open subscheme of $\sU_X$ of those $\pi$ for which $\pi_X$ is \'etale over a neighborhood of $\pi(x)$.
 
Let $B'\to B$ be a $B$-scheme and $y:=(x,v):B'\to \sU_X$ a $B$-morphism. We write $L(y)$ for $L(x)$. The affine-linear projection  $\pi_y:\A^n_{B'}\to \A^m_{B'}$ extends to the linear projection
$\bar{\pi}_{y}:\P^n_{B'}\setminus L(x)\to \P^m_{B'}$.  Since $L(y)\cap \bar{X}_{B'}=\0$, the map  $\bar{\pi}_y$ restricts to a finite surjective morphism $\bar{X}_{B'}\to \P^m_{B'}$. Thus  the restriction of $\pi_y$ to
\[
\pi_{y,X}:X_{B'}\to\A^m_{B'}
\]
is finite and surjective and the maps
\[
\pi_{y,X}\times\id:X_{B'}\times\Delta^p\to\A^m_{B'}\times\Delta^p
\]
are also finite and surjective for each $p$.

We make the following basic assumptions:

\begin{equation}\label{eqn:Finiteness}
\vbox{
\begin{enumerate}
\item  The structure morphism  $\sU_X\to B$  is surjective.
\item For  $x\in X$, $\sU_X(x)$ is non-empty.
\end{enumerate}
}
\end{equation}

\smallskip
\begin{rems} (1) In case $\dim B\le 1$,   \eqref{eqn:Finiteness}(1) is equivalent to the assumption that $\P^{n-1}_B(\infty)\cap\bar{X}$ is equi-dimensional (of relative dimension $\dim_BX-1$) over $B$, equivalently, for each $b\in B^{(1)}$, 
\[
(\bar{X})_b=\overline{X_b}
\]
where $\overline{X_b}$ is the closure of $X_b$ in $\P^n_b$. In particular, \eqref{eqn:Finiteness}(1) holds if $A$ is a field.\\

(2) It is easy to construct $X\in\Sm/B$ which is affine and equi-dimensional over $B$, of dimension say $d$ over $B$, but does not admit {\em any} finite surjective morphism $X_B\to \A^d_B$, even if all residue fields of $B$ are infinite. However, we will show  (Theorem~\ref{thm:GoodCpt}) that this property holds locally on $(X,B)$ in the Nisnevic topology.\\

(3) Suppose $i_0:X\to \A^n_B$ satisfies \eqref{eqn:Finiteness}(1). Let $i_d:\A^n\to \A^N=\Spec \Z[Y_1,\ldots, Y_N]$ be the embedding with $\{i_d^*(Y_j)\}$ a $\Z$-basis of the monomials of degree $d$ for some $d>1$. Then $i:=(i_d\times_\id B)\circ i_0$ satisfies \eqref{eqn:Finiteness}(2). Also, if all residue fields of $B$ have characteristic zero, then \eqref{eqn:Finiteness}(2) is true for all embeddings $i$ which satisfy \eqref{eqn:Finiteness}(1). See for example \cite[XI 2.1]{SGA4}.
\end{rems}

\subsection{Notation and elementary properties}\label{sec:Notation}

We assume in this section that \eqref{eqn:Finiteness} holds.

In particular, $\sU_X$ is an open subscheme of an affine space $\A^r_B:=\Spec B[y_1,\ldots, y_r]$ over $B$, of finite type and faithfully flat over $B$.  Generalizing the construction of \S\ref{subsec:Action}, we let $A(y)$ denote the localization of $A[y_1,\ldots, y_r]$ with respect to the set of $f=\sum_Ia_Iy^I$ such that the ideal in $A$ generated by the $a_I$ is the unit ideal. This gives us the semi-local noetherian ring  $A(y)$, and the $B$-scheme $B(y):=\Spec A(y)$, together with a $B$-morphism $y:B(y)\to \sU_X$ having Zariski dense image; in fact, for each point $b$ of $B$, the restriction $y_b:b(y)\to\sU_X\times_Bb$ has Zariski dense image.  

\begin{note} For each $p$, let $f_p:X_{B(y)}\times\Delta^p\to \A^m_{B(y)}\times\Delta^p$ denote the morphism $\pi_{y,X}\times\id$. We let $R_y\subset X_{B(y)}$ denote the ramification locus of $\pi_{y,X}$.

 If $W\subset X\times\Delta^p$ is a pure codimension $q$ closed subset with irreducible components $W_1,\ldots, W_s$, let $\cyc(W):=\sum_{i=1}^sW_i$ and set
\[
W^c:=\Supp\big(f_p^*(f_{p*}(\cyc(W)))-\cyc(W)\big).
\]
Let $W^+\subset X_{B(y)}\times\Delta^p$ denote the closure of $f_p^{-1}(f_p(W_{B(y)}))\setminus W_{B(y)}$. 
\end{note}

\begin{rem} For $W\subset X\times\Delta^p$ a pure codimension $q$ closed subset,  $W^+$ is also a pure codimenion $q$ closed subset of $X\times\Delta^p$, and $W^+\subset W^c$. In addition $W^+\cup W=W^c\cup W$, and $W^+=W^c$ if and only if $W^c$ and $W$ have no irreducible components in common.
\end{rem}

\begin{lem} \label{lem:MovLem1} Suppose that $A$ is a field.  Let $X$ be in $\Sm/B$ with a closed embedding $i:X\to\A^n_B$, $\sC$ a finite set of irreducible locally closed subsets of $X$, and $e:\sC\to\N$ a function. Let $q>0$ be an integer.
\begin{enumerate}
\item  Let $Z$ be an effective codimension $q$ cycle on $X$ with support $|Z|\in \sS^{(q)}_{X,\sC,e}(0)$. Then 
 $\pi_{y,X}^*(\pi_{y,X*}(Z))=Z+Z'$, with $Z'$ effective. In addition, $Z'$ has no component in common with $Z$, and  $|Z'|$ is in $\sS^{(q)}_{X,\sC,e-1}(0)$.
\item Let $W\subset X\times\Delta^p$ be a closed subset in $\sS^{(q)}_{X,\sC,e}(p)$. Then
\begin{enumerate}
\item $W_{B(y)}\cap (R_y\times\Delta^p)$ is in $\sS_{X_{B(y)},\sC,e}^{(q+1)}(p)$.
\item  $W^c$ is in $\sS_{X_{B(y)},\sC,e-1}^{(q)}(p)$.
\end{enumerate}
\end{enumerate}
\end{lem}

\begin{proof}  We note that $B(y)$ is Zariski dense in $\sU_X$. The result is a direct consequence (after a translation in the notation) of \cite[Part I, Chap. II, Lemma 3.5.4 and Lemma 3.5.6]{MixMot}.
\end{proof}

\subsection{Birationality} We now proceed to analyze the maps $f_p:X\times\Delta^p\to\A^m_{B'}\times\Delta^p$. For most $W\in\sS^{(q)}_{X,\sC,e}(p)$, the map $f_p:W\to f_p(W)$ is birational, but there are some special cases that need to be excluded. From now on, we assume that $B=\Spec A$, with $A$ a semi-local PID.

\begin{Def} For a $B$-scheme $B'\to B$ and   a codimension $q$ closed subset $Z$ of $X_{B'}\times\Delta^p$, we call $Z$ {\em induced} if $Z\subset X\times Z'$ for some codimension $q$ closed subset $Z'$ of $\Delta^p_B$.  If $W\subset X_{B'}\times\Delta^p$ is a pure codimension $q$ closed subset,  let $\Ind(W)$ denote the union of the irreducible components $Z$ of $W$ which are induced.
\end{Def}

\begin{rems}\label{rem:Induced} (1) If $Z\subset X\times\Delta^p$ is induced, then $(\pi_{y,X}\times\id)(Z)$ is
an induced closed subset of $\A^m_{B(y)}\times\Delta^p$.\\

(2) Suppose $Z\subset X\times\Delta^p$ is induced and in $\sS^{(q)}_X(p)$. Then for each finite set of irreducible locally closed subsets $\sC$, $Z$ is in $\sS^{(q)}_{X,\sC}(p)$. Indeed, each irreducible component of $Z$ is an irreducible component of a subset of the form $p_2^{-1}(Z_0)$ for some $Z_0\in \sS^{(q)}_B(p)$, from which our assertion follows.\\

(3)  Let $W=\{W_p\subset X\times\Delta^p\}$ be a simplicial closed subset, with each $W_p\in \sS_{X}^{(q)}$. Then 
$p\mapsto \Ind(W_p)$ forms a simplicial closed subset $\Ind(W)$ in $\sS_X^{(q)}$ and $p\mapsto (W_p)^+$ and $p\mapsto (W_p)^c$ form simplicial closed subsets  $W^+$ and $W^c$  in $\sS_{X_{B(y)}}^{(q)}$. This follows directly from the finiteness of $\pi_{y,X}$.
\end{rems}

\begin{lem}\label{lem:Birat}   Let $W$ be a proper closed subscheme of $X\times\Delta^p$ of pure codimension $q>0$. The map
\[
f_p:W_{B(y)}\setminus f_p^{-1}(f_p(\Ind(W)_{B(y)}))\to
f_p(W_{B(y)}) \setminus f_p(\Ind(W)_{B(y)})
\]
is birational.
\end{lem}

\begin{proof}We delete the subscripts ${}_{B(y)}$ indicating base-change unless they are needed for clarity. We write $f$ for $f_p$.

 If $Z'\subset\A^m\times\Delta^p$ is induced, then clearly
$f^{-1}(Z')$ is induced, hence, if $W'$ is an irreducible component of $W$, then
$f(W')$ is induced if and only if $W'$ is induced. In particular, if $W'$ is an
irreducible component of $W$ which is not induced, then $f(W')$ is not contained in $f(\Ind(W))$. Thus, we may assume that $\Ind(W)=\0$.  

Write $W$ as 
\[
W=W_\gen\cup \amalg_{b\in B^{(1)}}W(b),
\]
where $W_\gen$ is the union of the irreducible components of $W$ which dominate some component of $B$, and $W(b)$ is the union of irreducible components of $W$ contained in $X_b\times\Delta^p$. As $f_p$ is a $B$-morphism, it suffices to prove the result for $W=W_\gen$ and $W=W(b)$. In the first case, it suffices to replace $W$ with $W_\eta$, $\eta$ a generic point of $B$ and in the second, we may replace $X$ with $X_b$. Thus, it suffices to prove the result for $A$ a field $F$.

Let $W_0\subset X$ be the closure of the image of $W$ under the projection $p_X:X\times\Delta^p\to X$. We
proceed by induction on the maximum $d$ of the dimension of an irreducible component of $W_0$, starting with $d=-1$ (i.e. $W_0=\0$). Suppose we know our result for $d-1$. Let $W'\subset W$ be the union of the irreducible
components of $W$ whose projections to $X$ have dimension at most $d-1$, and let $W''$ be the union of the remaining components.

Let $W_1$ be an irreducible component of $W''$ and $W_2$ an irreducible component of $W'$. Suppose $f(W_1)\subset
f(W_2)$ or $f(W_2)\subset f(W_1)$. Since the projection $\pi_{y,X}$ is finite, and each component of $W$ has the
same dimension, it follows that $f(W_1)=f(W_2)$. 

Let $p_{\A^m}:\A^m\times\Delta^p\to \A^m$ be the projection. Then
\[
p_{\A^m}(f(W_1))=p_{\A^m}(f(W_2)),
\]
which is impossible, since 
\begin{align*}
&\dim_B p_{\A^m}(f(W_1))=\dim_B \pi_{y,X}(p_X(W_1))=d,\\
&\dim_B p_{\A^m}(f(W_2))=\dim_B \pi_{y,X}(p_X(W_2))\le d-1.
\end{align*}
In particular, it follows from our induction hypothesis that
 \[
f:W'\setminus  f^{-1}(f(W'')) \to f(W')\setminus f(W'')
\]
is birational. Similarly, $f:W\to  f(W)$ is birational if and only if $f:W''\to f(W'')$ is birational.  Thus, we may assume that $W=W''$.

We first consider the case $p=0$. Let $W$ be a pure codimension $q$ proper closed subset of $X$, let $W_1,\ldots, W_r$ be the irreducible components of $W$ and let $Z=\sum_iW_i$.   By \lemref{lem:MovLem1}(1),   the cycle $f^*(f_*(Z))-Z$ is effective, and has no $W_i$ as component. 

 In particular, each component of $f_*(Z)$ has multiplicity one, hence the irreducible cycles $f_*(W_i)$ each have multiplicity one, and are pair-wise distinct. Since the multiplicity of $f_*(W_i)$ is the degree of the field extension $[F(y)(W_i):F(y)(f(W_i))]$, this implies that the map $W_{B(y)}\to f(W_{B(y)})$ is  birational, as desired. 
 
The case $p=0$ easily implies the general case, in case the image subset $W_0$ is a proper subset of $X$. Thus we have only to consider the case in which each component of $W$ dominates $X$. 

Fix an  $\bar{F}$-point $z$ of $X$.  Suppose that, for all $\bar{F}$-points $x$ in an open subset of $X$, the fiber $W_x\subset x\times\Delta^p=\Delta^p_{\bar{F}}$ of $W$ over $x$ has an irreducible component $Z$ in common with the fiber $W_z$. Then clearly $W_{\bar{F}}$ contains $X\times Z$, from which it follows that $W$ contains an induced component, namely, the  union of the conjugates of $X\times Z$ over $F$. Since we have assumed that this is not the case, it follows that, for all pairs $(z,x)$ in an open subset $U$ of $X\times X$, the fibers $W_z$ and $W_x$ have no components in common. 

Let $\eta$ be a geometric generic point of $\A^m$ over $F$, such that $\eta$ remains a generic point over $F(y)$. If $z$ and $x$ are distinct points in $\pi_{y,X}^{-1}(\eta)$, then a dimension count implies that $(z,x)$ is a $\overline{F}$-generic point
of $X\times_F X$ and hence $(z,x)$ is a geometric point of $U$. Thus, $W_z$ and $W_x$ have no common components. Since 
\[
f(W)_\eta=\cup_{x\in \pi_{y,X}^{-1}(\eta)} W_x\subset \Delta^p_{\overline{F(y,\eta)}}
\]
as closed subsets of $\Delta^p_{\overline{F(y,\eta)}}$,
it follows that the map $W_{\pi_{y,X}^{-1}(\eta)}\to f(W)_\eta$ is  is generically one to one. 

As $f:W_{F(y)}\to f(W_{F(y)})$ is a map over $\A^m_{F(y)}$, this map is  one to one over a dense open subset of $f(W_{F(y)})$. Since,  by \lemref{lem:MovLem1}(2),  $W_{F(y)}$ is generically \'etale over $f(W_{F(y)})$, $W_{F(y)}\to f(W_{F(y)})$ is thus birational, completing the proof.
\end{proof}

\subsection{Behavior in codimension one} We now extend our analysis to behavior at the codimension one points of $B$. For $b\in B$, we let $b(y)=B(y)\times_Bb$, and write $f_{p,b}$ for the fiber of $f_p$ over $b(y)$. 

\begin{lem} \label{lem:Ram} Let $X$ be in $\Sm/B$ with a closed embedding $i:X\to\A^n_B$ satisfying \eqref{eqn:Finiteness}. Let $W\subset X\times\Delta^p$ be a closed subset of pure codimension $q>0$. Then
\begin{enumerate}
\item $R_y$ is a pure codimension one subscheme of $X_{B(y)}$ and is flat over $B(y)$.
\item For all $b\in B$, $R_y\times\Delta^p\cap W_{b(y)}$ has pure codimension $1$ on $W_{b(y)}$.
\end{enumerate}
\end{lem}

\begin{proof} Since $X_{B(y)}$ and $\A^m_{B(y)}$ are smooth over $B(y)$, $R_y$ is the subscheme defined by the vanishing of 
\[
\det(d\pi_{y,X}):\det\Omega_{X_{B(y)}/B(y)}\to \pi_{y,X}^*(\det\Omega_{\A^m_{B(y)}}).
\]
In particular, $R_y$ is a Cartier divisor  on $X$. From the case of fields, we know (\cite{Roberts}, \cite[Part I, Chap. II, Lemma 3.5.4]{MixMot}) that $\pi_{y_b,X_b}:X_{b(y)}\to\A^m_{b(y)}$ is generically \'etale for each point $b\in B$; since $R_{y_b}=R_y\times_{B(y)}b(y)$, it follows that $R_y$ is flat over $B$. This proves (1).

For (2), it suffices to show that $R_{y_b}\times\Delta^p\cap W_{b(y)}$ has pure codimension $1$ on $W_{b(y)}$ for all $b\in B$, which reduces us to the case of fields. This follows from \cite[{\it loc. cit.}]{MixMot}.
\end{proof}

\begin{rem} \label{rem:PushPull} Let $W\subset X\times\Delta^p$ be an irreducible closed subset. If $W$ is not induced, it follows from Lemma~\ref{lem:Birat} and Lemma~\ref{lem:Ram} that  $f_{p*}(W)=1\cdot f_p(W)$, thus
\[
f_p^*(f_{p*}(W))=1\cdot W+ Z^+,
\]
where $Z^+$ is an effective cycle with support $W^+$.
Similarly, if $W$ is induced, then $f_{p*}(W)=d_W\cdot f_p(W)$ for some integer $d_W\ge1$, and
\[
f_p^*(f_{p*}(W))= d_W\cdot W+ Z^+.
\]
where $Z^+$ as above  is an effective cycle with support $W^+$, and with $Z^+\ge d_W\cdot \cyc(W^+)$. Thus if $W$ is an arbitrary closed subset of  $X\times\Delta^p$ of pure codimension, then
\[
\cyc(W^c)-\cyc(W^+)=\sum_Z1\cdot Z,
\]
where the sum is over the induced irreducible components $Z$ of $W$ with $d_Z>1$.

\end{rem}

\begin{lem} Let $W\subset X\times\Delta^p$ be a pure codimension $q$ closed subset, and let $Z$ be an irreducible  codimension $q$ component of $W\cap W^c$.   Then $Z$ is induced.
\end{lem}

\begin{proof} This follows directly from Remark~\ref{rem:PushPull}.
\end{proof}

\begin{lem}\label{lem:base-change} Let $W\subset X\times\Delta^p$ be a pure codimension $q$ closed subset. Suppose that each irreducible component of $W$ dominates $B$. Then for each $b\in B$ we have the identity of closed subsets
\[
(W^c)_{b(y)} =(W_b)^c.
\]

\end{lem}

\begin{proof} Let $W_b^*\subset X_b\times\Delta^p$ be an effective cycle with support equal to $W_b$. It follows from Remark~\ref{rem:PushPull} (applied with $b$ replacing $B$) that $(W_b)^c$ is the support of $f_{p,c}^*(f_{p,c*}(W_b^*))-W_b^*$. For instance, we may take $W^*_b=i^*(\cyc\, W)$, where
$i:X_{b(y)}\times\Delta^p\to X_{B(y)}\times\Delta^p$ is the inclusion.

We have
\[
i^*(f_p^*(f_{p*}(\cyc\, W))-\cyc\, W)=f_{p,b}^*(f_{p,b*}(i^*(\cyc\, W))-i^*(\cyc\, W);
\]
by our remarks above,  the right-hand side has support $(W_b)^c$;  the left-hand side has support $(W^c)_{b(y)}$ by Remark~\ref{rem:PushPull}.
\end{proof}

We can now verify the fundamental properties of the projecting cone construction  for a one-dimensional base $B$.

\begin{prop}\label{prop:ChowMovLem} Let $X$ be in $\Sm/B$ with a closed embedding $i:X\to\A^n_B$ satisfying \eqref{eqn:Finiteness}, , $\sC$ a finite set of
irreducible locally closed subsets of $X$ such that each $C\in\sC$ dominates $B$, and $e:\sC\to\N$ a
function. Take $q>0$ and $W\in\sS_{X,\sC,e}^{(q)}(p)$.  Then  $W^c$  is in  $\sS_{X_{B(y)},\sC,e-1}^{(q)}(p)$. 
\end{prop}

\begin{proof}  We may assume that $W$ is irreducible. If $W\subset X_b$ for some closed point $b\in B$, the result follows from the case of a fields ({\it cf.} Lemma~\ref{lem:MovLem1}), so we may suppose that $W$ dominates $B$.  By Remark~\ref{rem:LocSuppCond}, it suffices to show
\begin{enumerate}
\item $(W^c)_{b(y)}\in \sS_{X_{b(y)},i_b^*\sC,i_b^*e-1}^{(q)}(p)$ for $b$ a generic point of $B$.
\item $(W^c)_{b(y)}\in \sS_{X_{b(y)},i_b^*\sC,i_b^*e-1}^{(q-1)}(p)$ for $b\in B^{(1)}$.
\end{enumerate}

Take $b\in B$. If $b$ is the generic point of $B$, then $W_b$ is in $\sS_{X_b,i_b^*\sC,i_b^*e}^{(q)}(p)$, so by Lemma~\ref{lem:MovLem1}, $(W_b)^c$ is in 
$\sS_{X_{b(y)},i_b^*\sC,i_b^*e-1}^{(q)}(p)$. If $b$ is a codimension one point of $B$, then by Remark~\ref{rem:LocSuppCond}, $W_b$ is in $\sS_{X_b,i_b^*\sC,i_b^*e}^{(q-1)}(p)$, so as above, $(W_b)^c$ is in $\sS_{X_{b(y)},i_b^*\sC,i_b^*e-1}^{(q-1)}(p)$. By Lemma~\ref{lem:base-change}, $(W^c)_{b(y)}=(W_b)^c$ (as closed subsets), which completes the proof.
\end{proof}

\begin{prop}\label{prop:iso} Let $X$ be in $\Sm/B$ with a closed embedding $i:X\to\A^n_B$ satisfying \eqref{eqn:Finiteness}. Take $q>0$ and $W\in \sS_X^{(q)}(p)$. Then
\begin{enumerate}
\item $f_p$ is \'etale in a neighborhood of $W_{B(y)}\setminus(\Ind(W)_{B(y)}\cup W^+)$,
\item the map
\begin{multline*}
W_{B(y)}\setminus(\Ind(W)_{B(y)}\cup W^+)\\
\xrightarrow{f_p} f(W_{B(y)})\setminus f_p((\Ind(W)_{B(y)}\cup W^+)\cap W_{B(y)})
\end{multline*}
is an isomorphism.
\end{enumerate}
\end{prop}

\begin{proof} Set
\[
W^0:=W_{B(y)}\setminus(\Ind(W)_{B(y)}\cup W^+)
\]
and
\[
W^1:=f(W_{B(y)})\setminus f_p((\Ind(W)_{B(y)}\cup W^+)\cap W_{B(y)}).
\]
Let $R_y\subset X_{B(y)}$ be the ramification locus of $\pi_{y,X}$.

We note that
\[
f_p^{-1}(f_p(\Ind(W))\subset \Ind(W)\cup W^+.
\]
Thus, as in the proof of Lemma~\ref{lem:Birat}, we may assume that $\Ind(W)=\0$. 

Take a point $x\in W$ and let $x'=f_p(x)$. The assertions are local for the \'etale topology on $\A^m\times\Delta^p$, so, after changing notation, we may replace $\A^m_{B(y)}\times\Delta^p$ with its Henselization $V$ at $x'$.  We set
\[
U:=X_{B(y)}\times\Delta^p\times_{\A^m_{B(y)}\times\Delta^p}V.
\]
Let $x_1=x$ and let $x_2,\ldots, x_s$ be the other points of $f_p^{-1}(x')$. We order the $x_i$ so that $x_2,\ldots, x_r$ are not in $R_y\times\Delta^p$ and $x_{r+1},\ldots, x_s$ are. Since $V$ is Hensel, we have a decomposition of $U$ as a disjoint union
\[
U=\amalg_{i=1}^sU_i
\]
with $x_i\in U_i$.  Let $W_i=W_{B(y)}\times_{X\times\Delta^p} U_i$, $W^+_i=W^+\times_{X\times\Delta^p} U_i$.

We claim that, if $x_i\in W_{B(y)}$ for some $i>1$, or if $x$ is in $R_y\times\Delta^p$, then $x$ is in $W^+$. We temporarily write $f_p$ for the map $f_p\times\id_V$. Suppose first that $x$ is not in $R_y\times\Delta^p$.  By Lemma~\ref{lem:Birat}, $W_{B(y)}$ maps birationally to its image $f_p(W_{B(y)})$, so the same is true after pull-back by $V$. Thus, $W_i$ maps birationally to its image $f_p(W_i)$, and $f_p(W_i)$ has no component in common with $f_p(W_1)$. Therefore $f_p^{-1}(W_i)\cap U_1$ is a union of components of $W^+_1$, so $x=x_1$ is in $W^+$.

If $x$ is in $R_y\times\Delta^p$, let $d$ be the degree of $f_p:U_1\to V$. Since $x$ is a point of ramification for $f_{p|U_1}$, $d$ is at least 2. Since $R_y$ is a divisor on $X$, $f_{p|U_1}$ is generically \'etale; since $f_p:W_1\to f_p(W_1)$ is birational, $f_p^{-1}(f_p(W_1))$ must strictly contain $W_1$. Since $f_p^{-1}(f_p(W_1))=W_1\cup W^+_1$, it follows that  $W^+_1$ is non-empty; since $U_1$ is local with closed point $x$, we must have $x\in W_1^+$.

From this, it follows that
\begin{enumera}
\item $f_p(W_{B(y)}\setminus W^+)\subset f_p(W_{B(y)})\setminus f_p(W^+\cap W_{B(y)})$,
\item $W_{B(y)}\cap R_y\times\Delta^p\subset W_{B(y)}\cap W^+$.
\item $f_p:W_{B(y)}\setminus W^+\to f_p(W_{B(y)})\setminus f_p(W^+\cap W_{B(y)})$ is set-theoretically one-to-one.
\end{enumera}
Since $f_p$ is flat, (b) implies that $f_p$ is \'etale in a neighborhood of  $W_{B(y)}\setminus W^+$, verifying (1). The statements (a), (c) and (1) together imply (2).
\end{proof}

We need one last technical result.

\begin{note} Let $\Br_y=\pi_{y,X}(R_y)$. Since $\codim_{X_{B(y)}}R_y=1$, $\Br_y$ is a divisor on $\A^m_{B(y)}$.
\end{note}

\begin{lem} Let $W$ be in $\sS^{(q)}_X(p)$. Then $f_p(W)$ is in  $\sS^{(q)}_{\A^m_{B(y)},\{\Br_y\}}(p)$.
\end{lem}

\begin{proof} Noting that each face $F$ of $\Delta^p$ is a $\Delta^s$ for some $s$, it suffices to prove that 
\[
\codim_{f_p(W)}(f_p(W)\cap \Br_y\times \Delta^p)\ge1.
\]
for each pure codimension $q$ closed subset $W$ of $X\times\Delta^p$. We may assume that $W$ is irreducible.

Let $W^0\subset X$ be the closure of the image of $W$. It suffices to show that $\pi_{y,X}(W^0_{B(y)})$ is not contained in $\Br_y$. Let $w$ be the generic point of $W^0$. Since the image of $y\times\id:B(y)\times_Bk(w)\to \sU_X\times_Bk(w)$ is Zariski dense, it follows from 
\eqref{eqn:Finiteness}(2) that $\pi_{y,X}$ is \'etale at all the points of $\pi_{y,X}^{-1}(\pi_{y,X}(w))$. Thus $\pi_{y,X}(w)$ is not in $\Br_y$, as desired.
\end{proof}

\section{Chow's moving lemma} We can now present the proof of the main moving lemma Theorem~\ref{thm:main}. Part of the proof relies on Theorem~\ref{thm:GoodCpt}; as this latter is a technical result which is completely disjoint from the main line of argument, we postpone the statement and proof of this result to \S\ref{sec:GoodCpt}.  

In this section, we fix a smooth irreducible $B$-scheme $X$ with a closed embedding $i:X\to \A^n_B$. We let $m=\dim_BX$, and we assume that $i$ satisfies \eqref{eqn:Finiteness}.  It follows from Theorem~\ref{thm:GoodCpt} that this case suffices to prove Theorem~\ref{thm:main}.

We fix a presheaf $E\in\Spt(B^{(q)})$. We extend the definition of the simplicial spectra $E^{(q)}(X,-)$ to $X\in \Sm^\ess/B$ by taking a homotopy colimit, as in \S\ref{sec:EasyMov}. We assume throughout this section that $E$ satisfies the axioms  \eqref{eqn:Ax1} and
\eqref{eqn:Ax2}.

\subsection{Excision and its consequences}

\begin{lem} \label{lem:PushPull} Let $B'$ be a $B$-scheme and $x:B'\to \sU_X$ a $B$-morphism. Let $W\supset W'$ be  simplicial closed subsets of $X\times\Delta^*$. We denote $\pi_{x,X}\times\id_{\Delta}$ by $f$.  Suppose   that for each $p\ge0$
\begin{enumerate}
\item $f_p$ is \'etale along $W_p\setminus W_p'$
\item $f_p: W_p\setminus W_p'\to f_p(W_p\setminus W_p')$ is an isomorphism (of reduced schemes).
\end{enumerate}
Let $V\subset X_{B'}\times\Delta^*$ be the simplicial closed subset $f^{-1}(f(W))$ and let $W''=f^{-1}(f(W'))$. Then the composition  (which we denote by $f^*$)
\begin{align*}
E^{f(W\setminus W')}&(\A^m_{B'}\times\Delta^*\setminus  f(W'))\\
&\xrightarrow{f^*}E^{V\setminus W''}(X_{B'}\times\Delta^*\setminus W'')\\
&=
E^{W\setminus W'}(X_{B'}\times\Delta^*\setminus W')\oplus E^{V\setminus(W\cup W'')}(X_{B'}\times\Delta^*\setminus W'')\\
&\xrightarrow{p}E^{W\setminus W'}(X_{B'}\times\Delta^*\setminus W')
\end{align*}
is an isomorphism in $\SH$. Here $p$ is the projection, and the direct sum decomposition of 
$E^{V\setminus W''}(X_{B'}\times\Delta^*\setminus W'')$ arises from the decomposition of $V\setminus W''$ into disjoint closed simplicial subsets
\[
V\setminus W''= W\setminus W'\amalg V\setminus(W\cup W'')
\]
and Lemma~\ref{lem:Add}.
\end{lem}

\begin{proof}  This follows from  axiom  \eqref{eqn:Ax2}.
\end{proof}

\begin{rem}\label{rem:UsefulSupport} For $B'=B(y)$, $x:B'\to\sU_X$ the canonical morphism $y:B(y)\to \sU_X$, $W\in\sS_X^{(q)}$ and $W'$ the simplicial  closed subset   defined by
\[
W'_p:= W_p\cap (W^+_p\cup\Ind(W)_p),
\]
the pair $(W,W')$
satisfies the conditions of Lemma~\ref{lem:PushPull}. This follows from Proposition~\ref{prop:iso}
\end{rem}

\begin{Def} Let $E$, $X$, $x:B'\to\sU_X$, $W$ and $W'$, $W''$, $V$ and $p$ be as in Lemma~\ref{lem:PushPull}. We define the push-forward morphism
\[
f_*:E^{W\setminus W'}(X_{B'}\times\Delta^*\setminus W')\to
E^{f(W\setminus W')}(\A^m_{B'}\times\Delta^*\setminus  f(W'))
\]
to be the inverse (in $\SH$) of the pull-back
\[
f^*:E^{f(W\setminus W')}(\A^m_{B'}\times\Delta^*\setminus  f(W'))\to
E^{W\setminus W'}(X_{B'}\times\Delta^*\setminus W').
\]
\end{Def}

\begin{lem}\label{lem:Replacement} Let $E$, $x:B'\to \sU_X$, $W$ and $W'$ be as in Lemma~\ref{lem:PushPull}. Then
\[
\id-f^*\circ f_*:E^{W\setminus W'}(X_{B'}\times\Delta^*\setminus W')\to
E^{W\setminus W'}(X_{B'}\times\Delta^*\setminus W')
\]
is the zero map (in $\SH$).
\end{lem}

\begin{proof} This follows directly from Lemma~\ref{lem:PushPull} and the definition of $f_*$.
\end{proof}

We will require the following notation:

\begin{Def}\label{Def:finitepush} Let $\pi:X\to Y$ be a finite morphism in $\Sm/B$, $\sC$ a finite set of irreducible locally closed subsets of $X$ and $e:\sC\to\N$ a function. We let $\sS_{Y,\sC,e}^{(q)}(p)$ be the set of $W\in\sS_Y^{(q)}(p)$ with $(\pi\times\id_{\Delta^p})^{-1}(W)\in \sS_{X,\sC,e}^{(q)}(p)$. 
\end{Def}
The $ \sS_{Y,\sC,e}^{(q)}(p)$ form a support condition on $Y$; we set 
\[
E^{(q)}(Y,-)_{\sC,e}:=E^{\sS_{Y,\sC,e}^{(q)}}(Y\times\Delta^*).
\]

\begin{rem}\label{rem:relabel} For each $C_i\in\sC$, $\pi(C_i)$ is an irreducible constructible subset of $Y$.  Setting $\sD:=\{\pi(C), C\in\sC\}$,  we define $d(\pi(C_i))$ to be the minimum of the $e(C)$ with $\pi(C)=D_i$. Then $\sS_{Y,\sC,e}^{(q)}=\sS_{Y,\sD,d}^{(q)}$. 
\end{rem}

\subsection{The proof of Theorem~\ref{thm:main}} Let $X\to\Spec B$ and $i:X\to\A^n_B$ be as above.  Let $\sC$ be a finite set of irreducible locally closed subsets of $X\in\Sm/B$ such that each $C\in\sC$ dominates $B$ and  $e:\sC\to\N$ a function. In this section, we assume that $B=\Spec A$, with $A$ a semi-local PID.  We will show that
\begin{equation}\label{eqn:reduction1}
E^{(q)}(X,-)_{\sC,e}\to E^{(q)}(X,-)
\end{equation}
is a weak equivalence. This and Theorem~\ref{thm:GoodCpt}  implies Theorem~\ref{thm:main}. 

As above, we have the $B$-morphism $y:B(y)\to\sU_X$, and the morphism of cosimplicial $B(y)$-schemes $f:X_{B(y)}\times\Delta^*\to \A^m_{B(y)}\times\Delta^*$.

If $e(C)\ge\dim_BX$ for all $C$, then $\sS_{X,\sC,e}^{(q)}=\sS_X^{(q)}$, so by descending induction on $e$ it suffices to show
\begin{equation}\label{eqn:reduction2}
E^{(q)}(X,-)_{\sC,e-1}\to E^{(q)}(X,-)_{\sC,e}
\end{equation}
is a weak equivalence. Let $E^{(q)}(X,-)_{\sC,e/e-1}$ denote the homotopy cofiber of $E^{(q)}(X,-)_{\sC,e-1}\to E^{(q)}(X,-)_{\sC,e}$; it thus suffices to show 
\begin{equation}\label{eqn:reduction3}
E^{(q)}(X,-)_{\sC,e/e-1}\text{ is weakly contractible.}
\end{equation}

\begin{lem}  The base-change
\[
\pi^*:E^{(q)}(X,-)_{\sC,e/e-1}\to E^{(q)}(X_{B(y)},-)_{\pi^*\sC,e/e-1}
\]
induces the zero-map on homotopy groups.
\end{lem}
\begin{proof} By axiom  \eqref{eqn:Ax1}, we have the weak equivalence of $E^{(q)}(X,-)_{\sC,e/e-1}$ with the homotopy colimit
\begin{equation}\label{eqn:hocolim}
E^{(q)}(X,-)_{\sC,e/e-1}\sim \hocolim_{W_1, W_2}E^{W_1\setminus W_1\cap W_2}(X\times\Delta^*\setminus W_1\cap W_2)
\end{equation}
where the limit is over simplicial closed subsets $W_1\in \sS_{X,\sC,e}^{(q)}$, $W_2\in \sS_{X,\sC,e-1}^{(q)}$. We call  a such pair $(W_1, W_2)$ {\em allowable} for $X,\sC,e$. Similarly, we call a pair $(W_1, W_2)$ of closed simplicial subsets of $\A^m\times\Delta^*$ allowable for $\A^m,\sC,e$ if $W_1$ is in $\sS_{\A^m,\sC,e}^{(q)}$ and $W_2$ is in $\sS_{\A^m,\sC,e-1}^{(q)}$ (see Definition~\ref{Def:finitepush} for the notation) .

Let $(W_1, W_2)$ be an allowable pair for $X,\sC,e$. Set 
\[
E^{(W_1,W_2)}(X,-):=E^{W_1\setminus W_1\cap W_2}(X\times\Delta^*\setminus W_1\cap W_2). 
\]
We set 
\[
\pi^*(W_1, W_2):=(W_{1B(y)}, W_{2B(y)}\cup W_1^+\cup \Ind(W_1)_{B(y)}),
\]
and
\begin{multline*}
f(\pi^*(W_1, W_2))\\
:=(f(W_{1B(y)}), f(W_{2B(y)}\cup [(W_1^+\cup \Ind(W_1)_{B(y)})\cap W_{1B(y)}])).
\end{multline*}
By Proposition~\ref{prop:ChowMovLem}, Remark~\ref{rem:PushPull} and Remark~\ref{rem:Induced}, $\pi^*(W_1, W_2)$ is allowable for $X_{B(y)}, \pi^*\sC,e$ and $f(\pi^*(W_1, W_2))$ is allowable for $\A^m_{B(y)},\sC,e$.

To complete the proof, it suffices to show that for each allowable pair $(W_1,W_2)$, the composition
\[
E^{(W_1,W_2)}(X,-)\xrightarrow{\iota_{W1,W_2}} E^{(q)}(X,-)_{\sC,e/e-1}\xrightarrow{\pi^*} E^{(q)}(X_{B(y)},-)_{\pi^*\sC,e/e-1}
\]
induces zero on the homotopy groups, where $\iota_{W1,W_2}$ is the map defined  by the description of
$E^{(q)}(X,-)_{\sC,e/e-1}$ in \eqref{eqn:hocolim}.

By Lemma~\ref{rem:UsefulSupport}, Remark~\ref{rem:UsefulSupport} and Lemma~\ref{lem:Replacement}, we may enlarge $W_2$ so that $f(W_1\setminus W_2)= f(W_1)\setminus f(W_2)$, $f$ is \'etale on a neighborhood of $W_1\setminus W_2$ and 
$f:W_1\setminus W_2\to f(W_1)\setminus f(W_2)$ is an isomorphism of reduced schemes.

Thus, using Lemma~\ref{lem:PushPull}, the   pull-back
\[
f^*:E^{f(\pi^*(W_1,W_2))}(\A^m_{B(y)},-)\to E^{\pi^*(W_1,W_2)}(X_{B(y)},-)
\]
admits the inverse
\[
f_*:E^{\pi^*(W_1,W_2)}(X_{B(y)},-)\to
E^{f(\pi^*(W_1,W_2))}(\A^m_{B(y)},-)
\]
in $\SH$.

 We consider the composition
\begin{multline*}
E^{(W_1,W_2)}(X,-)\xrightarrow{\pi^*} E^{(W_{1B(y)},W_{2B(y)})}(X_{B(y)},-)\\
\to E^{\pi^*(W_1,W_2)}(X_{B(y)},-)\xrightarrow{f_*}E^{f(\pi^*(W_1,W_2))}(\A^m_{B(y)},-)\\
\xrightarrow{f^*}E^{\pi^*(W_1,W_2)}(X_{B(y)},-)\xrightarrow{\iota}
E^{(q)}(X_{B(y)},-)_{\pi^*\sC,e/e-1};
\end{multline*}
by Lemma~\ref{lem:Replacement}, this composition induces the same map on homotopy groups as $\pi^*\circ\iota_{W_1,W_2}$. Thus it suffices to show that $\iota\circ f^*$ induces zero on homotopy groups.

We can factor $\iota\circ f^*$  as
\begin{multline*}
E^{f(\pi^*(W_1,W_2))}(\A^m_{B(y)},-)\xrightarrow{\iota'}
E^{(q)}(\A^m_{B(y)},-)_{\pi^*\sC,e/e-1}\\
\xrightarrow{f^*} E^{(q)}(X_{B(y)},-)_{\pi^*\sC,e/e-1}.
\end{multline*}
By Proposition~\ref{TransMovProp} and Remark~\ref{rem:relabel}  the middle term is weakly contractible, which completes the proof. 
\end{proof}

The argument we used in the proof of Proposition~\ref{TransMovProp} shows that the base-change $\pi^*$ is injective on homotopy groups. This verifies \eqref{eqn:reduction3}, completing the proof of
Theorem~\ref{thm:main}.

\section{Funtoriality} We complete the program by showing how Theorem~\ref{thm:main} yields functoriality for the presheaves $E^{(q)}(X_\Nis,-)$ and $E^{(q)}(X_\Zar,-)$. We actually prove two different levels of functoriality, a weak form, being the extension to a functor
\[
E^{(q)}:\Sm/B^\op\to \SH,
\]
and a strong form, which gives a functor
\[
E^{(q)}:\Sm/B^\op\to \Spt.
\]
For the strong form, we require that $B=\Spec F$, $F$ a field.

Let $E$ be in $\Spt(B^{(q)})$. Recall from \S\ref{subsec:FirstProperties} the subcategory $\Sm\ds B$ of $\Sm/B$ with the same objects, but with morphisms being the smooth morphisms in $\Sm/B$, the
presheaf
\[
E^{(q)}_\sm:\Sm\ds B^\op\to\Spt
\]
and the homotopy coniveau tower \eqref{eqn:HCTowerSm}. This is the starting point for full functoriality.

Throughout this section, we assume that $B$ is a regular noetherian separated scheme of Krull dimension one.

\subsection{Functoriality in $\SH$}  In order to give a natural setting for the form of functoriality one achieves over a general base $B$, we introduce some notation.

Let $f:Y\to X$ be a morphism in $\Sm/B$, $F$ in $\Spt_\Nis(Y)$. We have the functor $f_*:\Spt_\Nis(Y)\to\Spt_\Nis(X)$ with left adjoint $f^*$. Checking stalkwise, one sees that $f^*$ preserves cofibrations and weak equivalences, hence $f_*$ preserves fibrations.

Since $\sH\Spt_\Nis(Y)$ is equivalent to the homotopy category of fibrant objects in $\Spt_\Nis(Y)$, we have the functor 
\[
f_*:\sH\Spt_\Nis(Y)\to \sH\Spt_\Nis(X)
\]
sending a fibrant object $F\in \Spt_\Nis(Y)$ to $f_*F\in \sH\Spt_\Nis(X)$. 

For morphisms $f:Y\to X$, $g:Z\to Y$, there is a canonical isomorphism $\theta_{f,g}:f_*g_*\to (fg)_*$, satisfying the standard cocycle condition. Indeed, if $F$ is fibrant in $\Spt_\Nis(Z)$, then $f_*F$ is fibrant in $\Spt_\Nis(Y)$, so the canonical isomorphism $g_*f_*F\to (gf)_*F$ in $\Spt_\Nis(X)$  defines the natural isomorphism $\theta_{f,g}$. The cocycle condition is checked in a similar fashion.

Form the category $\sH\Spt_\Nis/B$ with objects $(X,\sE)$, with $X\in\Sm/B$ and with $\sE\in \sH\Spt_\Nis(X)$. A morphism $\alpha:(X,\sE)\to (Y,\sF)$ is a morphism $f:Y\to X$ in $\Sm/B$ together with a morphism $\phi:\sE\to f_*\sF$ in $\sH\Spt_\Nis(X)$. Composition is given by
\[
(g,\psi)\circ(f,\phi):=(fg,\theta_{g,f}\circ g_*(\psi)\circ \phi).
\]
Sending $(X,\sE)$ to $X$ defines the functor $\pi:\sH\Spt_\Nis/B\to \Sm/B^\op$. 

The functor  $E^{(q)}_\sm$ defines a section to $\pi$ over the subcategory $\Sm\ds B^\op$:
\[
\tilde{E}^{(q)}_\sm:\Sm\ds B^\op\to \sH\Spt_\Nis/B,
\]
with $\tilde{E}^{(q)}_\sm(X)$ the presheaf $E^{(q)}(X_\Nis,-)$.

With these notations, we can state our first main result on functoriality:
\begin{thm}\label{thm:Funct1} Let $E$ be in $\Spt(B^{(q)})$, satisfying axioms  \eqref{eqn:Ax1},
\eqref{eqn:Ax2} and \eqref{eqn:Ax3}. Then the section $\tilde{E}^{(q)}_\sm$ extends (up to isomorphism) to a section
\[
\tilde{E}^{(q)}:\Sm/B^\op\to \sH\Spt_\Nis/B.
\] 
\end{thm}

The remainder of this section is devoted to the proof of this result. As a matter of notation, if we have composable morphisms $f:Y\to X$, $g:Z\to Y$, objects $E_X\in\Spt_\Nis(X)$, $E_Y\in\Spt_\Nis(Y)$ and $E_Z\in\Spt_\Nis(Z)$ and morphisms $\phi:E_X\to f_*E_Y$, $\psi:E_Y\to g_*E_Z$, we let $\psi\circ\phi:E_X\to (fg)_*E_Z$ denote the composition $\theta_{f,g}\circ g_*(\psi)\circ\phi$.

We consider the following situation: Suppose $B$ is irreducible. Let $f:Y\to X$ be a morphism in $\Sm/B$ with $X$ and $Y$ equi-dimensional over $B$. Let $X^f\subset X$ be the subset of $X$ of points over which $f$ is equi-dimensional with fiber dimension $\dim_BY-\dim_BX$. Let $\sC$ be a finite set of irreducible locally closed subsets of $X$ such that there is a subset $\sC_0\subset \sC$ satisfying
\begin{enumerate}
\item $X\setminus X^f\subset \cup_{C\in \sC_0}C$
\item For each $C\in\sC_0$, the restriction of $f$ to $f:f^{-1}(C)\to C$ is equi-dimensional.
\item Each $C\in \sC_0$ dominates a component of $B$.
\end{enumerate}
We call such a $\sC$ {\em good for $f$}. We extend this notion to the general case by taking disjoint unions.

Suppose that $\sC$ is good for $f$.  Then
\[
W\in \sS^{(q)}_{X,\sC}(p)\Longrightarrow (f\times\id)^{-1}(W)\in \sS^{(q)}_Y(p),
\]
so we have the well-defined map of simplicial spectra
\[
f^*:E^{(q)}(X,-)_\sC\to E^{(q)}(Y,-).
\]
This construction is compatible with pull-back by a smooth morphism $Z\to X$, so we have the map of presheaves on $X_\Nis$:
\[
f^*:E^{(q)}(X_\Nis,-)_\sC\to f_*E^{(q)}(Y_\Nis,-).
\]
If we vary $q$, we have the well-defined map of towers of presheaves
\[
f^*:E^{(-)}(X_\Nis,-)_\sC\to f_*E^{(-)}(Y_\Nis,-).
\]

By Theorem~\ref{thm:main}(1), the  map $E^{(q)}(X_\Nis,-)_\sC\to E^{(q)}(X_\Nis,-)$ is a local weak equivalence in $\Spt_\Nis(X)$, so the diagram
\[
E^{(q)}(X_\Nis,-)\leftarrow E^{(q)}(X_\Nis,-)_\sC\xrightarrow{f^*}
f_*E^{(-)}(Y_\Nis,-)
\]
defines the map in the homotopy category $\sH\Spt_\Nis(X)$
\begin{equation}\label{eqn:pulback1}
f^*:E^{(q)}(X_\Nis,-)\to f_*E^{(q)}(Y_\Nis,-).
\end{equation}

\begin{lem} The map \eqref{eqn:pulback1} is independent of the choice of good $\sC$.
\end{lem}

\begin{proof}  Let $\sC'$ be another collection which is good for $f$. Then $\sC\cup\sC'$ is clearly good for $f$, so we may assume that $\sC'\supset \sC$. This yields the commutative diagram of presheaves
\[
\xymatrix{
E^{(q)}(X_\Nis,-)&E^{(q)}(X_\Nis,-)_\sC\ar[l]\ar[d]^{f^*}\\
E^{(q)}(X_\Nis,-)_{\sC'}\ar[u]\ar[ur]\ar[r]_{f^*}&f_*E^{(-)}(Y_\Nis,-),
}
\]
whence the result.
\end{proof}

For a general $f$, there may be no $\sC$ which is good for $f$. However, if $f$ is smooth, then every $\sC$    is good for $f$. If $i:Y\to X$ is a closed embedding, and $\sC$ is any collection containing $i(Y)$, then $\sC$ is good for $i$. Thus $f^*$ is defined for smooth morphisms and for closed embeddings. This justifies the following 

\begin{Def} Let $f:Y\to X$ be a morphism in $\Sm/B$. Factor $f$ as $f=p\circ i$, where $p:Y\times_BX\to X$ is the projection, and $i:Y\to Y\times_BX$ is the graph, $i=(\id_Y,f)$. Define the map $f^*$ in $\sH\Spt_\Nis(X)$ by 
\[
f^*:=i^*\circ p^*:E^{(q)}(X_\Nis,-)\to f_*E^{(q)}(Y_\Nis,-).
\]
\end{Def}

We proceed to show the functoriality $(fg)^*=g^*f^*$.

\begin{lem}\label{lem:basechng} Let 
\[\xymatrix{
T\ar[r]^{i_2}\ar[d]_{f_2}&Y\ar[d]^{f_1}\\
Z\ar[r]_{i_1}&X}
\]
be a commutative diagram in $\Sm/B$, with $f_1, f_2$ smooth morphisms, and $i_1, i_2$ closed embeddings. Then  $i_2^*\circ f_1^*=f_2^*\circ i_1^*$ in $\sH\Spt_\Nis(X)$.
\end{lem}

\begin{proof} Let $g=f_1i_2=f_2i_1$. Let 
\[
\sC_X:= \{i_1(Z)\},\ \sC_Y:=\{i_2(T)\}.
\]
One  checks that
\[
W\in \sS^{(q)}_{X,\sC_X}(p)\Longrightarrow (f_1\times\id)^{-1}(W)\in \sS^{(q)}_{Y,\sC_Y}(p).
\]

Thus, we have the pull-back 
\[
f_1^*:E^{(q)}(X_\Nis,-)_{\sC_X}\to f_{1*}E^{(q)}(Y_\Nis,-)_{\sC_Y},
\]
giving the commutative diagram
\[
\xymatrix{
&f_{1*}E^{(q)}(Y_\Nis,-)\\
E^{(q)}(X_\Nis,-)_{\sC_X}\ar[r]_{f_1^*}\ar[ru]^{f_1^*}\ar[d]_{i_1^*}& f_{1*}E^{(q)}(Y_\Nis,-)_{\sC_Y}
\ar[d]^{f_{1*}(i_2^*)}\ar[u]\\
 i_{1*}E^{(q)}(Z_\Nis,-)\ar[r]_{i_{1*}(f_2^*)}&g_*E^{(q)}(T_\Nis,-)
}
\]
The result follows easily from this.
\end{proof}

\begin{lem}\label{lem:funct2} Let $f:Y\to X$, $g:Z\to Y$ be morphisms in $\Sm/B$. If $f$ and $g$ are both smooth morphisms, or both closed embeddings, then 
\[
(fg)^*=g^*f^*:E^{(q)}(X_\Nis,-)\to (fg)_*E^{(q)}(Z_\Nis,-).
\]
in $\sH\Spt_\Nis(X)$.
\end{lem}

\begin{proof} For smooth morphisms $f$ and $g$,  $(fg)^*=g^*f^*$ in $\Spt_\Nis(X)$. If $f$ and $g$ are closed embeddings, take $\sC_X=\{f(Y), f(g(Z))\}$ and $\sC_Y=\{g(Z)\}$. Then
\[
W\in\sS_{X,\sC_X}^{(q)}(p)\Longrightarrow
(f\times\id)^{-1}(W)\in \sS_{Y,\sC_Y}(p),
\]
 giving us the pull-back map
\[
f^*:E^{(q)}(X_\Nis,-)_{\sC_X}\to f_*E^{(q)}(Y_\Nis,-)_{\sC_Y}.
\]

In addition, $\sC_X$ is good for $f$ and for $fg$ and we have the commutative diagram
\[
\xymatrix{
&f_*E^{(q)}(Y_\Nis,-)\\
E^{(q)}(X_\Nis,-)_{\sC_X}\ar[r]_{f^*}\ar[ur]^{f^*}\ar[dr]_{(fg)^*}&f_*E^{(q)}(Y_\Nis,-)_{\sC_Y}\ar[u]
\ar[d]^{f_*(g^*)}\\
&(fg)_*E^{(q)}(Z_\Nis,-).\\
}
\]
\end{proof}

\begin{prop}\label{prop:Funct1} Let $f:Y\to X$, $g:Z\to Y$ be morphisms in $\Sm/B$. Then 
\[
(fg)^*=g^*f^*:E^{(q)}(X_\Nis,-)\to (fg)_*E^{(q)}(Z_\Nis,-).
\]
in $\sH\Spt_\Nis(X)$.
\end{prop}

\begin{proof} Consider the commutative diagram
\[
\xymatrix{
Z\ar[rr]^-{(\id_Z,g)}&&Y\times_BZ\ar[d]_{p_1}\ar[rr]^-{(\id_Y,f)\times\id_Z}&&X\times_BY\times_BZ
\ar[d]^{p_{12}}\\
&&Y\ar[rr]_{(\id_Y,f)}&&X\times_BY\ar[d]^{p_1}\\
&&&&X}
\]
Let $i=[(\id_Y,f)\times\id_Z]\circ(\id_Z,g)$ and $p=p_1\circ p_{12}$. Applying Lemma~\ref{lem:basechng} and Lemma~\ref{lem:funct2} to this diagram shows that
\[
g^*\circ f^*=i^*\circ p^*.
\]

Applying Lemma~\ref{lem:basechng} and Lemma~\ref{lem:funct2} to the commutative diagram
\[
\xymatrix{
Z\ar[r]^-i\ar@{=}[d]&X\times_BY\times_BZ\ar[d]_{p_{13}}\ar[dr]^p\\
Z\ar[r]_-{(\id_Z,g)}&X\times_BZ\ar[r]_{p_1}&X}
\]
yields
\[
i^*\circ p^*=(fg)^*.
\]

\end{proof}

To construct the section  $\tilde{E}^{(q)}$, send $X$ to $E^{(q)}(X_\Nis,-)\in\sH\Spt_\Nis(X)$ and $f:Y\to X$ to $(f,f^*):(X, E^{(q)}(X_\Nis,-))\to (Y, E^{(q)}(Y_\Nis,-))$. By Proposition~\ref{prop:Funct1}, this gives a well-defined section to $\pi$, extending $\tilde{E}^{(q)}_\sm$. Since the construciton of $f^*$ is natural with respect to change in $q$, we have the tower of functors asserted in Theorem~\ref{thm:Funct1}. This completes the proof of Theorem~\ref{thm:Funct1}.

Sending $E\in\sH\Spt_\Nis(X)$ to the global sections of a fibrant model defines the functor
\[
R\Gamma(X,-):\sH\Spt_\Nis(X)\to \SH;
\]
varying $X$, we have the functor
\[
R\Gamma:\sH\Spt_\Nis/B\to \SH.
\]
We have for example the functor 
\[
R\Gamma E^{(q)}_\sm:\Sm\ds B^\op\to \SH
\]
$R\Gamma E^{(q)}_\sm:=R\Gamma\circ\tilde{E}^{(q)}_\sm$.

\begin{cor} \label{cor:Funct1}Let $E$ be in $\Spt(B)$, satisfying   \eqref{eqn:Ax1},
\eqref{eqn:Ax2} and \eqref{eqn:Ax3}. The functor $R\Gamma E^{(q)}_\sm$ extends to the functor
\[
R\Gamma E^{(q)}:\Sm/B^\op\to \SH
\]
by setting $R\Gamma E^{(q)}:=R\Gamma\circ\tilde{E}^{(q)}$, and we have the tower of functors
\[
\ldots\to R\Gamma E^{(q+1)}\to R\Gamma E^{(q)}\to\ldots\to R\Gamma E^{(0)}.
\]
\end{cor}
This result suffices to make the sheafified version of the coniveau spectral sequence functorial; see the next section \S\ref{sec:HomConSS} for details.

\subsection{The homotopy coniveau spectral sequence}\label{sec:HomConSS}\ \\

Fix $E\in \Spt(B)^*$. Corollary~\ref{cor:Funct1} gives us the tower of functors on $\Sm/B$
\begin{equation}\label{eqn:HomConFibr}
\ldots\to R\Gamma E^{(p+1)}\to R\Gamma E^{(p)}\to\ldots\to R\Gamma E^{(0)}.
\end{equation}
from $\Sm/B^\op$ to $\SH$. We let $R\Gamma E^{(p/p+1)}$ denote the functor
\[
X\mapsto \cofib[R\Gamma E^{(p+1)}(X)\to R\Gamma E^{(p)}(X)]
\]
and set
\begin{align*}
&\H^q(X,E^{(p)}):=\pi_{-q}R\Gamma E^{(p)}(X)\\
&\H^q(X,E^{(p/p+1)}):=\pi_{-q}R\Gamma E^{(p/p+1)}(X).
\end{align*}

Evaluating the tower \eqref{eqn:HomConFibr} at $X$ and taking the spectral sequence of a tower gives us the spectral sequence 
\begin{equation}\label{eqn:HomConHypSS}
E_1^{p,q}:=\H^{p+q}(X,E^{(p/p+1)})\Longrightarrow \H^{p+q}(X,E^{(0)}).
\end{equation}
which is functorial on $\Sm/B$. In general, this sequence has no good convergence properties.

Similarly, set 
\[
E^{(p/p+1)}(X,-):=\cofib[E^{(p+1)}(X,-)\to E^{(p)}(X,-)]. 
\]
Applying the same spectral sequence construction to the homotopy coniveau tower \eqref{eqn:HCTowerSm} yields the  spectral sequence 
\begin{equation}\label{eqn:HomConSS}
E_1^{p,q}:=\pi_{-p-q}(E^{(p/p+1)}(X,-))\Longrightarrow \pi_{-p-q}(E^{(0)}(X,-)).
\end{equation}
functorial on $\Sm\ds B$. 

The canonical maps $E^{(p)}(X,-)\to R\Gamma E^{(p)}(X)$ give a map of spectral sequences
\eqref{eqn:HomConHypSS}$\to$\eqref{eqn:HomConSS}, natural in $X$ over $\Sm\ds B$. 

\begin{rem}\label{rem:HConSSFunct} The main results of \cite{Loc} can be applied to $E\in\Spt(B)^*$. One result that can be proved is that, for $X\in\Sm/B$, the presheaf $E^{(p)}(X_\Nis,-)$ satisfies a Mayer-Vietoris property, namely, for $U\subset X$ an open subscheme $V\to X$ \'etale such that $\{U\to X, V\to X\}$ is a Nisnevic cover, the square of restriction maps
\[
\xymatrix{
E^{(p)}(X,-)\ar[r]\ar[d]&E^{(p)}(U,-)\ar[d]\\
E^{(p)}(V,-)\ar[r]&E^{(p)}(U\times_X V,-)
}
\]
is homotopy cartesian. By a basic result of Thomason \cite{Thomason, Jardine1}, this property implies that the canonical map $E^{(p)}(X,-)\to R\Gamma E^{(p)}(X)$ is a weak equivalence.

In particular, the map of spectral sequences \eqref{eqn:HomConHypSS}$\to$\eqref{eqn:HomConSS} is an isomorphism of spectral sequences. We have therefore extended the functoriality of the spectral sequence \eqref{eqn:HomConSS} from $\Sm\ds B$ to $\Sm/B$.

Finally, if $E$ satisfies the additional property of homotopy invariance:
\[
p_1^*:E(X)\to E(X\times\A^1)
\]
is a weak equivalence for all $X\in\Sm/B$, then the canonical map
\[
E(X)\to E^{(0)}(X,-)
\]
is well-known to be a weak equivalence. Thus, the spectral sequence \eqref{eqn:HomConSS}  becomes
\begin{equation}\label{eqn:HomConSS2}
E_1^{p,q}:=\pi_{-p-q}(E^{(p/p+1)}(X,-))\Longrightarrow \pi_{-p-q}(E(X)).
\end{equation}

The proof of the localization theorem for $E\in\Spt(B)^*$ is sketched in our preprint \cite{HomotopyConiveau}; details will appear in a forthcoming paper.
\end{rem}

\subsection{Straightening a homotopy limit}\label{Sec:Straight} In order to pass from functoriality in the homotopy category to a strict functoriality, we first take a certain homotopy limit. This yields a functor ``up to homotopy and all higher homotopies"; we then apply a construction of Dwyer-Kan to rectify this to a strict functor. In this  section, we present these constructions in the categorical setting.

As the first ingredient in our construction, we recall the notations and results of Dwyer-Kan \cite{DwyerKan} on which the construction is based.

Recall that a {\em simplicial category} is a category enriched in simplicial sets. We make the category of simplicial sets $\sSets$ a simplicial category as usual by setting $\Hom_\sSets(A,B)_n:=\Hom_\sSets(A\times\Delta[n],B)$.  For simplicial categories $\sA$, $\sB$, (with $\sA$ small) we have the category of simplicial functors $\sA\to \sB$, denoted $\sB^\sA$.  Given $X, Y\in \sSets^\sA$, a morphism $\theta:X\to Y$ is called a weak equivalence if 
$\theta(a):X(a)\to Y(a)$ is a weak equivalence of simplicial sets for each $a\in \sA$. 

A simplicial functor $f:\sA\to \sB$ is a {\em weak equivalence} if 
\begin{enumerate}
\item For each pair of objects of $\sA$, the map of simplicial sets $f:\Hom_\sA(x,y)\to\Hom_\sB(f(x), f(y))$ is a weak equivalence.
\item  $\pi_0(f):\pi_0\sA\to \pi_0\sB$ is surjective.
\end{enumerate}

Let $f:\sA\to\sB$ be a simplicial functor. We have the pull-back functor $f^*:\sSets^\sB\to\sSets^\sA$, i.e., $f^*X(a):=X(f(a))$. Dwyer-Kan define the {\em homotopy push-forward} $\hocolim  f_*:\sSets^\sA\to\sSets^\sB$ as follows: Take $X\in \sSets^\sA$. For $b\in\sB$, $\hocolim  f_*X(b)$ is the diagonal simplicial set associated to the bisimplicial set $(\hocolim  f_*X)_{**}(b)$ with
\begin{align}\label{multline:HocolimDef}
(&\hocolim  f_*X)_{n,*}(b):=\\
&\amalg X(a_0)\times \Hom_\sA(a_0,a_1)\times\ldots\times\Hom_\sA(a_{n-1},a_n)\times \Hom_\sB(f(a_n),b)\notag
\end{align}
and the evident face and degeneracy maps, and the evident functoriality in $b$. For $a\in\sA$, we have the evident map
\[
i_a:X(a)\to \hocolim  f_*X(f(a)).
\]
This map is natural in $X$, but not in $a$.

The functor $\hocolim  f_*$ admits a right adjoint
\[
 \hocolim  f^*:\sB^\sSets\to \sA^\sSets
\]
and there is a natural transformation $f^*\to  \hocolim  f^*$. For $Y\in \sB^\sSets$, we have the adjunction $\hocolim f_*(\hocolim  f^*Y)\to Y$; we let 
\[
r_Y:\hocolim  f_* (f^*Y)\to Y
\]
 denote the composition
\[
\hocolim  f_* (f^*Y)\to \hocolim f_*(\hocolim  f^*Y)\to Y.
\]

The main result of \cite{DwyerKan} is
\begin{thm}[\hbox{\cite[Theorem 2.1]{DwyerKan}}]\label{thm:DwyerKan}  Suppose that $f:\sA\to\sB$ is a weak equivalence of small simplicial categories. Then for each $X\in\sSets^\sA$ and $a\in\sA$, the map
\[
i_a:X(a)\to \hocolim  f_*X(f(a))
\]
is a weak equivalence. Also, for each $Y\in \sSets^\sB$, the map 
\[
r_Y:\hocolim  f_* f^*Y\to Y
\]
is a weak equivalence, natural in $Y$, 
\end{thm}

Let $I$ be a 2-category. We form the simplicial category $\sN_*I$ by setting
\[
\Hom_{\sN_*I}(a,b)_*:=\sN_*\Hom_I(a,b)
\]
for each pair of objects $a,b\in I$. A simplicial functor $X:\sN_*I\to \sSets$ is called a {\em continuous} functor $X:I\to\sSets$.

Let $B$ be a small category. Form the 2-category $\Fac B$ as follows: $\Fac B$ has the same objects as $B$. A morphism $b\to b'$ in $\Fac B$ consists of a sequence of composable morphisms in $B$,
\[
b=b_0\xrightarrow{f_1}b_1\xrightarrow{f_2}\ldots\xrightarrow{f_n}b_n=b',
\]
where if some $f_i=\id$, we identify the sequence with the one obtained by deleting $f_i$ (unless $n=1$). Composition is just by concatenation. We write such a sequence as $f_n*\ldots*f_1$. Given such a sequence, we can form a new (shorter) sequence by changing some of the $*$'s to compositions in the category $B$; this operation defines a 2-morphism. We have the functor $\pi:\Fac B\to B$ which is the identity on objects, changes the concatenation $*$ to composition in $B$, and sends all 2-morphisms to the identity. This yields the functor of simplicial categories $\pi:\sN_*\Fac B\to B$ (with $B$ discrete).

The category $\Hom_{\Fac B}(b,b')$ is a partially ordered set with maximal elements $\Hom_B(b,b')$, and each $\tilde{f}\in \Hom_{\Fac B}(b,b')$ is dominated by a unique maximal element, namely $\pi(\tilde{f})$. Thus, we have

\begin{lem} $\pi:\sN_*\Fac B\to B$ is a weak equivalence of simplicial categories.
\end{lem}

We may therefore apply the Dwyer-Kan theorem: For a continuous functor $X:\Fac B\to\sSets$, we denote $\hocolim\pi_*X:B\to\sSets$ by $C_BX$, or just $CX$ if the context makes the meaning clear. We  have for each $b\in B$ the weak equivalence
\[
i_b: X(b)\to  C_BX(b).
\]
Thus, the Dwyer-Kan $\hocolim$-construction transforms a "functor up to homotopy and all higher homotopies" on $B$, i.e., a continuous functor on $\Fac B$, to an honest functor on $B$.

The second ingredient in our construction involves some aspects of the $\holim$ construction, which we now recall (for details, see \cite{BousfieldKan}). 

Let $\sA$ be a small category. For $a\in\sA$, we have the category of objects over $a$, $\sA/a$; sending $a$ to the nerve $\sN(\sA/a)$ defines the functor
\[
\sN(\sA/-):\sA\to \sSets.
\]
Given a functor $F:\sA\to \sSets$, we have the homotopy limit $\holim_\sA F$, defined by
 \[
\holim_\sA F:=\sHom(\sN(\sA/-),F), 
\]
i.e., $\holim_\sA F$ is the simplicial set with $n$-simplicies the set of natural transformations $\sigma:\sN(\sA/-)\times\Delta[n]\to F$.

From the definition, it is clear that  $\holim_\sA F$ is natural in $F$, i.e., if $\theta:F\to G$ is a natural transformation of functors, we have the map of  spaces $\theta_*:\holim_\sA F\to \holim_\sA G$. Similarly, if $i:\sB\to \sA$ is a functor of small categories, we have the map of   spaces $i^*: \holim_\sA F\to \holim_\sB(F\circ i)$. 

\begin{rem}\label{rem:HolimTriv} Suppose $F:\sA\to\sSets$ has the properties
\begin{enumerate}
\item $F(a)$ is fibrant for all $a\in\sA$.
\item $F(a)\to F(a')$ is a weak equivalence for all morphisms $a\to a'$ in $\sA$.
\item $\sA$ has a final object $a_0$.
\end{enumerate}
Then the canonical map $\holim_\sA F\to F(a_0)$ is a weak equivalence. Indeed, by \cite[XI(5.6)]{BousfieldKan} the natural transformation of $F$ to the constant functor $c_{F(a_0}$ with value $F(a_0)$ induces a weak equivalence
\[
\holim_\sA F\to \holim_\sA c_{F(a_0)}.
\]
By \cite[IX(4.2)]{BousfieldKan} 
\[
 \holim_\sA c_{F(a_0)}=\sHom(\sN(\sA),F(a_0)).
 \]
 But since $a_0$ is the final object of $\sA$, the inclusion $a_0\to \sA$ induces a homotopy equivalence $pt=\sN(a_0)\to\sN(\sA)$, hence a weak equivalence
 \[
 \sHom(\sN(\sA),F(a_0))\to F(a_0).
 \]
 \end{rem}

For a category $\sA$, we have the 2-category $\Cat/\sA$, with objects functors $\alpha:\mathfrak{a}\to \sA$, with $\mathfrak{a}$ small. A morphism $[\alpha:\mathfrak{a}\to \sA]\to [\alpha':\mathfrak{a}'\to \sA]$ is a pair $(f,\theta)$ with $f:\mathfrak{a}\to\mathfrak{a}'$ a functor, and $\theta:\alpha'\circ f\to \alpha$ a natural transformation. The composition is 
\[
(f',\theta')\circ(f,\theta):=(f'f,\theta\circ (\theta'\circ f)).
\]
A 2-morphism $\vartheta:(f,\theta)\to(f',\theta')$ is a natural transformation $\vartheta:f\to f'$ such that $\theta'\circ (\alpha'\circ\vartheta)=\theta$.

Let $F:\sA\to\sSets$ be a functor. For $\alpha:\mathfrak{a}\to \sA$, we have the space $\holim_\mathfrak{a}F\circ \alpha$; for a morphism 
\[
(f,\theta):[\alpha:\mathfrak{a}\to \sA]\to [\alpha':\mathfrak{a}'\to \sA]
\]
we have the map
\[
\holim F(f,\theta):=\theta_*\circ f^*:\holim_{\mathfrak{a}'}F\circ \alpha'\to \holim_\mathfrak{a}F\circ \alpha.
\]
This defines a functor  (on the underlying category) 
\[
\Holim F:(\Cat/\sA)_0^\op\to \sSets.
\]

Let $\Cat/\sA^\op$ be the 2-category with underlying category $(\Cat/\sA)_0^\op$, and with 
\[
\Hom_{\Cat/\sA^\op}(\alpha,\alpha') :=\Hom_{\Cat/\sA}(\alpha',\alpha).
\]
The functor  $\Holim F$ extends to a continuous functor
\[
\Holim F:\Cat/\sA^\op\to \sSets.
\]
Indeed,  let $\underline{n}$ be the category associated to the ordered set $\{0<1<\ldots<n\}$. If $G:\sB\to\sSets$ is a functor, we have the functor $G\circ p_2:\underline{n}\times\sB\to \sSets$ and the natural isomorphism of simplicial sets
\[
\sHom(\Delta^n, \holim_\sB G)\cong \holim_{\underline{n}\times\sB} G\circ p_2,
\]
where $\sHom$ denotes the simplicial set of morphisms in $\sSets$. Indeed, it follows directly from the definition of $\holim$ that 
\begin{align*}
 \holim_{\underline{n}\times\sB} G\circ p_2&=\lim_{\substack{\leftarrow\\i\in\underline{n}}}
 \sHom(\sN(\underline{n}/i),\holim_\sB G)\\
&\cong \sHom(\lim_{\substack{\to\\i\in\underline{n}}}\sN(\underline{n}/i),\holim_\sB G)\\
&\cong \sHom(\Delta[n],\holim_\sB G).
\end{align*}

Now suppose we have functors $\alpha:\mathfrak{a}\to\sA$, $\beta:\mathfrak{b}\to\sA$ and an $n$-simplex
\[
\sigma:=[(f_0,\theta_0)\xrightarrow{\vartheta_1}(f_1,\theta_1)\xrightarrow{\vartheta_2}\ldots
\xrightarrow{\vartheta_n}(f_n,\theta_n)].
\]
in $\Hom_{\Cat/\sA}(\alpha,\beta)$.
Define the morphism
\[
(f_\bullet,\theta_\bullet):[\alpha\circ p_2:\underline{n}\times\mathfrak{a}\to\sA]\to
[\beta:\mathfrak{b}\to\sA]
\]
by
\begin{align*}
&(f_\bullet)_{|i\times\mathfrak{a}}:=f_i\\
&f_\bullet((i-1,x)\to (i,x)):=\vartheta_i(x):f_{i-1}(x) \to  f_i(x)\\
&\theta_\bullet(i,x):=\theta_i(x)\\
\end{align*}
This yields the map
\[
\Holim F(f_\bullet,\theta_\bullet):\Holim F(\beta)\to \sHom(\Delta[n], \Holim F(\alpha));
\]
taking the adjoint gives
\[
F(\sigma):\Delta[n]\times\Holim F(\beta)\to \  \Holim F(\alpha).
\]
One checks this respects the simplicial structure in $\sN_*\Hom_{\Cat/\sA}(\alpha,\beta)$ and the composition in  $\sN_*\Cat/\sA$.

For the last ingredient, we recall the notion of a  lax fibering of categories.

 \begin{Def} Let $\pi:\sA\to B$ be a functor of small categories. For each $b\in B$, we have the subcategory $\pi^{-1}(b)$ of $\sA$ consisting of objects over $b$ and morphisms over $\id_b$; let $\alpha_b:\pi^{-1}(b)\to\sA$ be the inclusion. A {\em lax fibering} of $\pi$ consists of the assignment:
\begin{enumerate}
\item for each morphism $h:b\to b'$ in $B$, a functor $h^*:\pi^{-1}(b')\to \pi^{-1}(b)$, and a natural transformation $\theta_h:\alpha_b\circ h^*\to \alpha_{b'}$
\item for each pair of composable morphisms $b\xrightarrow{h}b'\xrightarrow{g}b''$ in $B$, a natural transformation $\vartheta_{g,h}:h^*g^*\to (gh)^*$
\end{enumerate}
which satisfies
\begin{enumerate}
\item[a)] $\pi(\theta_h(x))=h$ for $h:b\to b'$ in $B$, and  $x\in\pi^{-1}(b')$.
\item[b)] for composable morphisms $h:b\to b'$, $g:b'\to b''$, and for $x\in \pi^{-1}(b'')$,
\[
\theta_g(x)\circ\theta_h(g^*x)=\theta_{gh}(x)\circ\vartheta_{g,h}(x).
\]
\item[c)] for composable morphisms $b\xrightarrow{h}b'\xrightarrow{g}b''\xrightarrow{f}b'''$
\[
\vartheta_{h,gf}\circ (h^*\circ\vartheta_{g,f})=   \vartheta_{hg,f}\circ (\vartheta_{h,g}\circ f^*).
\]
\end{enumerate}

If  the natural transformations $\vartheta_{g,h}$ are all identity maps, we call $(\pi:\sA\to B, \theta)$ a {\em strict fibering} of $\pi$.
\end{Def}

We apply these constructions to the following situation: 
Given a lax fibering $(\pi:\sA\to B, h\mapsto h^*,\theta,\vartheta)$, we have the 2-functor $s$ from $\Fac B$ to $\Cat/\sA^\op$ defined by sending $b$ to $\alpha_b:\pi^{-1}(b)\to\sA$, $h:b\to b'$ in $B$ to $(h^*,\theta_h)$, and the 2-morphism $h*g\to hg$ to $\vartheta_{h,g}$. Thus, if we have a functor $F:\sA\to\sSets$, 
we can pull back $\Holim F$ by $s$, giving the continuous functor
\[
\holim F:\Fac B\to \sSets.
\]
Applying the Dwyer-Kan machinery thus gives us the functor
\[
C{\holim F}:B\to \sSets,
\]
and for each $b\in B$ the weak equivalence
\begin{equation}\label{eqn:RectDiag}
i_b:{\holim F}(b)\to   C{\holim F}(b).
\end{equation}
 
\begin{rem}\label{rem:Strict1} Suppose we have a subcategory $i_0:B_0\to B$ of $B$ such that the restriction of $\vartheta$ to pairs of maps in $B_0$ is the identity transformation, i.e., the restriction of $\pi$ to $\pi_{B_0}:\pi^{-1}(B_0)\to B_0$ is a strict fibering. Then $b\mapsto (\holim F)(b)$ extends to a functor $\overline{\holim} F_{|B_0}:B_0\to\sSets$. Applying Theorem~\ref{thm:DwyerKan} to $Y:=
\overline{\holim} F_{|B_0}$,  and noting that $\holim F_{|B_0}=\pi_{B_0}^*Y$, we have the weak equivalence
\[
r_{B_0}:C_{B_0}\holim F_{|B_0}\to \overline{\holim} F_{|B_0}.
\]
The evident inclusion of $n,*$-simplices (see \eqref{multline:HocolimDef}) gives the weak equivalence
\[
\iota_{B_0}:C_{B_0}\holim F_{|B_0}\to(C_B\holim F)_{|B_0},
\]
which finally yields the isomorphism in $\sH \sSets^{B_0}$
\[
\psi:\overline{\holim} F_{|B_0}\to (C_B\holim F)_{|B_0}
\]
by setting $\psi:=\iota_{B_0}\circ r_{B_0}^{-1}$.
\end{rem}

\subsection{Strict functoriality}  We apply the categorical constructions of \S\ref{Sec:Straight} to the problem of functoriality of the homotopy coniveau tower.

We consider the case $B=\Spec k$, $k$ a field, and show how to extend the object $E^{(q)}_\sm$ of $\sH\Spt_\Zar(\Sm\ds B)$ to an object $E^{(q)}$ of $\sH\Spt_\Zar(B)$. The rough idea is to replace $E^{(q)}(X,-)$ with the homotopy limit of the spectra $E^{(q)}(X,-)_{\sC,e}$ over a good index category and use Theorem~\ref{thm:main} to conclude that the resulting homotopy limit is locally weakly equivalent to $E^{(q)}(X,-)$. The index category of the $(\sC,e)$ is chosen so that there are well-defined pull-back maps between the respective homotopy limits for each morphism $f:Y\to X$, however, due to a lack of strict functoriality in the indexing categories, one does not have a strict functoriality on the resulting pull-back maps.  We then use the results of \S\ref{Sec:Straight} to complete the construction.

The main result of this section is

\begin{thm}\label{thm:Funct3} Let $k$ be a field, and take $E\in\Spt(k)$ satisfying axioms  \eqref{eqn:Ax1},
\eqref{eqn:Ax2} and \eqref{eqn:Ax3}. Then the image of $E^{(q)}_\sm:\Sm\ds k\to\Spt$ in $\sH\Spt_\Zar(\Sm\ds k)$ extends (up to isomorphism) to an object $E^{(q)}_\Zar\in \sH\Spt_\Zar(k)$.
\end{thm}

We begin the construction by introducing the category $\sL(\Sm/k)$ (see \cite{MixMot} for slightly modified version). Objects are morphisms $f:X'\to X$ in $\Sm/k$  such that $X'$ can be written as a disjoint union $X'=X'_0\amalg X'_1$ with $f:X'_0\to X$ an isomorphism. The choice of this decomposition is not part of the data. We identify $f:X'\to X$ and $f\amalg p:X'\amalg X''\to X$ if $p:X''\to X$ is smooth, and we identify $f_1:X_1\to X$ and $f_2:X_2\to X$ if there is an $X$-isomorphism $X_1\cong X_2$. We note that a choice of component $X'_0$ of $X'$ isomorphic to $X$ via $f$ determines a smooth section $s:X\to X'$.

A morphism $g:(f_X:X'\to X)\to (f_Y:Y'\to Y)$ is a morphism $g:X\to Y$ in $\Sm/k$ such that there exists a smooth morphism $g':X'\to Y'$ with $f_Y\circ g'=g\circ f_X$. Using the smooth section $s_X:X\to X'$, one easily sees that this condition also respects the identifications 
$f_X\sim f_X\amalg p$ for $p:X''\to X$ smooth: we extend $g':X'\to Y'$ to $g'':X'\amalg X''\to Y'$ by the map $g's_Xp$ on $X''$. Similarly, this condition is unaltered by adding a smooth morphism to $Y'\to Y$ and also respects the identification of isomorphisms over $X$ or over $Y$.  The choice of $g'$ is not part of the data, so 
\[
\Hom_{\sL(\Sm/k)}((f_X:X'\to X),(f_Y:Y'\to Y))\subset \Hom_{\Sm/k}(X,Y). 
\]
Composition is induced by the composition in $\Sm/k$; this is well-defined since smooth morphisms are closed under composition. 

Sending $f_X:X'\to X$ to $X$ defines the faithful functor 
\[
\pi:\sL(\Sm/k)^\op\to \Sm/k^\op. 
\]
We make $\sL(\Sm/k)^\op$ a lax fibered category over $\Sm/k^\op$ as follows: For each morphism $f:X\to Y$ in $\Sm/k^\op$(i.e. $f:Y\to X$ in $\Sm/k$), let $f^*:\pi^{-1}(Y)\to \pi^{-1}(X)$ be the functor
\[
f^*(g:Y'\to Y):= (fg\amalg\id_X:Y'\amalg X\to X).
\]
The map $\theta_f(g):f^*(g:Y'\to Y)\to (g:Y'\to Y)$ in $\sL(\Sm/k)^\op$ is the opposite of the map
$f:(g:Y'\to Y)\to (fg\amalg\id_X:Y'\amalg X\to X)$; we see that this latter is a map in $\sL(\Sm/k)$ by using the inclusion $Y'\to Y'\amalg X$ as the required smooth morphism. Given maps $f:X\to Y$ and $f':Y\to Z$ in $\Sm/k^\op$, and $h:Z'\to Z$ in $\sL(\Sm/k)$, set $\vartheta_{f',f}$ to be the map
\[
f^*(f^{\prime*}(h:Z'\to Z))\\\to
(f'f)^*(h:Z'\to Z) 
\]
being the opposite of the map 
\begin{multline*}
(ff'g\amalg  \id_X:Z'\amalg X\to X)\\\to
(ff'g\amalg f\amalg \id_X:Z'\amalg Y\amalg X\to X)
\end{multline*}
 induced by the identity on $X$. One easily checks that this data does indeed define a lax fibering.

\begin{rem}\label{rem:Strict2} Let $\rho:\Sm\ds k\to \Sm/k$ be the inclusion. It is easy to see that the pull-back  $\rho^{-1}\sL(\Sm/k)^\op\to (\Sm\ds k)^\op$ is a strict fibering.
\end{rem}

Let $f:Y\to X$ be a morphism in $\Sm/k$. We let
\[
\sS^{(q)}_f(p)=\{W\in\sS_X^{(q)}(p)\ | (f\times\id_{\Delta^p})^{-1}(W)\text{ is in }\sS_Y^{(q)}\}.
\]
For $E\in\Spt(k)$, set
\[
E^{(q)}(X,-)_f:=E^{\sS^{(q)}_f}(X\times\Delta^*).
\]

\begin{lem}\label{lem:GoodPos} Let $g:(f_X:X'\to X)\to (f_Y:Y'\to Y)$ be a morphism in $\sL(\Sm/k)$. Then, for $W\in\sS^{(p)}_{f_Y}(n)$, $(g\times\id_{\Delta^n})^{-1}(W)$ is in $\sS^{(p)}_{f_X}(n)$.
\end{lem}

\begin{proof} Choose a component $X'_0$ of $X'$ for which $f_X:X'_0\to X$ is an isomorphism, and let $s:X\to X'$ be the composition
\[
X\xrightarrow{f_X^{-1}}X'_0\subset X'.
\]
Then $s$ is a smooth section to $f_X:X'\to X$, and we can factor $g$ as $g=f_Y\circ g'\circ s$, with $g':X'\to Y'$ smooth. Since $(f_Y\times\id)^{-1}$ maps $\sS^{(p)}_{f_Y}(n)$ to $\sS^{(p)}_{Y'}(n)$, by definition, and $g'\circ s$ is smooth, it follows that $W':=(g\times\id_{\Delta^n})^{-1}(W)$ is in $\sS^{(p)}_X(n)$. Similarly, since $g\circ f_X=f_Y\circ g'$ and 
\begin{multline*}
(f_X\times\id)^{-1}(W')=((f_X\circ g)\times\id)^{-1}(W)\\=((g'\circ f_Y)\times\id)^{-1}(W)=
(g'\times\id)^{-1}((f_Y\times\id)^{-1}(W)),
\end{multline*}
we see that $W'$ is in $\sS^{(p)}_{f_X}(n)$.
\end{proof}

The use of Theorem~\ref{thm:main} for the problem of functoriality is based on the following elementary result:

\begin{lem}\label{lem:FunctSupp} Let $f:Y\to X$ be a morphism in $\Sm/k$. Then there is a finite set $\sC$ of irreducible locally closed subsets of $X$ and a function $e:\sC\to\N$ such that $\sS^{(q)}_f=\sS^{(q)}_{\sC,e}$.  
\end{lem}

\begin{proof} We may assume that $X$ and $Y$ are irreducible. For  a subset $S$ of $Y$, we write
\[
\codim_YS\le d
\]
if there is a point $s\in S$ whose closure $\bar{s}\subset Y$ has $\codim_Y\bar{y}\le d$.

For each $n\in\N$, let
\[
D_n:=\{x\in X\ |\ \codim_Yf^{-1}(x)\le \codim_X\bar{x}-n\}.
\]
By Chevalley's theorem \cite{Jouanolou}, each $D_n$ is a constructible subset of $X$. Write $D_n\setminus D_{n+1}$ as a disjoint union of irreducible locally closed subsets
\[
D_n\setminus D_{n+1}=\amalg_j C^{(n)}_j,
\]
and let $e(C^{(n)}_j)=\codim_X\sC^{(n)}_j-n$. If $x$ is the generic point of $\sC^{(n)}_j$, 
then $e(\sC^{(n)}_j)$ is the codimension of some irreducible component of $f^{-1}(x)$, so 
$e(\sC^{(n)}_j)\ge0$. Letting $\sC:=\{C^{(n)}_j\}$, it is an easy exercise to show that $\sS^{(q)}_f=\sS^{(q)}_{\sC,e}$.  
\end{proof}

We define the functor 
\[
E^{(q)}(?,-)_?:\sL(\Sm/k)^\op\to \Spt
\]
by sending $(f:X'\to X)$ to the functorial fibrant model (in $\Spt$) $\widetilde{E^{(q)}(X,-)_f}$ of $E^{(q)}(X,-)_f$. It follows from Lemma~\ref{lem:GoodPos} that this really does define a functor. Let $hE^{(q)}(X)$ be the homotopy limit
\[
hE^{(q)}(X):=\holim_{(X,f)\in\pi^{-1}(X)}E^{(q)}(X,-)_f.
\]

By the construction of \S\ref{Sec:Straight}, we have the continuous functor
\[
hE^{(q)}:\Fac \Sm/k^\op\to \Spt,
\]
 the functor
\begin{align*}
&ChE^{(q)}: \Sm/k^\op\to \Spt,\\ &X\mapsto ChE^{(q)}(X)
\end{align*}
and for each $X\in\Sm/k$ the weak equivalence 
\[
i_X:hE^{(q)}(X)\to ChE^{(q)}(X).
\]

Thus, for fixed $X\in\Sm/k$, we have the diagram  
\begin{equation}\label{eqn:Diagram1}
\xymatrix{
&hE^{(q)}(X)\ar[r]^{i_X}\ar[d]_{w_X}&ChE^{(q)}(X)\\
E^{(q)}(X,-)\ar[r]_{j_X}&\widetilde{E^{(q)}(X,-)_\id}
}
\end{equation}

Using Remarks~\ref{rem:Strict2} and \ref{rem:Strict1} we have the
diagram of functors on $\Sm\ds k$:
\begin{equation}\label{eqn:Diagram2}
\xymatrix{
&\bar{h}E^{(q)}_\sm  \ar[d]_{w}&C(hE^{(q)}_\sm)\ar[l]_-r\ar[r]^-\iota&ChE^{(q)}\circ\rho\\
E^{(q)}_\sm\ar[r]_{j}&\widetilde{E^{(q)}_\sm}
}
\end{equation}
 $r$, $\iota$  and $j$ are pointwise weak equivalences. 
 
 Note that the fiber category $\pi^{-1}(X)\subset \sL(\Sm/k)^\op$ has the final object $\id_X$.
By Theorem~\ref{thm:main}, Lemma~\ref{lem:FunctSupp} and Remark~\ref{rem:HolimTriv}, $w_X$ is a weak equivalence for each $X$, hence $w$ is a weak equivalence in $\Spt_\Zar(\Sm\ds k)$. Thus, we have isomorphism in $\sH\Spt_\Zar(\Sm\ds k)$
\[
ChE^{(q)}\circ\rho\to E^{(q)}_\sm.
\]

One easily checks that the construction is natural in $q$. Setting $E^{(q)}:=ChE^{(q)}$  completes the proof of Theorem~\ref{thm:Funct3}.

\subsection{Naturality} We conclude our discussion of functoriality with various naturality issues. 

The following result is obvious from the construction of $E^{(q)}(X,-)$.

\begin{lem} The assignment $E\mapsto  E^{(q)}_\sm$ defines a functor $?^{(q)}_\sm:\Spt(\Sm/B^{(q)})\to \Spt(\Sm\ds B)$. $?^{(q)}_\sm$ preserves weak equivalences and homotopy cofiber sequences.
\end{lem}

Similarly, we have
\begin{lem} Fix $X\in\Sm/B$, a finite set of irreducible locally closed subsets $\sC$ of $X$ and a function $e:\sC\to\N$. The assignment $E\mapsto  E^{(q)}(X,-)_{\sC,e}$ defines a functor $?^{(q)}(X,-)_{\sC,e}:\Spt(\Sm/B^{(q)})\to \Spt$, preserving weak equivalences and homotopy cofiber sequences.
\end{lem}
and
\begin{lem} Let $k$ be a field. The assignment 
\[
E\mapsto[(X,f:X'\to X)\mapsto  E^{(q)}(X,-)_f]
\]
 defines a functor $?^{(q)}:\Spt(\Sm/k^{(q)})\to \Spt(\sL(\Sm/k))$, preserving weak equivalences and homotopy cofiber sequences.
\end{lem}

Let $\text{\bf Sec}_\pi(\Sm/B^\op,\sH\Spt_\Nis/B)$ be the category of sections to $\pi:\sH\Spt_\Nis/B\to\Sm/B^\op$. This category inherits a structure of homotopy fiber sequences from $\sH\Spt_\Nis/B$. Let $\Spt(B^{(q)})^*$ be the full subcategory of $\Spt(B^{(q)})$ of objects satisfying the axioms  \eqref{eqn:Ax1},
\eqref{eqn:Ax2} and \eqref{eqn:Ax3}. A weak equivalence (resp. homotopy fiber sequence) in $\Spt(B^{(q)})^*$ is a weak equivalence (resp. homotopy fiber sequence) in $\Spt(B^{(q)})$ involving objects in $\Spt(B^{(q)})^*$. Note that we do not require compatibility with the data in the various axioms, for instance, the contraction in axiom  \eqref{eqn:Ax1}.

From the three lemmas  above follows easily
\begin{thm} \label{thm:Nat1} The assignment $E\mapsto  \tilde{E}^{(q)}$ defines a functor
\[
\tilde{?}^{(q)}:\Spt(B^{(q)})^*\to \text{\bf Sec}_\pi(\Sm/B^\op,\sH\Spt_\Nis/B^*).
\]
$\tilde{?}^{(q)}$ sends weak equivalences to isomorphisms and homotopy fiber sequences to pointwise distinguished triangles.  In case $B=\Spec k$, $k$ a field, the assignment  $E\mapsto  {E}^{(q)}$ defines a functor
\[
{?}^{(q)}:\Spt(B^{(q)})^*\to  \Spt(B)^*.
\]
${?}^{(q)}$ sends weak equivalences to isomorphisms and homotopy fiber sequences to 
weak homotopy fiber sequences.
\end{thm}

\begin{rem} In case of an infinite field, we have similar functor ${?}^{(q)}$ on the category $\Spt(B^{(q)})^*_2$ of presheaves satisfying axioms  \eqref{eqn:Ax1} and
\eqref{eqn:Ax2}.
\end{rem}

We also have naturality with respect to change in $q$. Take $q'>q$, giving us the functor $\rho_{q,q'}:\Sm/B^{(q')}\to\Sm/B^{(q)}$, $(Y,W)\mapsto (Y,W)$. The inclusion of support conditions $\sS_X^{(q')}\subset \sS_X^{(q)}$ gives the natural transformation
\[
\phi_{q,q'}(E):(E\circ\rho_{q,q'})_\sm^{(q')}\to E_\sm^{(q)},
\]
natural in $E$, i.e., the natural transformation $\phi_{q,q'}:{?}^{(q')}_\sm\circ\rho^*_{q,q'}\to {?}^{(q)}_\sm$. We have the identity 
\[
\phi_{q,q'}\circ(\phi_{q',q''}\circ\rho_{q,q'}^*)=\phi_{q,q''}.
\]

\begin{thm} For $q'>q$   there is a natural transformation
$\tilde{\Phi}_{q,q'}:\tilde{?}^{(q')}\circ\rho_{q,q'}^*\to \tilde{?}^{(q)}$ extending $\phi_{q,q'}$. For $q''>q'>q$, we have 
$\tilde{\Phi}_{q,q'}\circ(\tilde{\Phi}_{q',q''}\circ\rho_{q,q'}^*)=\tilde{\Phi}_{q,q''}$. In case $B=\Spec k$, $k$ a field, we have a natural transformation ${\Phi}_{q,q'}:{?}^{(q')}\circ\rho_{q,q'}^*\to {?}^{(q)}$ extending $\phi_{q,q'}$, with the same compatibility. 
\end{thm}
\begin{proof} For each $X,\sC, e$, we have the natural inclusion of support conidtions $\sS^{(q')}_{X,\sC,e}\subset 
\sS^{(q)}_{X,\sC,e}$, giving the natural map
\[
\rho^{X,\sC,e}_{q,q'}:E^{(q')}(X,-)_{\sC,e}\to E^{(q)}(X,-)_{\sC,e}
\]
with $\rho^{X,\sC,e}_{q,q'}\circ \rho^{X,\sC,e}_{q',q''}=\rho^{X,\sC,e}_{q,q''}$.
These maps feed through the construction of $\tilde{E}^{(q)}$ and $E^{(q)}$, respectively, to give the result.
\end{proof}

\begin{cor} The homotopy coniveau tower \eqref{eqn:HCTowerSm} extends to the tower
\[
\ldots\to\tilde{?}^{(q+r+1)}\to\tilde{?}^{(q+r)}\to\ldots\to\tilde{?}^{(q)}
\]
of functors $\Spt(B^{(q)})^*\to \text{\bf Sec}_\pi(\Sm/B^\op,\sH\Spt_\Nis/B)$
 In case $B=\Spec k$, $k$ a field, there is a tower
 \[
\ldots\to{?}^{(q+r+1)}\to {?}^{(q+r)}\to\ldots\to {?}^{(q)}
\]
of functors $\Spt(B^{(q)})^*\to  \Spt_\Nis(B)$ whose composition with the canonical functor
$\Spt_\Nis(B)\to \sH\Spt_\Nis(B)$ extends the tower  \eqref{eqn:HCTowerSm}, considered as a tower of functors
$\Spt(B^{(q)})^*\to  \sH\Spt_\Nis(\Sm\ds B)$.
\end{cor}

\section{Equi-dimensional supports} Over a base-scheme $B$ of positive Krull dimension, especially in the case of mixed characteristic, it is often useful to consider an {\em equi-dimensional} version $\sS_{X/B}^{(q)}$ of the support conditions $\sS_X^{(q)}$. These have two main advantages:
\begin{enumerate}
\item The support condtions $\sS_{X/B}^{(q)}$ admit external products.
\item The support condtions $\sS_{X/B}^{(q)}$ admit an extended functoriality.
\end{enumerate}
Fortunately, one needs to make only minor modifications to the arguments used for the support conditions  $\sS_X^{(q)}$ to cover the equi-dimensional case. In this section, we describe the equi-dimensional theory, and explain these minor changes. We will not go into detail explaining the application to products in this paper, but we will discuss the extended functoriality.

Of course, it would be best if the two theories agreed. Unfortunately, this is at present not known in the most important case, that of mixed characteristic. For this reason, we need to present both theories.

We fix a noetherian base-scheme $B$.

\subsection{Basic definitions} 
\begin{Def} Let $X$ be in $\Sm/B$. \\
(1) Let $\sS_{X/B}^{(q)}(p)\subset \sS_X(p)$ be the set of closed subsets $W\subset X\times\Delta^p$ such that
\[
\text{for all }b\in B, W_b\subset X_b\times\Delta^p\text{ is in }\sS_{X_b}^{(q)}(p).
\]
(2) Let $\sC$ be a finite set of irreducible closed subsets of $X$, $e:\sC\to\N$ a function. Let $\sS_{X/B,\sC,e}^{(q)}(p)\subset \sS_{X/B}^{(q)}(p)$ be the set of closed subsets $W\subset X\times\Delta^p$ such that
\[
\text{for all }b\in B, W_b\subset X_b\times\Delta^p\text{ is in }\sS_{X_b, i_b^*\sC,i_b^*e}^{(q)}(p).
\]
\end{Def}

It is clear that $p\mapsto \sS_{X/B,\sC,e}^{(q)}(p)$ defines a support condition on $X\times\Delta^*$ and that $\sS_{X/B,\sC,e}^{(q)}(p)=\sS_{X/B,\0}^{(q)}(p)=\sS_{X/B}^{(q)}(p)$ if $e(C)\ge \dim_BX$ for all $C\in\sC$.

\begin{Def} Let $E$ be a presheaf of spectra on $\Sm/B^{(q)}$, $X\in\Sm/B$. We set
\begin{align*}
E^{(q)}(X/B,-)&:=E^{\sS_{X/B}^{(q)}}(X\times\Delta^*)\\
E^{(q)}(X/B,-)_{\sC,e}&:=E^{\sS_{X/B,\sC,e}^{(q)}}(X\times\Delta^*)
\end{align*}
\end{Def}

We have as well the equi-dimensional versions of the definitions presented in \S\ref{sec:ExtraGen}:
\begin{Def} Take $X\in\Sch/B$ and $q\in \Z$. \\
(1) The category $X\ds B_{(q)}$ is defined to have objects $(Y,W)$, with $Y\in\Sm/B$ and $W\subset X\times_BY$ a closed subset such that
\[
(Y_b,W_b)\text{ is in } X_b/b_{(q)}
\]
for all $b\in B$. A morphism $(Y,W)\to (Y',W')$ is a $B$-morphism $f:Y\to Y'$ with $W\supset f^{-1}(W')$. \\
(2) Let $\sS^{X/B}_{(q)}(p)\subset \sS_X(p)$ be the set of closed subsets $W\subset X\times\Delta^p$ such that
\[
\text{for all }b\in B, W_b\subset X_b\times\Delta^p\text{ is in }\sS^{X_b}_{(q)}(p).
\]
(3)  Let $\sC$ be a finite set of locally closed subsets of $X$, $e:\sC\to\N$ a function. Suppose that, for each $b\in B$, $C\in\sC$, the locally closed subset $C_b$ of $X_b$ has pure codimension (which may depend on $b$). Let $\sS^{X/B,\sC,e}_{(q)}(p)\subset \sS^{X/B}_{(q)}(p)$ be the set of closed subsets $W\subset X\times\Delta^p$ such that
\[
\text{for all }b\in B, W_b\subset X_b\times\Delta^p\text{ is in }\sS^{X_b, i_b^*\sC,i_b^*e}_{(q)}(p).
\]
\end{Def}

\begin{Def} Take $X\in\Sch/B$, $E\in\Spt(X\ds B^{(q)}$. We set
\begin{align*}
E_{(q)}(X/B,-)&:=E^{\sS^{X/B}_{(q)}}(X\times\Delta^*)\\
E_{(q)}(X/B,-)^{\sC,e}&:=E^{\sS^{X/B,\sC,e}_{(q)}}(X\times\Delta^*)
\end{align*}
\end{Def}

\subsection{The basic principle} The transfer of our main results to the setting of equi-dimensional support conditions is based on the obvious fact that the support conditions $\sS_{X/B,\sC,e}^{(q)}$ and $\sS^{X/B,\sC,e}_{(q)}$  are defined point-wise over $B$, and these support conditions agree with the support conditions $\sS_{X,\sC,e}^{(q)}$ and $\sS^{X,\sC,e}_{(q)}$, respectively, in case $B=\Spec F$, $F$ a field. 

Next,   our two basic moving operations: translation of $\A^n$ by $\A^n$ and taking the cone with respect to a projection, are defined by generic translations or generic projections. In order to make the necessary generic extension $B(x)$ of $B$, we assume that $B=\Spec A$, $A$ a noetherian semi-local ring. For polynomial variables $x_1,\ldots, x_n$, we let $A(x)$ be the localization of $A[x_1,\ldots, x_n]$ with respect to the mulitplicatively closed subset of polynomials $f=\sum_Ia_Ix^I$ such that the ideal of $A$ generated by the $a_I$ is the unit ideal, and we set $B(x):=\Spec A(x)$, as in \S\ref{sec:Notation}.

With this notation, we have as in \S\ref{subsec:Action}  the  generic translation $\psi_x:\A^1_{B(x)}\to \A^n_{B(x)}$, $\psi_x(t)=t\cdot x$, which restricts to the generic translation $\psi_{x,b}:\A^1_{b(x)}\to \A^n_{b(x)}$ for all $b\in B$. Similarly, if $X\subset \A^n_B$ is a closed subscheme, smooth over $B$ and satisfying the conditions \eqref{eqn:Finiteness}, then we have the generic linear projection $\pi_X:X_{B(x)}\to \A^m_{B(x)}$, such that the restriction $\pi_{X_b}:X_{b(x)}\to \A^m_{b(x)}$ remains generic for all $b\in B$.

\subsection{The main results} We conclude this section with a statement of the main results for equi-dimensional support conditions, along with a sketch of their proofs.

\begin{prop}[Moving by translation] Let $B$ be a semi-local noetherian ring with infinite residue fields. Let $\sC_0$ be a finite set of irreducible locally closed subsets of $\A^n_B$, $e:\sC_0 \to\N$ a
function and let $Y$ be in $\Sch/B$. Let $X=Y\times_B\A^n$, $\sC=p_2^*\sC_0$ and $e=p_2^*(e_0)$. Take $E\in\Spt(X\ds B_{(q)})$. Then  the map
\[
E_{(q)}(X/B,-)^{\sC,e}\to E_{(q)}(X/B,-)
\]
is a weak equivalence.
\end{prop}

\begin{thm}[Homotopy invariance] Let $B$ be a semi-local noetherian ring with infinite residue fields, take $X\in\Sch/B$ and $E\in \Spt(X\ds B_{(q)})$. Then the map
\[
p^*:E_{(q)}(X/B,-)\to E_{(q+1)}(X\times\A^1/B,-)
\]
is a weak equivalence.
\end{thm}

To prove these two results, one just repeats the arguments of \S\ref{sec:EasyMov}, using the equi-dimensional support conditions instead of the ones used in \S\ref{sec:ExtraGen}. Since one only needs to verify the various properties of these support conditions point-wise, the case $B=\Spec F$, $F$ a field, applied to the support conditions used in \S\ref{sec:ExtraGen}, suffices to prove all the results we need for the equi-dimensional support conditions. Note that one does not need to allow constructible subsets instead of locally closed subsets as elements of $\sC$, as one simply decomposes a constructible subset into a disjoint union of locally closed subsets and modifies a given function $e$ as in Remark~\ref{rem:LocSuppCond}.

Finally, we have the main result in the equi-dimensional case:

\begin{thm}\label{thm:mainEqui} Let $B$ be a regular  noetherian separated scheme, let $X$ be in $\Sm/B$, let $\sC$ be a finite set of irreducible locally closed subsets of $X$, $e:\sC\to\N$ a function and $q\ge0$ an integer. Take $E$ in $\Spt(B^{(q)})$ satisfying axioms  \eqref{eqn:Ax1},
\eqref{eqn:Ax2} and \eqref{eqn:Ax3}. \\

(1)  Let $x$ be a point of $X$. Then there is a basis of Nisnevic neighborhoods $p:U\to X$ of $x$, depending only on $x$ and $X$,  such that $E^{(q)}(U/C,-)_{p^*\sC,p^*e}\to E^{(q)}(U/B,-)$ is a weak equivalence for each $U$. In particular,  the map of presheaves 
\[
E^{(q)}(X_\Nis/B,-)_{\sC,e}\to E^{(q)}(X_\Nis/B,-)
\]
is a local weak equivalence.\\

(2)  Suppose that $B=\Spec k$, $k$ a field, and let $j:U\to X$ be an affine open subscheme. Then the map of spectra
\[
E^{(q)}(U/B,-)_{j^*\sC,j^*e}\to E^{(q)}(U/B,-)
\]
is a weak equivalence. In particular,  the map of presheaves 
\[
E^{(q)}(X_\Zar/B,-)_{\sC,e}\to E^{(q)}(X_\Zar/B,-)
\]
is a local weak equivalence.\\
\end{thm}

The proof of this result is accomplished by using the argument of \S\ref{sec:Cone}, in the case of $B=\Spec F$, $F$ a field, and our basic principle, to show that, with $f:X\to\A^m_{B(x)}$ the generic finite projection, we have
\[
W\in\sS^{(q)}_{X/B,\sC,e}(p)\Longrightarrow W^+\cup\Ind W\text{ is in } \sS^{(q)}_{X/B,\sC,e-1}(p).
\]
In addition, the map 
\[
X\times\Delta^p\setminus  W^+\cup\Ind W\to \A^m_{B(x)}\times\Delta^p\setminus f( W^+\cup\Ind W)
\]
 induced by $f\times\id$ is \'etale on a neighborhood of $|W|\setminus (W^+\cup\Ind W)$, and the restriction
\[
|W|\setminus( W^+\cup\Ind W)\to f(|W|\setminus (W^+\cup\Ind W))
\]
is an isomorphism. Once we have this, the argument used to prove Theorem~\ref{thm:main} proves the equi-dimensional analog above.

\subsection{Functoriality} Theorem~\ref{thm:mainEqui} yields the same functoriality for the spectra 
 $E^{(q)}(X/B,-)$ (over an arbitrary regular base $B$, or even over a category of such $B$) that we have for the spectra $E^{(q)}(X,-)$ in case $B=\Spec k$, $k$ a field (Theorem~\ref{thm:Nat1}), as well as the naturality of the homotopy coniveau tower. 

Indeed, let $\sV\subset \Sch$ be a small category of regular separated noetherian schemes. Form the category $\Sm/\sV$ with objects being $Y\to B$ in $\Sm/B$ with $B\in\sV$, and morphisms $(f,f_0):(Y\to B)\to (Y'\to B')$ being a pair of morphisms $f_0:B\to B'$ in $\sV$ and $f:Y\to Y'$ a morphisms of schemes over $f_0$.

For an integer $q\ge0$, we have the category $\Sm/\sV^{(q)}$, consisting of pairs $(Y\to B, W)$ with $(Y\to B, W)\in\Sm/B^{(q)}$, where a morphism $(Y\to B, W)\to (Y'\to B', W')$ is a morphism $(f,f_0):(Y\to B)\to(Y'\to B')$ in $\Sm/\sV$, with $f^{-1}(W')\subset W$. Both $\Sm/\sV^{(q)}$ and 
$\Sm/\sV$ are fibered over $\sV$ by the evident forgetful functor, with fiber over $B\in\sV$ being $\Sm/B^{(q)}$ and $\Sm/B$, respectively.

\begin{Def} The category $\Spt(\Sm/\sV)$ of presheaves of spectra on $\Sm/\sV$ is the category of functors $E:\Sm/\sV^\op\to \Spt$ such that the restriction of $E$ to each fiber $\Sm/B^\op$ is a presheaf (i.e. additive). The caetgory $\Spt(\Sm/\sV^{(q)})$ of presheaves of spectra on $\Sm/\sV^{(q)}$ is defined similarly.

We let $\Spt(\Sm/\sV)^*$, resp.   $\Spt(\Sm/\sV^{(q)})^*$ be the full subcategory of $\Spt(\Sm/\sV)$, resp. $\Spt(\Sm/\sV^{(q)})$ consisting of $E$ satisfying axioms  \eqref{eqn:Ax1},
\eqref{eqn:Ax2} and \eqref{eqn:Ax3}.
\end{Def}

 Using equi-dimensional supports throughout, we have the support conditions $\sS^{(q)}_f/f_0\subset \sS_X/B^{(q)}$:  
\[
\sS^{(q)}_f/f_0(p)=\{W\in\sS_X/B^{(q)}(p)\ |\ (f\times\id)^{-1}(W)\text{ is in }
\sS_Y/B^{\prime(q)}\}
\]
We set $E^{(q)}(X/B,-)_f:=E^{\sS^{(q)}_f/f_0}(X\times\Delta^*)$, and we have the well-defined pull-back
\[
f^*:E^{(q)}(X/B,-)_f\to E^{(q)}(Y/B',-).
\]

Using the evident analog of the definition of the collection $\sC$ and the function $e$ from Lemma~\ref{lem:FunctSupp}, this lemma extends directly to the equi-dimensional setting, for an arbitrary morphism $f$ over a morphism $f_0$. The line of argument in Theorem~\ref{thm:Nat1} yields

\begin{thm} \label{thm:Nat1Equi} The assignment $E\mapsto  \tilde{E}^{(q)}$ defines a functor
\[
{?}^{(q)}:\Spt(\Sm/\sV^{(q)})^*\to \Spt(\Sm/\sV)^*.
\]
${?}^{(q)}$ preserves weak equivalences and homotopy fiber sequences.
 
 We have as well naturality with respect to change in $q$, yielding the  homotopy coniveau tower functor ${?}^{(*)}$.
\end{thm} 

\section{The examples} In this section,  we give conditions on a presheaf of spectra $E\in\Spt(B)$ which imply the axioms  \eqref{eqn:Ax1},
\eqref{eqn:Ax2} and \eqref{eqn:Ax3} in \S\ref{sec:axioms}.

Let $E$ be a presheaf of spectra on $\Sm/B$. Recall the extension of $E$ to $\Spt(\Sm/B^{(q)})$ by $E(Y,W):=E^W(Y):=\fib(E(Y)\to E(Y\setminus W))$. Note that this construction is compatible with the restriction functors $\rho_{q,q'}:\Sm/B^{(q')}\to\Sm/B^{(q)}$ for $q'\ge q$. Thus, it suffices to  consider the extension of $E$ to a presheaf on $\Sm/B^{(0)}$, and so we need not refer to the codimension of $W$ in $Y$.

\subsection{Localization} 

For $W=\0$, we have $E^\0(Y)=\cofib(\id:E(Y)\to E(Y))$, hence we have a canonical contraction of $E^\0(Y)$. Let $W'\subset W\subset Y$ be closed subsets of $Y\in\Sm/B$. The Quetzalcoatl lemma gives us the canonical isomorphism
\[
\cofib(i_{W',W}:E^{W'}(Y)\to E^W(Y))\cong E^{W\setminus W'}(Y\setminus W')
\]
which verifies the localization axiom  \eqref{eqn:Ax1}.

\subsection{Nisnevic excision} We need to impose this as an axiom:

\begin{equation}\label{eqn:axiomII}
\vbox{\ \\
Axiom II (Nisnevic excision). Let $f:Y\to X$ be an \'etale morphism in $\Sm/B$ and $W\subset X$ a closed subset. Let $W'=f^{-1}(W)$ and suppose that $f:W'\to W$ is an isomorphism of reduced schemes. Then $f^*:E^W(X)\to E^{W'}(Y)$ is a weak equivalence.
}
\end{equation}

\begin{rems} (1) The condition \eqref{eqn:axiomII} on some $E\in\Spt(B)$ is equivalent to requiring that, for  $\phi:E\to\widetilde{E}$ a fibrant model for the Nisnevic-local model structure, the map $\phi$ is a pointwise weak equivalence. In particular, if $E$ is already fibrant in $\Spt_\Nis(B)$, then 
\eqref{eqn:axiomII} is satisfied.\\
\\
(2) Quillen's localization theorem for $G$-theory together with the weak equivalence $K(X)\to G(X)$ for $X$ regular \cite[Thm. 5, \S 5]{Quillen}  implies that the $K$-theory functor $X\mapsto K(X)$ satisfies \eqref{eqn:axiomII}; indeed, the localization theorem yields the weak equivalence $G(W)\to K^W(X)$, natural with respect to pull-back by smooth morphisms.
\end{rems}

\subsection{Finite fields and $\P^1$-spectra} In case the base-scheme $B$ is $\Spec F$ for a finite field $F$, axiom  \eqref{eqn:Ax3} is valid in case $E$ is the 0-spectrum of a $\P^1$-spectrum; we proceed to explain. 

Let $\SH(F)$ denote the stable homotopy category of $\P^1$-spectra over $F$ (see \cite{Morel, MorelLec, MorelVoev} for details). An object $\sE$ of $\SH(F)$ is a sequence of objects $E_n\in\Spt(F)$, together with maps $\epsilon_n:E_n\to\Omega_TE_{n+1}$. If the $E_n$  are all fibrant objects in the category $\Spt(F)$ (for the pointwise model structure), satisfy the Nisnevic excision axiom II, satisfy homotopy invariance:
\[
p^*:E(X)\to E(X\times\A^1)
\]
is a weak equivalence for all $X$, and the maps $\epsilon_n$ are point-wise weak equivalences, then $E_0$ is the 0-spectrum of $\sE$; in general, the 0-spectrum $\sE_0$ is the colimit of $\Omega^n_T\tilde{E}_n$ where $\tilde{E}_n$ is a fibrant model for $E_n$ in the $\A^1$-local model structure. Here $\Omega_T$ is the operation
\[
\Omega_T(E)(X):=E^{X\times0}(X\times\A^1).
\]

In addition, the 0-spectrum $\sE_0$ of a $\P^1$-spectrum $\sE$ admits weak push-forward maps in the following situation: Take $X\in\Sm/F$ and let $Y\subset X\times_F\A^n$ be a closed subscheme, \'etale over $X$. Suppose that the normal bundle $N_Y$ of $Y$ in  $X\times_F\A^n$ is trivial; fix an isomorphism $\psi:\sO_Y^n\to N_Y$. 

Then there is a push-forward morphism
\[
f^\psi_*:\sE_0(Y)\to \sE_0(X),
\]
natural (in $\SH$) in the triple $(X,Y,\psi)$. 

\subsection{Norms in finite extensions} 
For a field $F$, we have the Grothen\-dieck-Witt ring $\GW(F)$ with augmentation $\dim:\GW(F)\to\Z$ and augmentation ideal $I(F):=\ker\dim$.

Let $F\to L$ be a separable field extension of degree $n$; let $f:\Spec L\to \Spec F$ be the map of schemes. We fix an embedding $i:\Spec L\to \A^1$ and a defining equation $g$ for $i(\Spec L)$. For $X\in \Sm_F$, we thus have the embedding $i_X:X_L\to X\times\A^1_F$ with defining equation $g$. the 1-form $dg$ defines a trivialization of the normal bundle of $X_L$ in $ X\times\A^1_F$, giving the push-forward map $f_*^g:E(X_L)\to E(X)$ for  $E$ the 0-spectrum of a $\P^1$-spectrum $\sE$.

\begin{lem}\label{lem:Invert} Let $F$ be a field such that each element $x\in I(F)$ is nilpotent. Suppose that $E$ is the $0$-spectrum of a $\P^1$-spectrum $\sE\in\SH(F)$. Then, for $X\in\Sm_F$, the map $f_*^gf^*:E(X)[1/n]\to E(X)[1/n]$ is an isomorphism in $\SH$.
\end{lem}

\begin{proof}  For an element $u\in L^\times$, we have the quadratic form $t(u)$ on the $F$-vector space $L$ defined by
\[
(x,y)\mapsto \Tr_{L/F}(uxy).
\]
This gives the class $[t(u)]$ in   $\GW(F)$. In particular, the defining equation $g(x)$ for $ i(\Spec L)$ gives the element $dg/dx\in L^\times$, and the class $[t(dg/dx)]\in \GW(F)$.

$\sE$ is canonically a module for the   sphere spectrum $\1$, so it suffices to prove the result for $\sE=\1$.   By results of Morel \cite{Morel, MorelLec} there is a ring homomorphism $\rho:\GW(F)\to \Hom_{\SH(F)}(\1,\1)$, and $f_*f^*:\1\to \1$ is given by $\rho([t(dg/dx)])$. Since $t(dg/dx)$ has dimension $n$, this element is invertible in $\GW(F)[1/n]$ by our assumption on $F$, hence $f_*f^*$ is invertible in 
$\Hom_{\SH(F)}(\1,\1)[1/n]$.
\end{proof}

\begin{cor}\label{cor:ExtToFiniteFld} Let $F$ be a finite field, $E$ the 0-spectrum of some $\sE\in\SH(F)$. Then $E$ satisfies axiom  \eqref{eqn:Ax3}.
\end{cor}

\begin{proof}
For a finite field $F$, one has $I(F)^2=0$.  Let $F\to L$ be a finite (hence Galois) extension of degree say $n$. By
Lemma~\ref{lem:Invert}, the map $f^*:E(X)\to E(X_L)$ is injective on homotopy groups after inverting $n$. After inverting $n$, the homotopy groups of $E(X_L)^G$ are just the $G$-invariants in $\pi_*(E(X_L))$, which is the same as the image of $f^*f_*$ (after inverting $n$). Thus $E\to \pi_*^G\pi^*E$ is a weak equivalence after inverting $n$.
\end{proof}
\begin{rem} As the 0-spectrum $\sE_0$ of a $\P^1$-spectrum $\sE$ is fibrant in $\Spt_\Nis(B)$, $\sE_0$ satisfies axiom \eqref{eqn:axiomII}. Thus, all our functoriality results are available for $\sE_0$.
\end{rem}

\begin{rem} The codimension $q$ cycles presheaf $z^q\in\Spt(B^{(q)})$   satisfies the axioms  \eqref{eqn:Ax1},
\eqref{eqn:Ax2} and \eqref{eqn:Ax3}. For the localization axiom, let $W'\subset W\subset Y$ be closed subsets of $Y\in\Sm/B$ with $\codim_YW\ge q$. The sequence (of abelian groups)
\[
0\to z^q_{W'}(Y)\xrightarrow{i_*}z^q_W(Y)\xrightarrow{j^*}z^q_{W\setminus W'}(Y\setminus W')\to 0
\]
is well-known to be exact, which verifies the localization property for the associated Eilenberg-Maclane spectra. 

The Nisnevic excision property follows from the isomorphism  
\[
f^*:z^q_W(X)\to z^q_{W'}(Y)
\]
for $f:Y\to X$ \'etale with $W'=f^{-1}(W)$ and $f:W'\to W$ an isomorphism.

The usual push-forward of cycles for proper morphisms gives a well-defined pushforward map $f_*:z^q_W(Y)\to z^q_{f(W)}(X)$ for $f:Y\to X$ a finite morphism in $\Sm/B$. If $f:Y\to X$ has degree $n$, then $f_*f^*$ is multiplication by $n$, from which one easily verifies axiom  \eqref{eqn:Ax3}.

We note that the simplicial spectrum $z^q(X,-)$ is the simplicial Eilenberg-Maclane spectrum associated to the simplicial abelian group $z^q_{ab}(X,-)$, whose associated chain complex $z^q(X,*)$ is just Bloch's cycle complex (defined in \cite{AlgCyc}, at least for $X$ quasi-projective over a field). We have as well the equi-dimensional version $z^q(X/B,*)$, being the chain complex associated to the 
the simplicial abelian group $z^q_{ab}(X/B,-)$, whose simplicial Eilenberg-Maclane spectrum is the simplicial abelian group $z^q(X,-)$.  

Thus,  our results are applicable to Bloch's cycle complexes $z^q(X,*)$, even in the case of $X\in\Sm/B$ with $B$ a regular scheme of Krull dimension one, as well as for the equi-dimensional version
$z^q(X/B,*)$ over an arbitrary neotherian base $B$. Since the presheaf $z^q(X_\Nis,*)$ satisfies localization (see \cite{MovLem}), this yields, in particular, functorial models for $z^q(X,*)$ on $\Sm/k$ for $k$ a field.
\end{rem}

\begin{rem}  As we have already remarked, the $K$-theory presheaf $X\mapsto K(X)$ on $\Sm/B$ satisfies localization, hence the canonical map $K(X)\to R\Gamma K(X)$ is a weak equivalence for all $X\in\Sm/B$. Also, for $p:Y\to X$ a finite \'etale map of degree $n$, we have $p_*:K(Y)\to K(X)$ and $p_*p^*$ induces $\times n$ on the $K$-groups $K_n(X)$. Thus, our full package of functoriality is available for the simplicial spectra $K^{(q)}(X,-)$ and the equi-dimensional version $K^{(q)}(X/B,-)$.
\end{rem}

\section{Good compactifications}\label{sec:GoodCpt}

This section is devoted to proving the following result:

\begin{thm}\label{thm:GoodCpt} Let $B$ be a noetherian  scheme with infinite residue fields, and let $\pi:X\to B$ be in $\Sm/B$, with $X$ equi-dimensional over $B$. For each point $x\in X$, there is a commutative diagram of pointed schemes
\[
\xymatrix{
(Y,x)\ar[r]^f\ar[d]_q&(X,x)\ar[d]^\pi\\
(B',q(x))\ar[r]_g&(B,p(x))
}
\]
and a closed embedding $i:Y\to \A^n_{B'}$ such that
\begin{enumerate}
\item $f$ and $g$ are \'etale
\item  $Y\to B'$ in $\Sm/{B'}$
\item $i$ satisfies \eqref{eqn:Finiteness}
\end{enumerate}
\end{thm}

\subsection{Families of curves} Let $d=\dim_BX$. The first step is to compactify a  neighborhood of $(X,x)$ to form a family of curves over $\P^{d-1}_B$. We then try to  find a good line bundle on the compactification which will contract undesirable components. We begin by reducing to the case in which the fiber $X_{p(x)}$ is geometrically irreducible.

\begin{lem} \label{lem:prep} Let $\pi:X\to B$ be in $\Sm/B$, $x\in X$. Then there is a commutative diagram of pointed schemes
 \[
\xymatrix{
(Y,x)\ar[r]^f\ar[d]_q&(X,x)\ar[d]^\pi\\
(B',b')\ar[r]_g&(B,b)
}
\]
such that
\begin{enumerate}
\item $f$ and $g$ are \'etale
\item  $Y\to B'$ in $\Sm/{B'}$
\item The fiber   $Y_{b'}$ is geometrically irreducible.
\end{enumerate}
\end{lem}

\begin{proof}   In general, since we are allowed to pass to a Zariski open neighborhood of $x$, we may assume that $B$ is local, $B=\Spec\sO$, and that $p(x)$ is the closed point $b$ of $B$. 

Let $L$ be the algebraic closure of $k(b)$ in $k(X_b)$. Then $X_b\to b$ factors through $X_b\to \Spec L$, and $X_b$ is geometrically irreducible over $L$. Now let $B'\to B$ be a local \'etale cover with closed point $b'$ having $k(b')=L$. Then $X\times_BB'\to B'$ has fiber  $X_b\times_{k(b)}L$ over $b'$, hence contains $X_b$ as a connected component and thus contains $x$ as a closed point. Therefore, the projection $X\times_BB'\to X$ defines a Nisnevic neighborhood of $(X,x)$. Removing the complement of $X_b$ in $X_b\times_{k(b)}L$ from $X\times_BB'$ and letting $Y\to B'$ be the resulting open subscheme completes the construction. 
 \end{proof}

We therefore assume that $B=\Spec\sO$, $\sO$ noetherian  and local, with closed point $b=\pi(x)$,  and that the fiber $X_b$ is  geometrically irreducible. We may also assume that $X$ is connected.
 
For a commutative ring $R$ and   an ideal $I\subset R$, we let   $V(I)$ denote the closed subscheme $\Spec R/I$ of $\Spec R$. Similarly, if $S$ is a graded $R$ algebra and $J\subset S$ a homogeneous ideal, we let $V_h(J)$ denote the closed subscheme $\Proj(S/J)$ of $\Proj(S)$.  For a coherent sheaf $\sF$ on a scheme $S$, we have the sheaf of symmetric algebras $\Sym^*_{\sO_S}\sF$ and the projective $S$-scheme $\P_S(\sF):=\Proj_S(\Sym^*_{\sO_S}\sF)$.

\begin{Def} Let $T$ be a $B$-scheme of finite type, $Y\to T$ a $T$-scheme. We call $Y$ a {\em hypersurface} over $T$ if $Y$ is flat over $T$ and there is a  locally free coherent sheaf $\sE$ on $T$, a closed embedding $Y\subset \P_T(\sE)$ and integer $d$ such that each $t\in T$ has an open neighborhood $j:U\to T$ such that $j^*Y\subset \P_U(j^*\sE)$ is the subscheme defined by a non-zero section of $\sO_\sE(d)$ over  $\P_U(j^*\sE)$.
\end{Def}

\begin{rem} A closed subscheme $Y$  of $P:=\P_U(j^*\sE)$ defined by a section $F\neq0$ of $\sO_\sE(d)$ over $P$ is flat over $U$ if and only if for each point $u\in U$, the restriction $i_u^*F$ to $P_u$ is non-zero. Indeed, $\Tor^{\sO_U}_i(\sO_Y,k(u))$ is the $i$th homology sheaf  of the complex
\[
\sO_{P_u}(-d)\xrightarrow{\times i_u^*F}\sO_{P_u}.
\]
\end{rem}

\begin{lem}\label{lem:Family} There is a Zariski open neighborhood $(U,x)$ of $(X,x)$, an open immersion $j:U\to \bar{U}$ and projective morphisms $p:\bar{U}\to\P^{d-1}_B$, $q:\bar{U}\to\P^{d+1}_B$ such that
\begin{enumerate}
\item $p$ has fibers of dimension one.
\item $\bar{U}$ is a hypersurface over $\P^{d-1}_B$.
\item For each $y\in \P^{d-1}_B$, the map $q_y:p^{-1}(y)\to \P^{d+1}_{k(y)}$ is a closed embedding.
\item $p_{|U}:U\to \P^{d-1}_B$ is smooth.
\item  Let $y=p(x)$. Then $U_y$ is geometrically irreducible. 
\end{enumerate}
\end{lem}

\begin{proof} The construction is more or less  standard: We first show that there is an open neighborhood of $x$ in $X$ which is isomorphic to an open subscheme of a hypersurface $V(F)$ in $\A^{d+1}_B$, flat over $B$, i.e. $F$ is in $\sO[X_1,\ldots, X_{d+1}]$ and $F$ is not zero modulo the maximal ideal $\fm$ of $\sO$.

For this, it suffices to replace $x$ with any closed point $x_0$ of $X$ in the closure of $x$. Thus, we may temporarily assume that $x$ is a closed point of $X$. Let $\fm_{X,x}$ denote the maximal ideal in $\sO_{X,x}$.  As the fiber $X_b$ is smooth of dimension $d$ over $k(b)$ , we can find $d$ elements $\bar{f}_1,\ldots,\bar{f}_d\in\fm_{X_b,x}$ which generate $\fm_{X_b,x}$ and form a regular sequence in $\fm_{X_b,x}$. Lifting the $\bar{f}_j$ to elements $f_j\in\fm_{X,x}$, we have the regular sequence $f_1,\ldots, f_d$. As the morphism
\[
\bar{f}:=(\bar{f}_1,\ldots, \bar{f}_d):\Spec\sO_{X_b,x}\to \A^d_{k(b)}
\]
is \'etale, the morphism
\[
f:=(f_1,\ldots, f_d):\Spec\sO_{X,x}\to \A^d_B
\]
is also \'etale. Note that $f(x)$ is the origin $0_b\in\A^d_b$.

By the presentation theorem for \'etale morphisms, there are elements $F, G\in \sO_{\A^d,0_b}[X_{d+1}]$ such that the localization $(\sO_{\A^d,0_b}[X_{d+1}]/F)_G$ is \'etale over $\sO_{\A^d,0_b}$, and a further localization of $(\sO_{\A^d,0_b}[X_{d+1}]/F)_G$ at some closed point is  isomorphic to $\sO_{X,x}$ (as an algebra over $\sO_{\A^d,0_b}$). We may assume that $F$ and $G$ are actually in $\sO[X_1,\ldots, X_{d+1}]$ and that the open subscheme $G\neq0$ of the hypersurface $V(F)$ in $\A^{d+1}_B$ is isomorphic to an open neighborhood of $x$ in $X$. Since $X_b$ has dimension $d$ over $k(b)$, it follows that $F$ is not in $\fm[X_1,\ldots, X_{d+1}]$, so the hypersurface $F=0$ is flat over $B$. Changing notation, we may assume that $X=V(F)\setminus V(G)$.

Let $\bar{F}$ be the homogenization of $F$ to a homogeneous polynomial in $\sO[X_0, X_1,\ldots, X_{d+1}]$  and $\bar{X}\subset \P^{d+1}_B$ the subscheme defined by $\bar{F}$. Clearly 
$\bar{X}$ is a hypersurface over $B$ which contains $X$ as an open subscheme. We extend $\pi:X\to B$ to $\pi:\bar{X}\to B$. Let  $\P^d_\infty\subset\P^{d+1}_B$ be the hyperplane in  defined by $X_0=0$.

Let $P\subset \P^{d+1}_b$ be a general $\P^2$ containing $x$. $P$ will intersect $\bar{X}_b$ transversely at $x$; by Bertini's theorem \cite{Jouanolou}, the intersection $P\cap X_b$ is geometrically irreducible for all sufficiently general $P$. For such a $P$, let  $\ell_b:= P\cap \P^d_\infty$. We may choose $P$ so that $\ell_b\cap \bar{X}_b$ has dimension 0; let $\ell\subset \P^{d}_\infty$ be any $\P^1_B$ with fiber $\ell_b$ over $b$. Then $\ell\cap\bar{X}\to B$ is finite. 
For $z\in \bar{X}\setminus \ell$, let $P_z$ be the $\P^2$ spanned by $z$ and $\ell_{\pi(z)}$.

Let $\mu:\tilde{\P}^{d+1}_B \to\P^{d+1}_B$ be the blow-up of $\P^{d+1}_B$ along $\ell$.   The projection $\P^{d+1}_B\setminus\ell\to \P^{d-1}_B$ with center $\ell$ defines a morphism $\Pi:\tilde{\P}^{d+1}_B\to  \P^{d-1}_B$ making $\tilde{\P}^{d+1}_B$ into a Zariski locally trivial $\P^2$-bundle over   $\P^{d-1}_B$, so $\tilde{\P}^{d+1}_B=\Proj_{\P^{d-1}_B}(\sE)$ for some locally free coherent sheaf $\sE$ of rank three on $\P^{d-1}_B$.

Let $\tilde{X}:=\mu^{-1}(\bar{X})$, and let $p:\tilde{X}\to\P^{d-1}_B$ be the restriction of $\Pi$. Since $\ell\cap\bar{X}\to B$ is finite, $p$ is equi-dimenionsal with fibers of dimension one. Also, $\bar{X}\subset \P^{d+1}_B$ is a hypersurface over $B$, so $\tilde{X}\subset \tilde{\P}^{d+1}_B$ is a Cartier divisor on $\tilde{\P}^{d+1}_B$.  As $p$ has fiber  dimension one, this implies that $\tilde{X}\subset \tilde{\P}^{d+1}_B$ is a hypersurface over $\P^{d-1}_B$.

By construction,  $p^{-1}(p(z))= P_z\cap \bar{X}$ for all $z\in \bar{X}\setminus \ell$. Thus, 
the fiber $p^{-1}(p(x))$ admits an open neighborhood  of $x$ which is smooth over $k(p(x))$ and is geometrically irreducible. Since $p$ is flat, there is a neighborhood $U$ of $x\in \tilde{X}$ which is smooth over $\P^{d-1}_B$, and with fiber $U_{p(x)}$ geometrically irreducible. 

Also, the identity  $p^{-1}(p(z))= P_z\cap \bar{X}$ implies that the  projection $\mu:p^{-1}(y)\to \P^{d+1}_y$ is a closed embedding for all $y\in \P^{d-1}_B$.  We may therefore take   $\bar{U}:=\tilde{X}$ and $q$ the restriction of $\mu$, which completes the proof.
\end{proof}

The next step is to construct a ``numerically contracting" invertible sheaf on $\bar{U}$. For a $B$-morphism $T\to \P^{d-1}_B$, we denote the fiber product $\bar{U}\times_{ \P^{d-1}_B}T$ by $\bar{U}_T$, and similarly use $U_T$ to denote $U\times_{ \P^{d-1}_B}T$.

\begin{lem} \label{lem:GeomIrr} Let $T$ be a noetherian scheme.  Let $p:\bar{U}\to T$ be a hypersurface over $T$ with fibers of dimension one, $0$ a point of $T$,  $x\in\bar{U}_0$. Suppose
\begin{enumera}
\item There is an open neighborhood $U$ of $x$ such that $p:U\to T$ is smooth.
\item There is a Cartier divisor $\sD\subset U$ such that $p:\sD\to T$ is   finite.     
\item The fiber $U_0$ is geometrically irreducible.
\end{enumera}
Then there is a neighborhood $V\subset U$ of $x$ and an open neighborhood $T_0$ of $0\in T$  such that 
\begin{enumerate}
\item $p(V)=T_0$,  
\item $\sD\cap p^{-1}(T_0)\subset  V$,
\item the fibers of $p:V\to T$ are geometrically irreducible. 
\end{enumerate}
Moreover, there is a closed subset $\bar{U}^\0\subset  \bar{U}\setminus V$ so that the fiber $\bar{U}^\0_t$ is a union of irreducible components of $\bar{U}_t$  and 
$\bar{U}_t=\bar{U}^\0_t\cup \overline{V}_t$  for all $t\in T_0$, where $\overline{V}_t$ is the closure of $V_t$ in $\bar{U}$.
\end{lem}

\begin{proof}  We may replace $T$ with any open neighborhood of $0$. Since $\bar{U}$ is a hypersurface over $T$ of relative dimension one,  we may assume that $\bar{U}$ is a closed subscheme of $\P^2_T$, defined by an $F\in H^0(\P^2_T,\sO(d))$ with the property that $i_t^*F\neq0$ for all points $i_t:t\to T$ of $T$. 

Let $d_1,\ldots, d_r$ be be positive integers with $\sum_jd_j=d$, and let $T_{d_1,\ldots, d_r}$ be the subset of $T$ consisting of those points $t\in T$ such that  $i_t^*F$ factors   in $\overline{k(t)}[X_0,X_1,X_2]$ as
\[
i_t^*F=\prod_jF_{t,j},
\]
with $\deg F_{t,j}=d_j$. Then $T_{d_1,\ldots, d_r}$ is a closed subset of $T$, as follows from the fact that the set of $F\in \P(H^0(\P^2,\sO(d)))$ which admit such a factorization is the closed subset given as the image of the multiplication map
\[
\prod_j\P(H^0(\P^2,\sO(d_j)))\to \P(H^0(\P^2,\sO(d))).
\]
Similarly, the set of $t$ for which $i_t^*F$ admits such a factorization with the $F_{t,j}$ irreducible in $\overline{k(t)}[X_0,X_1,X_2]$ is an open subset $T^{\text{irr}}_{d_1,\ldots, d_r}$ of $T_{d_1,\ldots, d_r}$. Clearly the $T^{\text{irr}}_{d_1,\ldots, d_r}$ form a stratification of $T$ by locally closed subsets.

Let $V\subset U$ be an open subscheme. Then the set $\text{Red}(V)$ of $t\in T$ such that the fiber $V_t$ is not geometrically irreducible is a constructible subset of $T$. This follows from the remark above, noticing that $t$ is not in $\text{Red}(V)\cap T^{\text{irr}}_{d_1,\ldots, d_r}$ if and only if $V(F_{t,j})\subset \bar{U}\setminus V$ for all but at most one index $j$, so $\text{Red}(V)\cap T^{\text{irr}}_{d_1,\ldots, d_r}$ is open in $T^{\text{irr}}_{d_1,\ldots, d_r}$.

Now let $W$ be an irreducible closed subset of $T$ containing $0$. Let $w$ be the generic point of $W$. Denote the geometric fiber of $\bar{U}\to T$ over $w$ by $\bar{U}_{\bar{w}}$. We claim that there is a unique irreducible component $\bar{U}_{\bar{w}}^\sD$ which has non-empty intersection with $\sD$. 

Indeed,  note that the hypotheses on $\bar{U}\to T$, $U\subset \bar{U}$ and $\sD$ are preserved under arbitrary base-change $(T',0')\to (T,0)$. Thus, we may assume that  each irreducible component of $\bar{U}_w$ is geometrically irreducible and similarly for $\bar{U}_0$. By assumption, there is a unique irreducible component $\bar{U}_0^\sD$ of $\bar{U}_0$ with non-empty intersection with $\sD$, namely, the closure of $U_0$. Suppose are two components, $\bar{U}_w^1$ and $\bar{U}_w^2$, of $\bar{U}_w$ which intersect $\sD$. Extend the specialization $w\to 0$ to a specialization of cycles
$\bar{U}_w^i\to Z^i\subset \bar{U}_0$, $i=1,2$. Then both $Z^1$ and $Z^2$ intersect $\sD$, so both $Z^1$ and $Z^2$ contain $\bar{U}^\sD_0$. Thus the cycle $\bar{U}_w$ specializes over $w\to0$ to a cycle $Z$ supported in $\bar{U}_0$, containing $\bar{U}^\sD_0$ with multiplicity at least two. This contradicts the smoothness of $U\to T$.

Since the irreducible component  $\bar{U}_{\bar{w}}^\sD$ is unique, this component must be the base-extension to $k(\bar{w})$ of a unique irreducible component  $\bar{U}_w^\sD$ of $\bar{U}_w$. Thus  there is a unique  irreducible component  $\bar{U}_w^\sD$ of $\bar{U}_w$ intersecting $\sD$, and $\bar{U}_w^\sD$ is geometrically irreducible. Let $\bar{U}_w^\0$ be the union of the other irreducible components of $\bar{U}_w$, and let 
$\bar{U}_W^\0$ be the closure of $\bar{U}_w^\0$ in $\bar{U}$.

Let $\bar{U}^\0$ be the union of the $\bar{U}_W^\0$, where $W$ runs over irreducible closed subsets of $T$ containing 0. We claim that $\bar{U}^\0$ is a closed subset of $\bar{U}$. To see this, take such a $W$ and let $w$ be the generic point. There is a unique sequence $d_1\ge\ldots\ge d_r>0$ with $\sum_jd_j=d$ and $w\in T^{\text{irr}}_{d_1,\ldots, d_r}$. There is thus an irreducible closed subset $W'$ of $T$ containing $W$ such that $W'$ is the closure of an irreducible component of $T^{\text{irr}}_{d_1,\ldots, d_r}$; it suffices to show that 
$\bar{U}_W^\0\subset \bar{U}_{W'}^\0$. This is clear, since our condition implies that $\bar{U}_{\bar{w}}$ and $\bar{U}_{\bar{w'}}$ have $r$ irreducible components (counted with mulitplicity), hence $\bar{U}_w^\0$ is in the closure of $\bar{U}_{w'}^\0$.

Let $V_1=U\setminus\bar{U}^\0$. By our construction, the closure $\overline{\text{Red}(V_1)}$ does not contain 0. Thus, taking $V=V_1\setminus p^{-1}(\overline{\text{Red}(V_1)})$, $T_0=T\setminus  \overline{\text{Red}(V_1)}$ proves the result.
\end{proof}

We suppose we have $U$, $p:\bar{U}\to\P^{d-1}_B$ and $q:\bar{U}\to\P^{d-1}_B$ as in Lemma~\ref{lem:Family}. For a $B$-morphism $T\to \P^{d-1}_B$, and a point $t\in T$ lying over $p(x)$, we let $C_t$ denote the closure of  the geometrically irreducible curve $V_t$ in $\bar{U}_t$.

\begin{lem} \label{lem:NumCont} Let $j:U\to\bar{U}$, $p:\bar{U}\to\P^{d-1}_B$, $x\in U$ and $q:\bar{U}\to \P^{d+1}_B$ be as in Lemma~\ref{lem:Family}.   There is a Nisnevic neighborhood $(T,0)\to (\P^{d-1}_B, p(x))$ and an effective Cartier divisor $\sD\subset \bar{U}_T$ such that
\begin{enumerate}
\item $T$ is affine and $\sD\to T$ is finite and \'etale
\item $\sD$ is supported in $U_T\subset\bar{U}_T$.
\item Let $\sL:=\sO_{\bar{U}}(\sD)$. Then  $\sL\otimes\sO_{C_0}\cong \sO_{C_0}(1)$.
\end{enumerate}
\end{lem}

\begin{proof}   Let $H$ be a general  hyperplane in $\P^{d+1}_B$ and let $D:=q^{-1}(H)$.  Choosing $H$ sufficiently general, we may assume that $D\cap C_{p(x)}$ is contained in $U$ and that the intersection is transverse. We may also assume that $D$ intersects the entire fiber $\bar{U}_{p(x)}$ properly.

Thus there is a neighborhood $U_0$ of $p(x)$ in $\P^{d-1}_B$ such that $D\cap\bar{U}_y$ is finite for all $y\in U_0$. Let  $D_0:=D\cap p^{-1}(U_0)$, giving us the finite $U_0$-scheme $p:D_0\to U_0$.

Let $\sO$ be the strict Henselization $\sO_{U_0,p(x)}^{sh}$. Since $D\cap C_{p(x)}$ is contained in $U$ and the intersection is transverse, it follows that $D_0\times_{U_0}\Spec(\sO)$ breaks up as a disjoint union
\[
D_0\times_{U_0}\Spec(\sO)=\amalg_{w\in D\cap C_{p(x)}}\Spec(\sO)\coprod \sD''
\]
where  $\sD''\cap C_{p(x)}=\0$ and  the $\Spec(\sO)$ indexed by $w\in D\cap C_{p(x)}$ is the unique component $Y$ of $D_0\times_{U_0}\Spec(\sO)$ with $Y\cap C_{p(x)}=w$. We label this component as $\Spec(\sO)_w$. 

Let $e=\deg_{\P^{d+1}}(q(C_{p(x)}))$. We thus have the section 
\[
s_1:\Spec(\sO)\to S^{e}(U/\P^{d-1}_B)
\]
with $s_1(y)=\sum_w p^{-1}(y)\cap\Spec(\sO)_w$. Here $S^{e}(U/\P^{d-1}_B)$ is the relative symmetric product:
\[
 S^{e}(U/\P^{d-1}_B):=U\times_{\P^{d-1}_B}\ldots\times_{\P^{d-1}_B}U/\Sigma_{e}.
\]

The zero-cycle $D\cap C_{p(x)}$ is defined over $k(p(x))$ and supported in the smooth locus $U$, so the section $s$ descends to a section on the Henselization $\sO_{U_0,p(x)}^h$:
\[
s_0:\Spec(\sO_{U_0,p(x)}^h)\to  S^{e}(U/\P^{d-1}_B).
\]
Since $U$ is of finite type, there is a Nisnevic neighborhood $(T,0)\to (U_0,p(x))$ such that $s_0$ descends to a section
\[
s: T\to  S^{e}(U/\P^{d-1}_B).
\]
Thus the Cartier divisor $D_0\times_{U_0}T\subset U\times_{U_0}T$ has an open and closed subscheme $\sD$ such that we have the identities of cycles
\[ 
\sD\cap U_t=s(t)
\]
for all $t\in T$. Since $D_0$ is finite over $U_0$, $\sD\to T$ is finite.

Since 
\[
\sD\times_T\Spec(\sO)=\amalg_{w\in D\cap C_{p(x)}}\Spec(\sO)_w,
\]
$\sD\to T$ is \'etale over a neighborhood of 0; shrinking $T$, we may assume that $T$ is affine and $\sD\to T$ is \'etale, completing the proof.
\end{proof}

\begin{lem}\label{lem:GlobalGen} We retain the notation from Lemma~\ref{lem:NumCont}. There is an integer $m_0$ such that $\sL^{\otimes m}$ is generated by global sections for all $m\ge m_0$.
\end{lem}

\begin{proof} The base-scheme $T$ is affine, so write $T=\Spec A$.

 Let $s\in H^0(\bar{U}_T, \sL)$ be the canonical section of $\sL=\sO_{\bar{U}_T}(\sD)$ with divisor $\sD$, and let $m>0$ be an integer. We have the standard sheaf sequence
\[
0\to \sL^{\otimes m-1}\xrightarrow{\times s} \sL^{\otimes m}\xrightarrow{i_\sD^*}
\sL^{\otimes m}\otimes\sO_{\sD}\to 0
\]
which yields the long exact cohomology sequence
\begin{multline}\label{mult:CohSeq}
0\to H^0(\bar{U}_T,\sL^{\otimes m-1})\to H^0(\bar{U}_T,\sL^{\otimes m})
\to\\
H^0(\sD,\sL^{\otimes m}\otimes\sO_{\sD})\to  H^1(\bar{U}_T,\sL^{\otimes m-1})\to \\
H^1(\bar{U}_T,\sL^{\otimes m})
\to H^1(\sD,\sL^{\otimes m}\otimes\sO_{\sD})\to\ldots.
\end{multline}

Since $\sD\to T$ is finite and $T$ is affine, $\sD$ is affine as well. Thus $H^1(\sD,\sL^{\otimes m}\otimes\sO_{\sD})=0$ for all $m\ge0$, and we therefore have the tower of surjections
\[
\ldots\to H^1(\bar{U}_T,\sL^{\otimes m-1})\to H^1(\bar{U}_T,\sL^{\otimes m})\to\ldots
\]
Since $\bar{U}_T\to T$ is projective, the cohomology groups $H^1(\bar{U}_T,\sL^{\otimes m})$ are finitely generated $A$-modules for all $m$, and hence this tower is eventually constant, say
\[
H^1(\bar{U}_T,\sL^{\otimes m-1})\xrightarrow{\sim} H^1(\bar{U}_T,\sL^{\otimes m})
\]
for all $m\ge m_0$. 

Looking back at our sequence \eqref{mult:CohSeq}, we thus  have the surjection
\[
H^0(\bar{U}_T,\sL^{\otimes m})
\to
H^0(\sD,\sL^{\otimes m}\otimes\sO_{\sD})
\]
for $m\ge m_0$. Since $\sD$ is affine, $\sL^{\otimes m}\otimes\sO_{\sD}$ is globally generated for all $m\ge0$, hence $\sL^{\otimes m}_x$ is generated by global sections of $\sL^{\otimes m}$ for all $x\in\sD$ and $m\ge m_0$.  Since $s^{\otimes m}$ generates $\sL^{\otimes m}_x$ for all $x\not\in\sD$, it follows that  $\sL^{\otimes m}$ is generated by global sections for all $m\ge m_0$, as desired.
\end{proof}

\subsection{Finishing the proof}

\begin{prop}\label{prop:MainStep} Let $q:(X,x)\to (B,b)$ be smooth. Then there is a Nisnevic neighborhood $(X_1,x)\to (X,x)$, an open immersion $j:X_1\to \bar{X}_1$, a smooth $B$-scheme $T\to B$ and a projective $B$-morphism $p:(\bar{X}_1,x)\to (T,0)$ such that
\begin{enumerate}
\item $p$ has geometrically irreducible fibers of dimension one.
\item $p\circ j$ is smooth.
\end{enumerate}
Furthermore, $\bar{X}_1$ admits a finite $T$-morphism $\bar{X}_1\to T\times_B\P^1$.
\end{prop}

\begin{proof} By Lemma~\ref{lem:prep}, we may assume that the fiber $X_b$ is geometrically irreducible. Take $U$, $q:\bar{U}\to\P^{d+1}_B$, $(T,0)$, $\pi:(\bar{U}_T,x)\to (T,0)$, $V$ and $\sL$ as given by Lemma~\ref{lem:NumCont}. By Lemma~\ref{lem:GlobalGen}, the sheaf $\sL^{\otimes m}$ is generated by global sections for some $m>0$. As $\bar{U}_T\to T$ has fibers of dimension one, we may find two global sections $s_0$, $s_1$ of $\sL^{\otimes m}$ which generate $\sL^{\otimes m}$, after shrinking $T$. This yields the $T$-morphism
\[
f:=(s_0:s_1):\bar{U}\to T\times_B\P^1.
\]
Since   $\sL^{\otimes m}=\sO_{\bar{U}}(m\sD)$, we may assume that $s_1$ is the
canonical section with divisor $m\sD$.

Let
\[
\xymatrix{
\bar{U}\ar[r]^g\ar[rd]_f&\bar{X_1}\ar[d]^h\\
&T\times_B\P^1
}
\]
be the Stein factorization of $f$ ({\it cf}. \cite[III, Corollary 11.5]{Hartshorne}).

By Lemma~\ref{lem:GeomIrr},  we may shrink $T$ and find an open neighborhood $V$ of $x$ in $U_T$ containing $\sD$ such the fibers of $V\to T$ are geometrically irreducible. For $t\in T$, we let $C_t$ denote the closure in $\bar{U}_t$ of the geometrically irreducible curve $V_t$. We let $C_t^\0$ denote the union of the remaining components of $\bar{U}_t$.  By
Lemma~\ref{lem:GeomIrr}, $\cup_{t\in T} C_t^\0$ forms a closed subset $\bar{U}_T^\0$ of $\bar{U}_T$ with $\sD$ is supported in $\bar{U}_T\setminus \bar{U}_T^\0$.

Take $t\in T$. Since $\sL^{\otimes m}=\sO_{\bar{U}}(m\sD)$ and $\sD$ is supported in $\bar{U}_T\setminus \bar{U}_T^\0$, $f(C_t^\0)$is  a finite subset of $t\times\P^1$, and so $f(\bar{U}_T^\0)$ is a closed subset $W$ of $T\times_B\P^1$, finite over $T$. On the other hand, $\sL$ restricted to $C_t$ has positive degree, so the restriction of $f$ to $f_t:C_t\to t\times \P^1$ is finite.  Furthermore $f(\sD)\subset T\times(1:0)$. Since $C^\0_t\cap\sD=\0$, $s_1$ in nowhere zero on $C^\0_t$. Thus $f(C^\0_t)$ does not contain $t\times(1:0)$, hence  $f(\sD)\cap W=\0$. Shrinking $V$ if necessary, we may assume that $V\subset \bar{U}\setminus f^{-1}(W)$.

It follows  that the restriction of $f$ to
\[
\tilde{f}:\bar{U}\setminus f^{-1}(W)\to T\times_B\P^1\setminus W
\]
is finite. By the universal property of the Stein factorization, the restriction of $g$ to $\bar{U}\setminus f^{-1}(W)$ is an open immersion, and thus the restriction of $g$ to $V$ is an open immersion. We have the morphism $p_1\circ h:\bar{X_1}\to T$.

Since $p:\bar{U}\to T$ is a hypersurface, the fibers of $p$ are all connected. Since $g$ is surjective, it follows that $C_t\to \bar{X}_{1t}$ is surjective for all $t\in T$. Since $C_t$ is geometrically irreducible for all $t$, it follows that the fibers of $p_1\circ h$ are geometrically irreducible. We take $X_1:=g(V)\subset \bar{X}_1$; since $g:V\to \bar{X}_1$ is an open immersion and $V\to T$ is smooth, $X_1\to T$ is smooth. Since $(V,x)\to (X,x)$ is a Nisnevic neighborhood of $(X,x)$, the proof is complete.
\end{proof}

We now proceed by induction.

\begin{thm}\label{thm:irred} Let $\pi:X\to B$ be a smooth $B$-scheme of relative dimension $d$, $x$ a point of $X$. Then there is a Nisnevic neighborhood $(Y,x)\to (X,x)$, a local, \'etale $B$-scheme $(B',0)\to (B,\pi(x))$, a projective $B'$-scheme $q:\bar{Y}\to B'$ and an open immersion $j:Y\to \bar{Y}$ such that
\begin{enumerate}
\item $\bar{Y}$ has geometrically irreducible fibers of dimension $d$.
\item $j(Y)$ is dense in $\bar{Y}$.
\item $q(j(x))=0$.
\end{enumerate}
\end{thm}

\begin{proof} We may assume from the start  that $B$ is local and that $x$ lies over the closed point $b$ of $B$. If $d=0$, then $X\to B$ is \'etale, so we just take $B'=Y=\bar{Y}=X$. In general, assume the result for $X'\to B$ smooth of dimension $<d$. By Proposition~\ref{prop:MainStep}, we have a Nisnevic neighborhood $(X_1,x)\to (X,x)$, a dense open immersion $j_1:X_1\to\bar{X}_1$, a smooth $B$-scheme $V\to B$ and a finite morphism $f:\bar{X}_1\to V\times\P^1$ such that $p_1\circ f$ has geometrically irreducible fibers of dimension one and $p_1\circ f\circ j_1$ is smooth. Let $0=p_1\circ f\circ j_1(x)$.

We apply our induction hypothesis to $(V,0)\to (B,b)$, giving the Nisnevic neighborhood $(V',0)\to (V,0)$, the dense open immersion $j':V'\to\bar{V}'$ and the projective morphism $\pi':\bar{V}'\to B'$ with geometrically irreducible fibers of dimension $d-1$. Taking a base-change with respect to $V'\to V$, we may simply replace $V$ with $V'$ and change notation, so we may assume that $V=V'$. In particular $V\to B'$ has geometrically irreducible fibers of dimension $d-1$ and so $\bar{X}_1\to B'$ has geometrically irreducible fibers of dimension $d$.

Since $f$ is finite, $f$ is projective, so there is a factorization of $f$ as a closed embedding $i:\bar{X}_1\to V\times\P^1\times\P^N$, followed by the projection $V\times\P^1\times\P^N\to
V\times\P^1$. Let $\bar{X}_1^c$ be the closure of $i(\bar{X}_1)$ in $\bar{V}\times\P^1\times\P^N$, giving us the open immersion $j_1^c:\bar{X}_1\to \bar{X}_1^c$ and the projective morphism $p_{12}:\bar{X}_1^c\to \bar{V}\times\P^1$.

Form the Stein factorization  of $p_{12}$:
\[
\xymatrix{
\bar{X}_1^c\ar[r]^q\ar[dr]_{p_{12}}&\bar{Y}\ar[d]^h\\
& \bar{V}\times\P^1
}
\]
Since $\bar{X}_1\to V\times\P^1$ is finite, and $\bar{X}_1= \bar{X}_1^c\times_{\bar{V}}V$, it follows that the composition 
\[
\bar{X}_1\to \bar{X}_1^c\to \bar{Y}
\]
is an open immersion with dense image; we identify $\bar{X}_1$ with its image in $\bar{Y}$. Let $\bar{Y}\to B'$ be the composition
\[
\bar{Y}\to \bar{V}\times\P^1\to \bar{V}\to B'.
\]
We claim that $\bar{Y}\to B'$ has geometrically irreducible fibers.

Indeed, it is easy to see that each irreducible component of $\bar{Y}$ dominates some irreducible component of $B'$, and that the  fiber of $\bar{Y}\to B'$ over a generic point $\eta\in B'$ has
pure dimension $d$ over $k(\eta)$. Take $b'\in B'$. By the upper-semi-continuity of the fiber dimension, each  irreducible component of $\bar{Y}_{b'}$ has dimension $\ge d$ over $k(b')$. Since $\bar{Y}_{b'}\to \bar{V}_{b'}\times\P^1$ is finite, each geometrically irreducible component $P$ of $\bar{Y}_{b'}$
must have
\[
\dim_{\overline{k(b')}}P=d.
\]
Set ${\bar{b}'}:=\Spec\overline{k(b')}$.   

 $\bar{V}_{\bar{b}'}$ is  irreducible of dimension $d-1$. Thus $P\to 
\bar{V}_{\bar{b}'}\times\P^1$ is surjective and the fibers of $P\to \bar{V}_{\bar{b}'}$ are pure dimension one. The fibers of $\bar{Y}\to \bar{V}$ agree with the fibers of $\bar{X}_1\to \bar{V}$ over $V$. Thus, if $\eta$ is the generic point of $\bar{V}_{\bar{b}'}$, we have $P_\eta = (\bar{X}_1)_\eta$, and hence $(\bar{X}_1)_{\bar{b}'}$ is dense in  $P$. Since  $P$ was arbitrary,  $\bar{Y}_{b'}$ is geometrically irreducible. Taking $Y:=X_1$ completes the proof.
\end{proof}

\begin{proof}[Proof of Theorem~\ref{thm:GoodCpt}] Take $(Y,x)\to (X,x)$, $j:Y\to\bar{Y}$ and $\bar{Y}\to B'\to B$ as given by Theorem~\ref{thm:irred}.

$\bar{Y}\to B'$ is projective, so there is a closed embedding $i_1:\bar{Y}\to \P^N_{B'}$. Let $W=\bar{Y}\setminus j(Y)$. Since the fibers of $Y$ are dense in the fibers of $\bar{Y}$, there is a   hypersurface $H\subset \P^N_{B'}$ over $B'$ of sufficiently high degree $m$ containing $W$ and avoiding the point $x$. We may replace $Y$ with $\bar{Y}\setminus H$; taking the $m$-fold Veronese embedding, we may assume that $H=\P^{N-1}_\infty$. The resulting closed embedding
\[
i:Y\to \A^N_{B'}=\P^N_{B'}\setminus\P^{N-1}_\infty
\]
satisfies the condition  \eqref{eqn:Finiteness} since the fibers of $\bar{Y}$ are irreducible.  
\end{proof}

\end{document}